\documentclass[11pt]{article}
\usepackage{pst-plot,pstricks,pst-node}
\usepackage[utf8]{inputenc}
\usepackage{xspace}
\usepackage{epsfig,graphicx}
\usepackage{amsmath,amssymb}
\usepackage{theorem}
\usepackage{setspace}
\usepackage{multirow}
\usepackage{fancyhdr} 
\usepackage[numbers,sort]{natbib} 
\usepackage{empheq}
\usepackage{cleveref}

\usepackage{tikz}
\usepackage{pgfplots}
\usepackage{cancel}
\usepackage{enumerate}
\usepackage{rotating}
\usepackage{pdflscape}
\usepackage{afterpage}
\usepackage{longtable}
\usepackage{url}
\usepackage{subcaption}
\onehalfspace

\usepackage[algo2e,lined,algoruled,linesnumbered,titlenumbered]{algorithm2e}
\allowdisplaybreaks[1]

\SetKwRepeat{Do}{do}{while}%

\setcitestyle{authoryear,open={(},close={)}}

\fancypagestyle{memo}{
	\pagestyle{fancy}
	\fancyhf{}%
	\fancyhead[RO]{\thepage}%
	\fancyhead[LO]{\rightmark}%
	
}%

\theoremstyle{plain}
\newtheorem{prop}{Proposition}
\newtheorem{lemma}{Lemma}

\newtheorem{example}{Example}

\theoremstyle{plain} {
	\theorembodyfont{\rmfamily}%
	
  \newtheorem{remark}{Remark}
}

\newcommand{\ProofNoNL}{{\bf \noindent Proof.}\xspace}

\newcommand{\EndProofNoNL}{\hfill $\Box$ \par \bigskip}

\newcommand{\wloge}{without loss of generality\xspace}

\makeatletter
\renewcommand*\env@matrix[1][c]{\hskip -\arraycolsep
	\let\@ifnextchar\new@ifnextchar
	\array{*\c@MaxMatrixCols #1}}
\makeatother 
\newcommand{\name}{\mbox{(Flow-Cov)}\xspace}
\newcommand{\named}{\mbox{(Path-Cov)}\xspace}
\newcommand{\namet}{\mbox{(Path)}\xspace}

\newcommand{\namedist}[1]{\mbox{$\left(P_{d(i,#1,\delta^{i#1})}\right)$}\xspace}

\newcommand{\var}{y}
\usepackage[normalem]{ulem}

\definecolor{verde}{rgb}{0,0.4,0}
\definecolor{verdeclaro}{rgb}{0.01, 0.75, 0.24}
\begin{document}

\title{
\vspace{-2cm}Upgrading edges in the Maximal Covering Location Problem}

\author{\smaller Marta Baldomero-Naranjo$^1$\footnote{ Corresponding author: marta.baldomero@uca.es (M. Baldomero-Naranjo)}, J\"org Kalcsics$^2$, Alfredo Mar\'in$^3$, Antonio M. Rodr\'iguez-Ch\'ia$^1$ \\[1ex]
\smaller $^1$ Departamento de Estad\'istica e Investigaci\'on Operativa, Universidad de C\'adiz, Facultad de Ciencias,  \\ \smaller 11510, Puerto Real (C\'adiz), Spain, marta.baldomero@uca.es, antonio.rodriguezchia@uca.es\\
\smaller $^2$ School of Mathematics, University of Edinburgh, Edinburgh EH9 3FD, UK, joerg.kalcsics@ed.ac.uk \\
\smaller $^3$  Departamento de Estad\'istica e Investigaci\'on Operativa, Universidad de Murcia, \\ \smaller Facultad de Matemáticas, 30100, Murcia, Spain, amarin@um.es\\
}
\date{}

\maketitle

\vspace{-1cm}
\begin{abstract}
We study the upgrading version of the maximal covering location problem with edge length modifications on networks. This problem aims at locating $p$ facilities on the vertices (of the network) so as to maximise coverage, considering that the length of the edges can be reduced at a cost, subject to a given budget. Hence, we have to decide on: the optimal location of $p$ facilities and the optimal edge length reductions. 

This problem is NP-hard on general graphs. To solve it, we propose three different mixed-integer formulations and a preprocessing phase for fixing variables and removing some of the constraints. Moreover, we strengthen the proposed formulations including valid inequalities. Finally, we compare the three formulations and their corresponding improvements by testing their performance over different datasets. 
\end{abstract}


\section{Introduction}%
\label{sec:Introduction}%
The maximal covering location problem was first introduced by \citet{ChuRev74}. Given a set of clients, each with their own demand, the aim is to locate a fixed number of facilities so as to maximise the amount of covered demand. A client is hereby considered to be covered if their distance to a facility is smaller than or equal to a given coverage radius. Since its origins, this model has been widely studied in the literature under different perspectives. One of the most distinguishing aspects is the solution domain of the problem: continuous \citep{Chu84,Pla02,BanKia17}, discrete \citep{ChuRev74,AveBocVas09,GarMar19,CorFurLju19}, or on networks \citep{ChuMea79,BerKalKrass16,FroMaiHam20}. 
Furthermore, the maximal covering location problem has been solved dealing with alternative coverage assumptions, like gradual coverage \citep{BerKra02} and cooperative coverage \citep{AveBerKraKalNic14,KarEri21}, and with uncertainty, for example uncertainty in the customer demand \citep{BerWan11,BalKalRod20}, in the availability of facilities to provide coverage \citep{Das83,MarMarRodSal18,VatJay21}, or in a combination of these and other parameters\citep{AraBlaFer20,ZhaPenLi17, GuzMas16}.

Common to all those problems is, however, that the parameters of the network and the problem are not decision variables of the model. In this work, we propose a different approach dealing with the maximal covering location problem on networks assuming that edges can be \emph{upgraded} and the total cost of all upgrades is subject to a budget constraint. Upgrading an edge hereby means reducing its length, usually within certain limits, at a given cost which is proportional to the extent of the upgrade. In what follows, we give an example to illustrate the proposed problem.

\begin{example}
Consider the single facility upgrading version of the  maximal covering location problem in the graph depicted in Figure~\ref{fig7:a}. The numbers next to the edges are their lengths, the coverage radius is 11, all the nodes have the same demand, all the edges have a reduction cost of 1 unit per unit, and the maximum reduction is 25\% of the edge length. The problem has been solved for different values of the budget: 0 (Figure~\ref{fig7:a}), 2.5 (Figure~\ref{fig7:b}), 5 (Figure~\ref{fig7:c}), and 10 (Figure~\ref{fig7:d}). The facility is represented as a red diamond, the covered nodes are colored in orange, and the upgraded edges are shown as thicker blue edges.
\begin{figure}[ht] 
  \begin{subfigure}[b]{0.5\linewidth}
    \centering
    \includegraphics[width=0.8\linewidth]{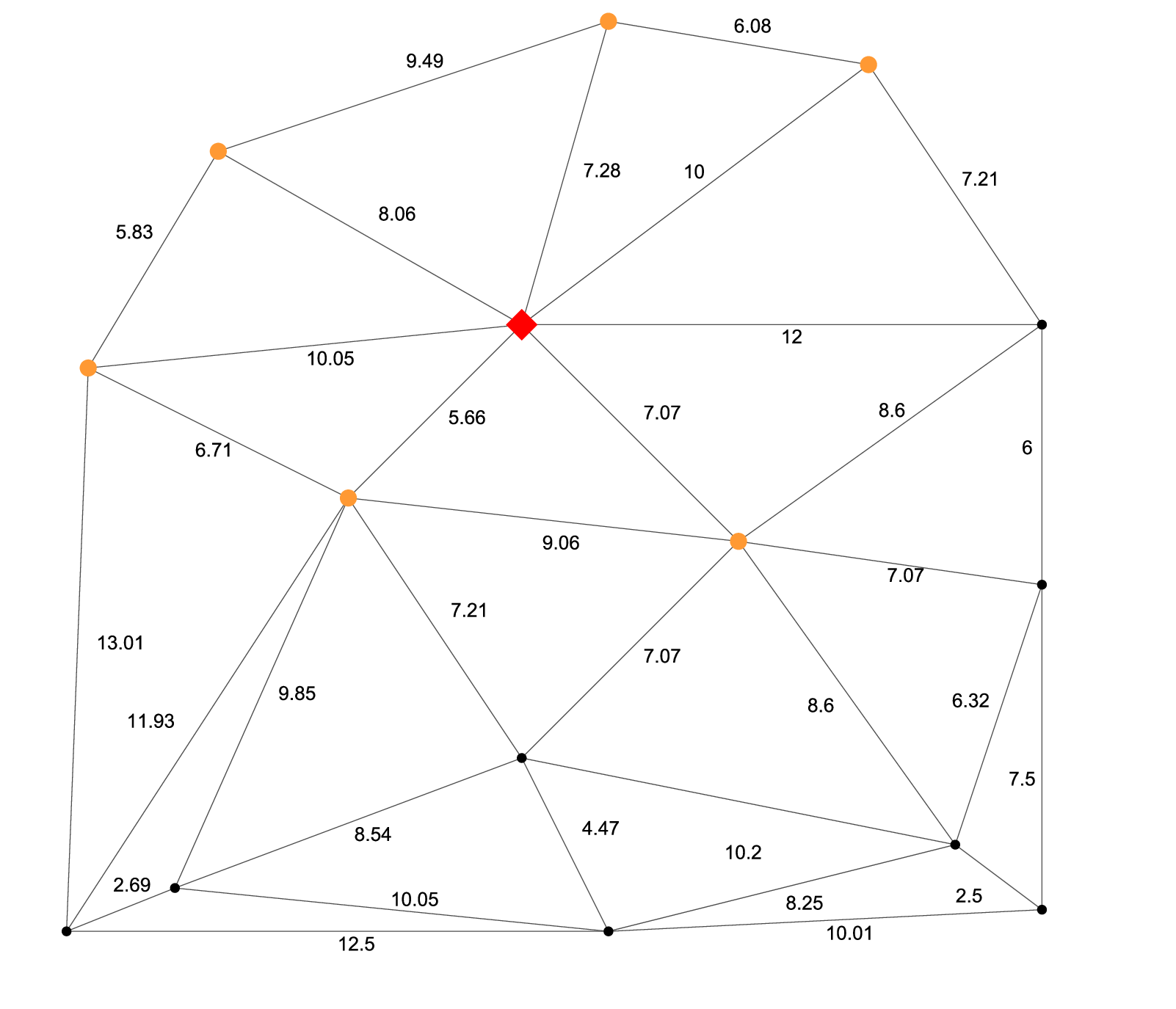} 
    \caption{Without upgrading (budget=0)}
    \label{fig7:a} 
    \vspace{4ex}
  \end{subfigure}
  \begin{subfigure}[b]{0.5\linewidth}
    \centering
    \includegraphics[width=0.8\linewidth]{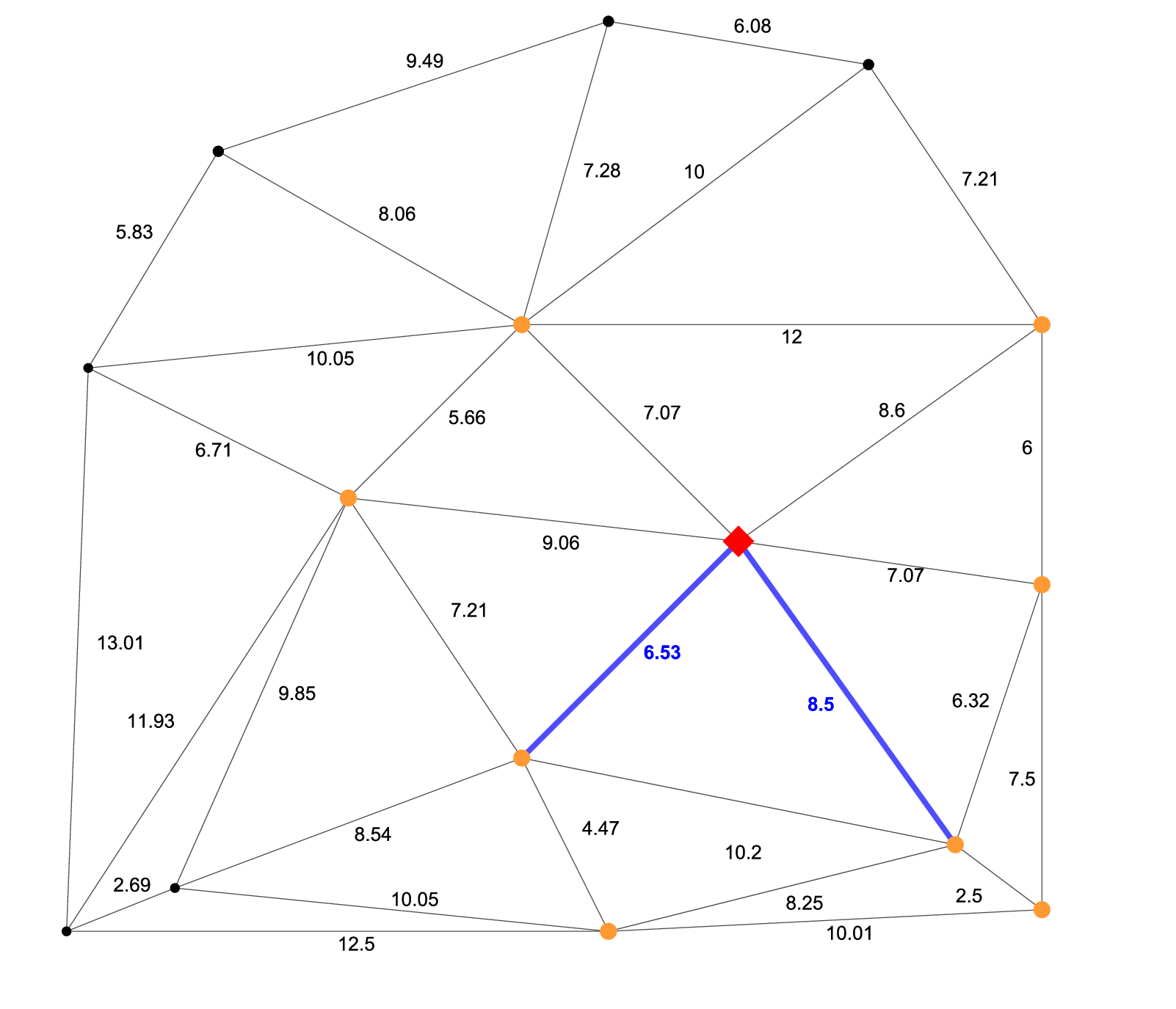}
    \caption{Small budget (budget=2.5)} 
    \label{fig7:b} 
    \vspace{4ex}
  \end{subfigure} 
  \begin{subfigure}[b]{0.5\linewidth}
    \centering
    \includegraphics[width=0.8\linewidth]{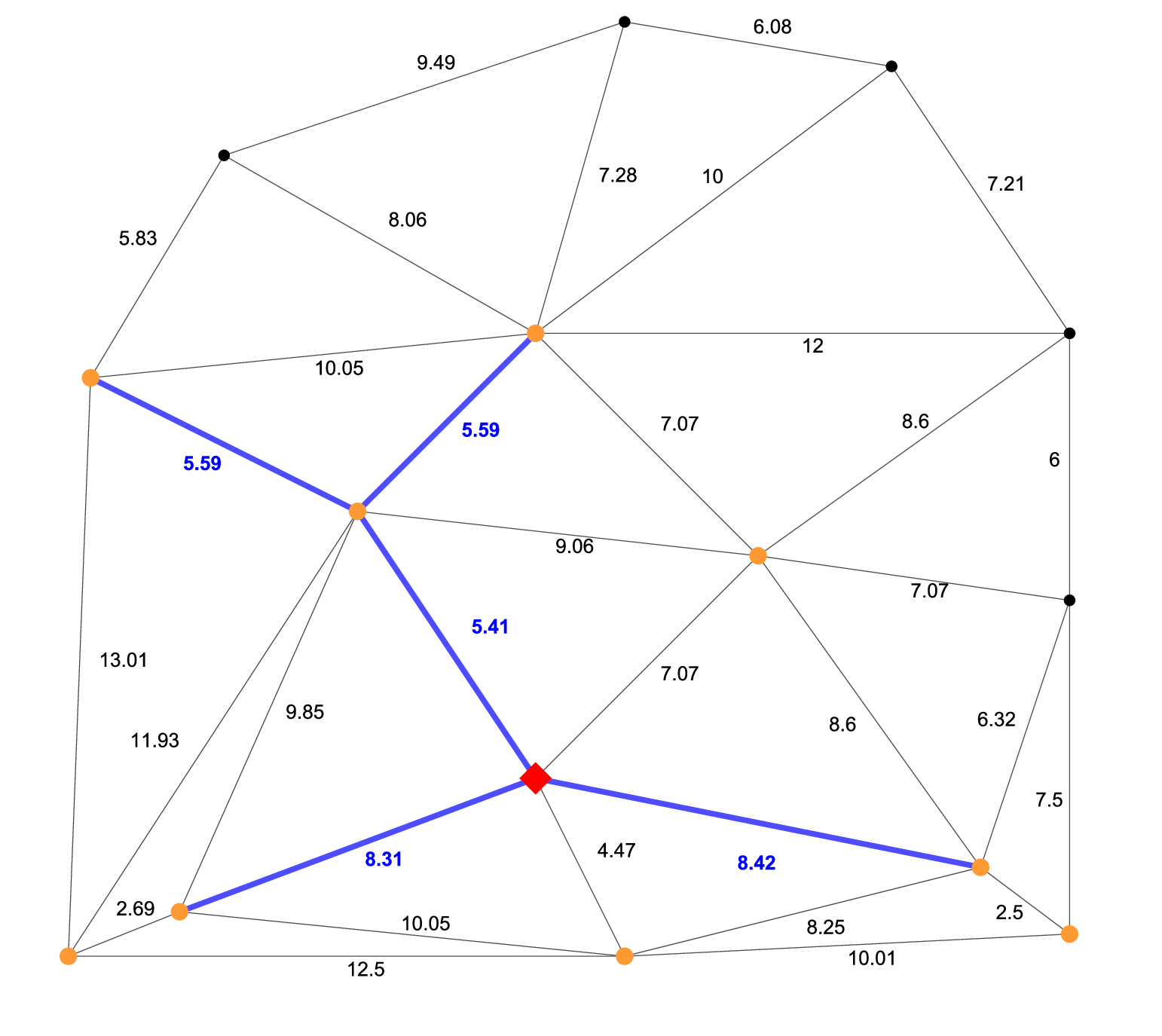} 
    \caption{Medium budget (budget=5)} 
    \label{fig7:c} 
  \end{subfigure}
  \begin{subfigure}[b]{0.5\linewidth}
    \centering
    \includegraphics[width=0.8\linewidth]{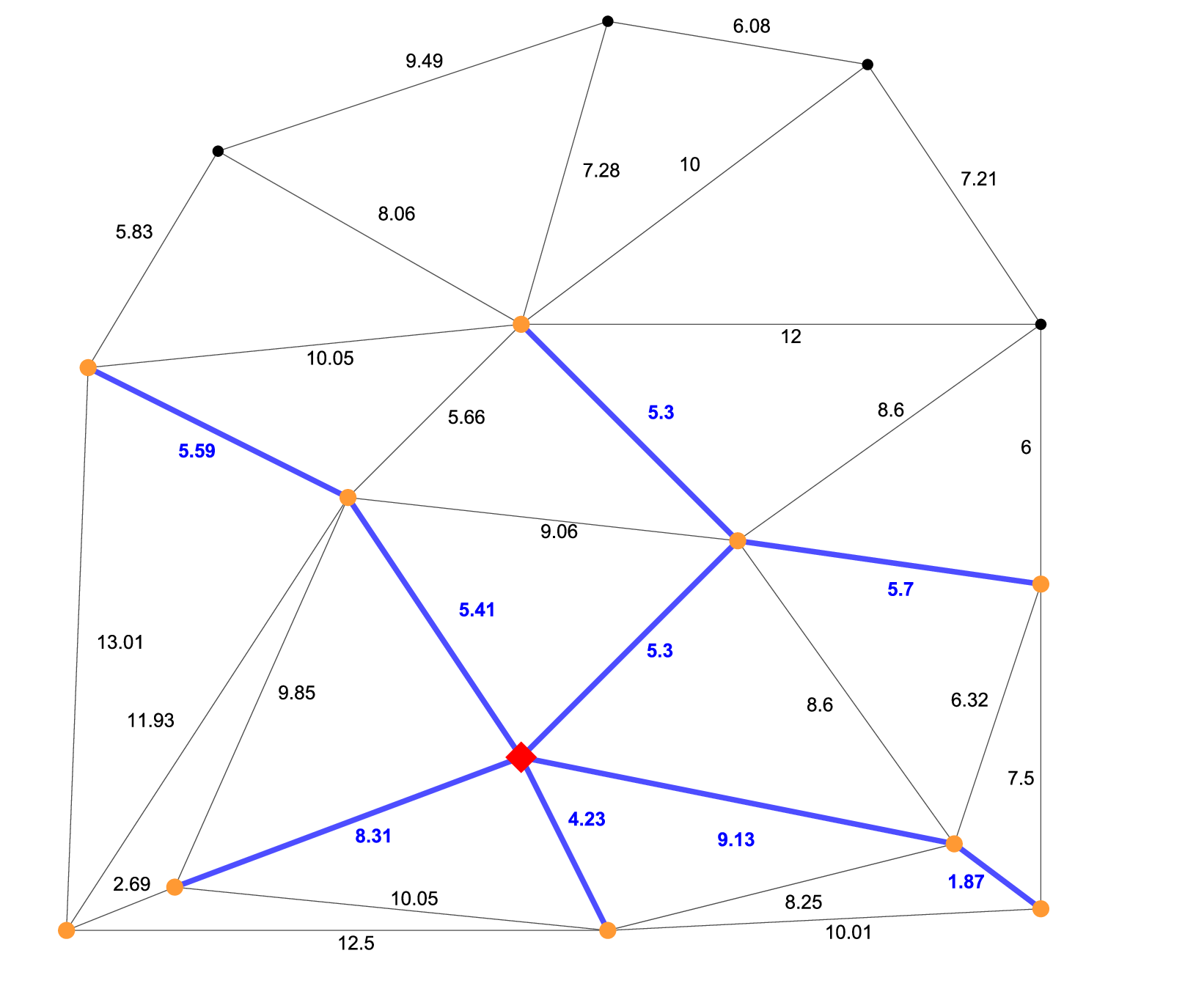} 
    \caption{Large budget (budget=10)} 
    \label{fig7:d} 
  \end{subfigure} 
  \caption{Illustrative example}
  \label{fig7} 
\end{figure}

The pictures show that the location of the facility changes when the budget grows, covering more demand with each increase. In fact, if we fix the optimal location of the problem without upgrading and apply the most beneficial upgrades (i.e., the location of the facility is given and only the upgrade of each edge should be determined), we get that this facility would cover one node less for each of the other three cases. These results illustrate the usefulness of the upgrading version of this problem.
\end{example}

\subsection{Related work}
There are three main types of problems in the literature in which two key parameters of the network, demand weights and edge lengths, can be adjusted, i.e., they are decision variables of the model: 
\begin{itemize}
    \item In \emph{inverse problems}, the objective is to modify one of the two key parameters at minimum cost such that a given feasible solution becomes optimal, see e.g. \citet{HeuSurveyInverse, BurPleZhaInversemedian,BarBurGasInversepmedian,WuLeeZhaWanInverse1median,AliEteReverseObCenter,YanZhaInversecenter,NguSep16,GasInverseOrdMed}.  
    \item In \emph{reverse problems}, the goal is to maximise/minimise the objective value of a given solution by modifying one of the key parameters subject to a given budget. Basically, in reverse problems the roles of variables and input parameters are interchanged, i.e., the variables are considered as parameters and the parameters become variables, see e.g. \citet{BurGasHatReverse1median, BurGasHatReverse2median,WanBaiReverse2median,ZhaYanCaiRervecenter}. 
    \item In \emph{up/downgrading problems},  an actor modifies the parameters of the network and then a reactor takes a decision. In upgrading problems, actor and reactor have the same goal; in  downgrading problems, their objectives are conflicting. 
\end{itemize}

Summing up, the main difference between inverse (reverse) problems and up/downgrading problems is that in the former there is a given solution that we want to improve, while in the latter there is not. 
In this paper, we will focus on the upgrading maximal covering location problem with variable edge lengths. In the following, we denote problems where the edge lengths (demand weights) can be changed as \emph{edge upgrading} (\emph{node upgrading}) problems.

Next, we briefly review the literature of upgrading problems.  
The upgrading version of many classical problems has been studied during the last decades, e.g. 
for the spanning tree problem \citep{ALVAREZMIRANDA201713}, for the min-max spanning tree problem \citep{SeMon17}, for bottleneck problems \citep{BurLinZhaUpDowngradingCombProb}, for minimum flow cost problems \citep{DemNolWirUpMinFlow}, for the shortest path problem \citep{DiLaGo11}, for the maximal shortest path interdiction problem \citep{ZhaGuaPar21}, or for communication and signal flow problems \citep{PaiSah95}. 

In the context of node upgrading location problems (the weight of the vertices can be modified subject to a prespecified budget), the following problems have been analysed: the 1-median problem \citep{Gassner07UpDown1median}, the 1-center problem \citep{Gassner09UpDown1center}, the Euclidean 1-median problem \citep{Plas16}, the $p$-median problem \citep{sepasian15uppmedian}, and the hub-location problem \citep{BlaMar19}, among others. 

In the context of edge upgrading location problems, we are aware of only two directly related publications: upgrading the 1-center problem \citep{Sepasian181center} and upgrading the obnoxious $p$-median problem on trees \citep{AfrAliBar20}. Somewhat related are the models in \citet{Daskin1,Daskin2} which consider the possibility of adding new edges to the network. This can be interpreted as an edge upgrading problem where an edge is in one of two states: non-upgraded with a length of infinity, and upgraded with a finite length. This, however, differs from typical upgrading models where an upgrade can be any fraction of the edge length (or node weight). The models consider the minimization of the overall cost. So far, no results are known for covering problems.
Therefore, the main aim of this paper is to fill this gap in the literature by studying the upgrading maximal covering location problem with variable edge lengths. 

For sake of clarity, we summarise the cited literature of upgrading problems in Table~\ref{tab:Upliter}.  We add two columns labelled N and C for network and continuous problems. 

\begin{table}[htbp]
    \centering
       \resizebox{\hsize}{!}{
    \begin{tabular}{c|l|l|l|l}
    \hline
    Field& Upgrading& N&C& Problem and reference\\
    \hline
         \hline
         \multirow{8}[0]{*}{\shortstack[l]{Location\\ problems}} & \multirow{5}[0]{*}{Nodes} 
  &  X& &1-center: \cite{Gassner09UpDown1center}. \\
        & & X& &  Hub-location: \cite{BlaMar19}.\\
         & & X& & 1-median: \cite{Gassner07UpDown1median}.\\
                 & & & X& Euclidean 1-median: \cite{Plas16}. \\
        & & X && $p$-median: \cite{sepasian15uppmedian}.\\
         \cline{2-5}          
         & \multirow{3}[0]{*}{Arcs/edges} &  X& &1-center: \cite{Sepasian181center}.\\
         & & X&& Maximal covering: our paper. \\
       &    & X& &  Obnoxious $p$-median:  \citet{AfrAliBar20}.\\
            \hline
          \multirow{6}[0]{*}{Others} & \multirow{3}[0]{*}{Nodes}   
            &  X& &Communication and signal flow problems: \cite{PaiSah95}.\\
        &    &X& & Shortest path: \cite{DiLaGo11}.\\
             &  &  X && Spanning tree: \cite{ALVAREZMIRANDA201713}. \\
            \cline{2-5}    
               &  \multirow{3}[0]{*}{Arcs/edges}  & X&& Maximal shortest path interdiction problem: \cite{ZhaGuaPar21}. \\
                & &  X&& Min–max spanning tree: \cite{SeMon17}.\\
       &  &X & & Minimum flow cost: \cite{DemNolWirUpMinFlow}.\\
         \hline
          \hline
    \end{tabular}
    }
    \caption{Summary of literature review of upgrading problems.}
    \label{tab:Upliter}
\end{table}

\subsection{Applications}
This problem has several interesting applications in real-life. Note that two decisions are made at the same time. On the one hand, decide where to locate the $p$ facilities, and on the other hand, determine which edges to upgrade and by how much.    

One application of this problem arises when a public administration wants to improve the accessibility of public services for citizens, e.g. for health centres, educational facilities or social welfare facilities. As the improvement is closely linked with distances \citep{EnsCoo04}, one way to achieve this is to invest in the infrastructure in order to reduce travel times to those services. Such an investment is often a combination of building new facilities and improving the means to get to them, for example by upgrading roads (developing a road into a highway, adding new lanes, etc.) and enhancing public transport (incorporating high-speed lines, adding dedicated bus lines, increasing the frequency of service along links, etc.).

An interesting application in the private sector is for telecommunication companies. To improve their transmission rates and broadband coverage, they will have to increase the bandwidth on existing network links as well as build new or extend existing switching centers. Similar problems are faced by gas and electricity companies who wish to increase their coverage.

Finally, we would like to highlight another useful application in shopping centers, airports, etc. The aim is to locate services such as defibrillators and information posts, in combination with building additional passenger conveyors or escalators to make sure that as many people as possible are within a fixed walking distance of these facilities.

\subsection{Overview}
In this work, we derive three mixed-integer linear programming formulations for the maximal covering location problem with edge upgrades. Furthermore, we develop an effective preprocessing phase that allows us to reduce the dimension of the proposed formulations, allowing us to solve instances faster and also solve larger instances than without preprocessing. Besides, we include several sets of valid inequalities in order to eliminate symmetries and even further improve the solution times of the formulations. 

The rest of the paper is structured as follows. In Section~\ref{sec:ProblemDescription} the problem is introduced. Section~\ref{sec:Flow} presents the first formulation for the problem based on flow variables. Moreover, a preprocessing phase and valid inequalities are developed. Next, in Section~\ref{sec:Fz} and Section~\ref{sec:fy} two new formulations are proposed. In addition, several valid inequalities to enhance them are presented. Section~\ref{sec:ComputationalResults}
contains computational experiments in which we compare the three formulations. We also test the efficiency of the developed valid inequalities.
Finally, our conclusions and some future research topics are included in Section~\ref{sec:Conclusions}.

\section{Definitions and Problem Description}
\label{sec:ProblemDescription}%

Let $N=(V,E,\ell)$ be an undirected network with node set $V=\{1, \ldots, n\}$ and edge set $E$, where $|E|=m$. 
Every edge $e=[k,q]=[q,k]\in E$, $k,q \in V,$ has a positive length $\ell_e=\ell_{[k,q]}$ and is assumed to be rectifiable. For $i,j\in V,$ $d(i,j)$ is the length of the shortest path connecting $i$ with $j$. Furthermore, we are given a fixed coverage radius $R > 0$. We say that a node $i \in V$ is \emph{covered} by a facility at node $j$ if $d(i,j) \leq R.$ Finally, for each node $i \in V$ we are given a non-negative amount $w_i$ that specifies the demand at the node. 

The length $\ell_e$ of each edge $e\in E$ can be reduced by an amount lower than or equal to $u_e$ $\in [0,\ell_e)$, $e\in E$. Without loss of generality, we assume that $\ell_e-u_e\leq R$, for $e\in E$ (if that were not the case, i.e., there were an edge $e \in E$ such that $\ell_{e}-u_{e}> R$, then $e$ can be removed from the network without affecting the optimal solution). Moreover, 
any unit of reduction of the length of the edge $e$ comes at a cost of $c_e$ and there is a budget constraint $B$ on the overall cost of reduction. Again without loss of generality, we assume that $c_eu_e\leq B$, for $e\in E$ (if that were not the case, i.e., there were a cost $c_e$ for $e\in E$ such that $c_{e}u_{e}> B$, then $u_e$ can be substituted by $u_{e}= B/c_{e}$ without affecting the optimal solution).  Finally, we assume that facilities can only be located at nodes. The upgrading maximal covering location problem (Up-MCLP) aims to locate $p$ service facilities covering the maximum demand taking into account that the total cost for the edge length reductions is within the given budget.

Let $\delta = (\delta_e)_{e \in E}$ denote a vector of edge length reductions, $0 \le \delta_e \le u_e$, for $e\in E$. Moreover, let $d(i,j,\delta)$ be the length of a shortest path between nodes $i$ and $j$ after the edge length reductions $\delta$ have been applied, i.e., a shortest path in the network $(V,E,\ell(\delta))$ where $\ell_e(\delta) = \ell_e - \delta_e,$ for $e \in E$. 
Finally, for $p \in \mathbb{N}$ let $X_p \subseteq V$ denote a set of $p$ nodes and let $C(X_p,\delta) = \{i \in V \mid \exists j \in X_p: d(i,j,\delta) \le R \}$ denote the set of all nodes covered by a facility in $X_p$ after the edge upgrades.
Then, Up-MCLP 
can be formulated as:
\[ \max \left\{ \sum_{i \in C(X_p,\delta)}\, w_i\: \Big|\: \sum_{e \in E} c_e \delta_e \le B, X_p \subseteq V, |X_p|=p, 0 \le \delta_e \le u_e, e \in E \right\}.
\]
Table \ref{tab:notation} summarizes the notation used in this paper.

\begin{table}[htb]
    \centering
   \small{
    \begin{tabular}{lp{0.80\textwidth}}
    \hline
    		$A$& Set of all arcs in the induced directed network.  \\
    	$B$& Budget.  \\
    	$c_e$& Unit cost of reducing the length of edge $e,$ $e\in E.$ \\ 
         $d(i,j)$ & Distance between nodes $i$ and $j$ before upgrading, $i,j\in V$. \\
         $d(i,j\delta)$ & Distance between nodes $i$ and $j$ after the edge length reductions $\delta_e$, $e \in E$. In particular, $d(i,j,\delta^{ij})$ represents the distance after the most favourable feasible edge length reductions in the path from $i$ to $j$ and $d(i,j,u)$ the distance in a network with edge lengths $\ell_e-u_e,$ $e\in E.$ \\
         $\Gamma_i$ & Set of edges incident to node $i$ for each $i\in V.$\\
         $\Gamma_i^{+}$ ($\Gamma_i^{-}$) & Set of outgoing (incoming)  arcs for each $i\in V.$ \\
         $m$&Number of edges.\\
         $n$&Number of nodes.\\
          $N=(V,E,\ell)$ & Network with node set $V$, edge set $E$, where $e \in E$ has length $\ell_e$.\\ 
          $p$& Number of facilities. \\ 
         $R$ & Coverage radius.\\
         	$u_e$& Maximum amount that edge $e$ can be reduced, $e\in E$.\\ 
        $\hat{V}_i$ & Set of nodes whose distance to $i$ before upgrading is lower than or equal to $R$, i.e., $\{j\in V\setminus\{i\} : d(i,j) \leq R\}.$\\
         $w_i$ & Demand of node $i,$ for $i\in V.$\\
         \hline
    \end{tabular}
   }
 \caption{Notation used in the paper}
    \label{tab:notation}
\end{table}

Observe that this problem is NP-hard because the maximal covering location problem (MCLP) is a particular case of Up-MCLP (setting $u_e=0$ for all $e\in E$). The NP-hardness of the maximal covering location problem is proved in \citet{Hochbaum97MCLPNPhard}.

\section{Flow coverage formulation} \label{sec:Flow}

In this section, we propose the first of our three  Mixed-Integer Programming (MIP) formulations for Up-MCLP. Using flow variables, the idea of this formulation is to model a path between those pairs of nodes for which the distance between them is smaller than or equal to $R$ after the edge length reductions have been applied. That is, if $d(i,j,\delta) \le R$, then this will be reflected in the formulation by a unit flow between nodes $i$ and $j$. If, however, $d(i,j,\delta) > R$, then the flow between $i$ and $j$ will be zero. We note that in the former case, any path of length $\le R$ will do to assert coverage of $i$ $(j)$ by a service facility located at site $j$ ($i$), so we do not insist on finding the shortest path.

To facilitate the use of flow variables, we consider a directed network $N_{D}= (V,A,\ell)$ with node set $V=\{1,\dots,n\}$ and arc set $A$ containing arcs $(i,j)$ and $(j,i)$ for each edge $[i,j]\in E$. We denote $e_a \in E$ the undirected edge corresponding to $a \in A$ and we define $\Gamma_i^{+}$ ($\Gamma_i^{-}$) as the set of outgoing (incoming)  arcs for each $i\in V.$ 
The set of variables used in the formulation is summarized below. 

{\small
\noindent\textbf{Decision variables}\\
\begin{tabular}{lp{0.85\textwidth}}
	$x_j$&   1, if there is a facility at node $j$, and $0$, otherwise, for $j\in V$.\\ 
	$\var_{ij}$&  1, if node $i$ is assigned to a facility at node $j$, and $0$, otherwise, for $i,j\in V,i\neq j$. \\ 
	$\delta_e$&  The amount of reduction of the length of edge $e$,  for $e\in E$.\\ 
	$f_a^{ij}$& 1, if a path of length $\le R$ from $i$ to $j$ traverses arc $a$, and $0$, otherwise, for $i,j\in V, i<j,$ $a\in A\setminus \left(\Gamma_i^{-} \cup \Gamma_j^{+}\right)$.   \\
	$\alpha_a^{ij}$&  The length of arc $a,$ if this arc belongs to a path of length $\le R$ from node $i$ to node $j$ ($\alpha_a^{ij}=0$ otherwise), for $i,j\in V, i<j,$ $a\in A\setminus \left(\Gamma_i^{-} \cup \Gamma_j^{+}\right)$.
\end{tabular} 
} 

\vspace{0.5cm}
\noindent
Observe that in the definition of the $y$-variables, we use the term ``assign to'' instead of ``covered by''. A node can potentially be covered by more than one service facility and we decided to resolve this ambiguity by explicitly assigning a node to a facility as this simplifies the explanations of the formulations.
Taking into account the notation presented above, the flow coverage formulation for Up-MCLP is:
\begin{alignat}{4} 
		&\hspace{-0.9cm} \name&&\mbox{ max } && \sum_{i\in V}w_i\left(\sum_{j\in V \setminus\{i\}}\var_{ij}+x_i\right) & & \nonumber \\
		 &&&\mbox{ s.t.} &   &  \sum_{j\in V}x_j=p,  & \quad &  \label{ec:pfacility}\\
		   &&& & & \sum_{j\in V\setminus\{i\}}\var_{ij} + x_i \leq 1,   & \quad  & i \in V, \label{ec:zyij}\\
 & & & &&\var_{ij}\leq x_j,  & \quad  & i,j\in V, i\neq j,  \label{ec:zxij}\\
  & & & &&\sum_{e\in E}c_e\delta_e\leq B, \label{ec:B}  & \quad  & \\
& & & &&0\leq \delta_e\leq u_e,   & \quad  & e\in E, \label{ec:ue}\\ 
		&& & & & \sum_{a\in A \setminus \left(\Gamma_i^{-} \cup \Gamma_j^{+}\right)}
		\alpha_a^{ij} \leq R, & & i,j \in V, i<j,  \label{ec:flow1}\\
		&& & & &\alpha_a^{ij} \geq f_{a}^{ij}\ell_{e_a}-\delta_{e_a}, & & i,j \in V, i<j, a\in A  \setminus \left(\Gamma_i^{-} \cup \Gamma_j^{+}\right), \label{ec:flow2}\\
		&& & & &  \sum_{a\in \Gamma_k^{+}, a\notin\Gamma_i^{-}} f_{a}^{ij} - \sum_{a\in \Gamma_k^{-}, a\notin\Gamma_j^{+}} f_{a}^{ij}=0,  & \quad  & i,j\in V, i<j, k\in V\setminus\{i,j\},  \label{ec:flow3}\\
		&& & & &  \sum_{a\in \Gamma_i^{+}} f_{a}^{ij} =\var_{ij}+\var_{ji},  & \quad  & i,j\in V, i<j, \label{ec:flow4} \\
		&& & & &  \sum_{a\in \Gamma_j^{-}} f_{a}^{ij} =\var_{ij}+\var_{ji},  & \quad  & i,j\in V, i<j, \label{ec:flow5} \\
		&& & & &0\leq \alpha_a^{ij}\leq \ell_{e_a}-\delta_{e_a},	& \quad  & i,j\in V, i<j, a\in A\setminus \left(\Gamma_i^{-} \cup \Gamma_j^{+}\right), \label{ec:alphabounds}\\
		& & & &&x_{j}\in\{0,1\},  & \quad  & j\in V,  \label{ec:xy_bin} \\ 
 & & & &&\var_{ij}\in \{0,1\},   & \quad  &  i,j\in V, i\neq j, \label{ec:ybinary} \\
 && & & & f_{a}^{ij}\in\{0,1\}, & & i,j \in V, i<j,  a\in A\setminus \left(\Gamma_i^{-} \cup \Gamma_j^{+}\right). \label{ec:int_f} 
	\end{alignat} 

The objective of the problem is to maximise the amount of covered demand. Constraint \eqref{ec:pfacility} fixes the number of located facilities. The family of constraints \labelcref{ec:zyij} guarantees that either node $i$ is itself a service facility or is assigned to at most one node. The family of constraints \labelcref{ec:zxij} ensures that a node is assigned to an open facility. The families of constraints \eqref{ec:B} and \eqref{ec:ue} force that the reduction on the length of the edges in the network is feasible.
The families of constraints \labelcref{ec:flow1,ec:flow2,ec:flow3,ec:flow4,ec:flow5} ensure that if $\var_{ij}$ $(\var_{ji})$ takes value one, there exists a path shorter than or equal to $R$ from $i$ to $j$ (from $j$ to $i$). Indeed, if $y_{ij}+y_{ji}=1,$ then the flow balance contraints \labelcref{ec:flow3,ec:flow4,ec:flow5} aim at building a path from $i$ to $j$ and consequently from $j$ to $i$. Futhermore, constraints \labelcref{ec:flow1,ec:flow2} ensure that this path is shorter than or equal to $R$. Note that constraints \eqref{ec:zyij} and \eqref{ec:zxij} imply $\var_{ij}+\var_{ji}\leq 1,$ for $i,j\in V, i\neq j.$

\noindent

Next, we introduce a result that proves that the integrality condition of some families of variables of \name  can be relaxed. 

	\begin{lemma}\label{lm:x_relax}
		An equivalent formulation of \name is obtained substituting the set of constraints~\eqref{ec:xy_bin} by: 
		\begin{equation}
		0\leq x_{j}\leq 1,    \quad   j\in V,  \label{ec:xy} 
		\end{equation}
		and the set of constraints~\eqref{ec:ybinary} by: 
		\begin{equation}
		0\leq y_{ij}\leq 1,    \quad   i,j\in V, i\neq j.  \label{ec:y_relax} 
		\end{equation}
	\end{lemma}
	\ProofNoNL
	Concerning the first part of the lemma, the $x$-variables inherit the integrality condition from $\var$-variables due to constraints \eqref{ec:zyij} and \eqref{ec:zxij}. 

	Let $x^*$ and $y^*$ be optimal values for the $x$- and $y$-variables, respectively, of formulation \name when constraints \eqref{ec:xy_bin} are substituted by \eqref{ec:xy}. For any $i\in V$ such that $\sum_{j\in V \setminus\{i\}}y_{ij}^*=1,$ constraint \eqref{ec:zyij} ensures that $x_i^*=0$. On the other hand, for any $i\in V,$ such that there exists $j_0\in V,$ $j_0\neq i,$ with $y_{j_0i}^*=1$, constraints \eqref{ec:zxij} guarantee that $x_i^*=1$. Finally, the model will choose to locate the remaining service facilities (up to a total of $p$) at the uncovered nodes with the largest demand. 
	
	Regarding the second part of the lemma, following a similar argument than before, we conclude that the $y$-variables inherit the integrality condition from the $f$-variables due to constraints \eqref{ec:flow4} and \eqref{ec:flow5} and the condition that $y_{ij}+y_{ji}\leq 1$ (derived by constraints \eqref{ec:zyij} and \eqref{ec:zxij}).
	\EndProofNoNL 
	
Observe that even though the $f$-variables are the intuitive candidates for relaxation (since their number is much larger than the number of $x$- and $y$-variables), there are examples where the optimal value differs when the integrality condition of these variables is relaxed.

An alternative formulation for Up-MCLP can be derived from the formulation \name by replacing constraints \eqref{ec:flow1}, \eqref{ec:flow2}, and \eqref{ec:alphabounds} with the following ones:  
\begin{alignat}{3}
&\sum_{a\in A \setminus \left(\Gamma_i^{-} \cup \Gamma_j^{+}\right)}\left(f_{a}^{ij}\ell_{e_{a}} -\gamma_{a}^{ij}\right)\leq R, & \quad &  i,j\in V, i<j, \label{ec:flow1_gamma}\\
&\gamma_{a}^{ij}\leq u_{e_{a}}f_{a}^{ij},  &&  i,j\in V, i< j, a\in A\setminus \left(\Gamma_i^{-} \cup \Gamma_j^{+}\right), \label{ec:flow2_gamma}\\
& 0\leq \gamma_{a}^{ij} \leq \delta_{e_{a}},  &&  i,j\in V, i< j, a\in A\setminus \left(\Gamma_i^{-} \cup \Gamma_j^{+}\right), \label{ec:gammabounds}
\end{alignat}
where $\gamma_{a}^{ij}$ represents the reduction on the length of arc $a$ if this arc belongs to a path of length $\le R$ from node $i$ to node $j,$  for $i,j\in V, i<j,$ $a\in A\setminus \left(\Gamma_i^{-} \cup \Gamma_j^{+}\right)$.
A preliminary computational analysis showed that the alternative formulation \name for Up-MCLP where constraints \eqref{ec:flow1}, \eqref{ec:flow2}, and \eqref{ec:alphabounds} are replaced with \eqref{ec:flow1_gamma}, \eqref{ec:flow2_gamma}, and  \eqref{ec:gammabounds} is better than the original \name. 
In this analysis, which was carried out using the data described in Subsection~\ref{sub:data}, we compared the formulations based on the number of instances  solved to optimality within the time limit, the time employed to obtain the optimal solutions, the MIP relative gaps, the best solution gaps, and the linear relaxation gaps.
More details about these performance values can be found in Section~\ref{sec:ComputationalResults}.

\subsection{Preprocessing phase}\label{subsec:preplemas}
Next, we present two results for preprocessing the model which reduce the number of constraints and variables of the above formulation and, subsequently, shorten the computational time required to solve them to optimality. The idea of the first is that if the distance between node $i$ and node $j$ is smaller than or equal to $R$ even before modifying the edge lengths of the network, then node $i$ can always be covered by a facility located at node $j$ (and vice versa) regardless of the edge length reductions made in the network.

\begin{prop}\label{lm:cov}
	If $d(i,j) \leq R,$ for $i, j \in V, i< j$ then it is not necessary to include either the $f^{ij}_a$-variables or the $\alpha_{a}^{ij}$-variables ($\gamma_{a}^{ij}$-variables) in formulation \name. Moreover, we can remove the constraints associated with these variables from the family of constraints \labelcref{ec:flow1,ec:flow2,ec:flow3,ec:flow4,ec:flow5,ec:int_f,ec:alphabounds} in the original \name formulation and \labelcref{ec:flow3,ec:flow4,ec:flow5,ec:int_f,ec:flow1_gamma,ec:flow2_gamma,ec:gammabounds} in the alternative \name formulation. 
\end{prop}

The second result analyses the opposite case, i.e., Proposition~\ref{prop:nocov} considers the situation in which the distance between node $i$ and node $j$ is greater than $R$ independently of the edge length reductions.

\begin{prop}  \label{prop:nocov}
If one of the following four conditions is fulfilled for $i,j\in V,$ $i< j$, the variables $y_{ij}, y_{ji}, f_a^{ij}, \alpha_{a}^{ij} (\gamma_{a}^{ij})$ for $a\in A$ can be removed from the \name formulation. Moreover, the constraints associated with this pair of nodes can be deleted, in particular \labelcref{ec:zyij,ec:zxij,ec:flow1,ec:flow2,ec:flow3,ec:flow4,ec:flow5,ec:int_f,ec:ybinary,ec:alphabounds} in the original \name formulation and  \labelcref{ec:zyij,ec:zxij,ec:flow1_gamma,ec:flow2_gamma,ec:flow3,ec:flow4,ec:flow5,ec:int_f,ec:ybinary,ec:gammabounds} in the alternative \name formulation. 
\begin{enumerate}[i)]
    \item $d(i,j)>R+\sum_{e\in E}u_e,$ for $i, j \in V, i<j.$
    \item $d(i,j,u)>R,$ for $i,j\in V, i<j,$ where $d(i,j,u)$ is the length of the shortest path from $i$ to $j$ in a graph with edge lengths $\ell_e-u_e,$ for $e\in E.$
    \item $\displaystyle d(i,j)>R+\sum_{k=1}^{\bar{k}}u_{e_{\sigma(k)}}+\dfrac{B-\displaystyle\sum_{k=1}^{\bar{k}}u_{e_{\sigma(k)}}}{c_{e_{\sigma(\bar{k}+1)}}}$ for $i, j \in V, i< j,$ where $\bar{k}$ is the largest index $k$ that satisfies the following condition:
\begin{equation}
    \sum_{h=1}^{k} u_{e_{\sigma(h)}}{c}_{e_{\sigma(h)}}\leq B, \label{ec:iiB}
\end{equation}
and $\sigma(\cdot)$ is a permutation of  $\{1,\dots,m\}$ that sorts the unit upgrade costs in non-decreasing order. 
      \item The optimal value of the following problem is greater than $R,$ for $i,j \in V, i<j,$
\begin{alignat}{4}
&\namedist{j}&\quad & \min&  \quad & \sum_{a\in A} \left(f_{a}\ell_{e_{a}} -\gamma_{a}\right) & & \nonumber \\
&& &\mbox{s.t.} &  \quad & \eqref{ec:B},\eqref{ec:ue},\nonumber\\
&&& & & \sum_{a\in \Gamma_k^{+}} f_{a} - \sum_{a\in \Gamma_k^{-}} f_{a}=g_k,  & \quad  & k\in V,\\
&&& & & \gamma_{a}\leq u_{e_{a}}f_{a},  & \quad  & a\in A,\\
&&& & & \gamma_{a} \leq \delta_{e_{a}},  & \quad  & a\in A, \\
&&& & & f_{a}\in\{0,1\}, & \quad  & a\in A, 
\end{alignat}
where the $f$-variables and $\gamma$-variables are defined as above (we dropped the indices $i$ and $j$ for the ease of exposition), and
\begin{equation*}
g_k=\begin{cases}
1, & \mbox{if } k=i, \\
-1, & \mbox{if } k=j, \\
0, & \mbox{otherwise}. 
\end{cases}
\end{equation*}
\end{enumerate}
\end{prop}
\ProofNoNL Each of the items of the proposition is proven below. 
\begin{enumerate}[\itshape i)]
    \item The first condition considers the case where the distance from $i$ to $j$ is greater than $R$ even when reducing the length of every edge by the maximum amount allowed. Therefore, it is straightforward to conclude that the distance between the two nodes cannot be less than or equal to $R$.
    \item In the previous condition, the maximum amount of reduction in the whole network was considered without taking into account the edges for which this reduction is made.  
    Now we compute the shortest path between two nodes in the network assuming an unlimited budget, i.e., the full discount is applied to all edges. 
    For each edge $e$, let $\ell_{e_u}=\ell_e-u_e$. For $i,j\in V,$ let $d(i,j,u)$ be the length of the shortest path connecting $i$ with $j$ where the length of the edges are $\ell_{e_u},$ for $e\in E$. Hence, even if reducing the maximum amount allowed on all edges, the distance is greater than $R$, clearly, node $i$ cannot be assigned to node $j,$ and vice versa. Although condition \textit{i)} is weaker than \textit{ii)}, \textit{i)} can be checked more efficiently.
  
    \item This condition is similar to the first one, but takes into account the budget constraint~\eqref{ec:B}. In this case, we calculate the maximum reduction in the network  allowed by the budget. For this purpose, we sort the upgrade costs $c_e$, for $e\in E$, in non-decreasing order. Let $\sigma$ be a permutation of $\{1,\dots,m\}$ such that $c_{e_{\sigma(1)}} \le c_{e_{\sigma(2)}} \le \ldots \le c_{e_{\sigma(m)}}.$
Then, we compute the maximum total length reduction over the network, i.e., we spend the budget on upgrading the cheapest edges. Let $\bar{k}$ be the largest index $k$ that satisfies  condition \eqref{ec:iiB}. 
Therefore, the right-hand side of \textit{iii)} minus $R$ is the maximum length reduction between any two nodes of the network.
Taking into account the above arguments, we conclude that node $i$ cannot be assigned to a facility located at node $j,$ and vice versa. 

\item Condition \textit{iii)} provides the maximal reduction without taking into account whether this reduction can be achieved in a path form $i$ to $j$. For this reason, that bound can be tightened, but it requires to solve a separate problem for each pair of vertices. Formulation \namedist{j} computes the shortest path between node $i$ and node $j$ assuming that all the budget can be spent just for the path between those two nodes. Therefore, the optimal value of this problem, named $d(i,j,\delta^{ij})$, is the minimal distance between node $i$ and node $j$ after the most favourable edges length reductions. Hence, if $d(i,j,\delta^{ij})$ is greater than $R$, node $i$ can never be assigned to a facility at node $j,$ and vice versa.  \hfill $\Box$

\end{enumerate}

As stated in \cite{DemNolWirUpMinFlow}, the shortest path problem where the length of the edges can be reduced, \namedist{j}, is NP-hard. 
However, the optimal value of the LP relaxation of formulation \namedist{j} provides a valid bound that can still be used instead, albeit yielding a weaker condition. If this value is greater than $R$, then $i$ can never cover $j$, and vice versa.

\subsection{Valid inequalities}
In the previous subsection, we have presented two results to preprocess the model, reducing the number of constraints and variables. In this one, we propose several families of valid inequalities to strengthen the \name formulation which help us to further shorten the computational times.

\begin{prop}
Let $\hat{V}_i:=\{j\in V\setminus\{i\} : d(i,j) \leq R\}$. The following families of constraints 
    are valid inequalities for \name:
\begin{alignat}{3}
		& f_{(k,q)}^{ij}+ f_{(q,k)}^{ij}\leq 1, &&[k,q]\in E, i,j\in V, i<j, k,q\neq i,j, \label{ec:flow6}\\
&\var_{kj}\geq \var_{ij}+f_{a}^{ij}-1, &\quad  & i,j,k\in V,  i<j, k\neq i, k\neq j, a\in\Gamma_k^{-}, \label{ec:ddvvconection_1}  \\
&\var_{ki}\geq\var_{ji}+f_{a}^{ij}-1, &\quad &i,j,k\in V, i<j, k\neq i, k\neq j, a\in\Gamma_k^{-}, \label{ec:ddvvconection_2} \\
   &x_j \leq \sum_{k: k\neq i, d(i,k)\leq d(i,j)} y_{ik} +x_i, &\quad &i\in V, j \in \hat{V}_{i}, \label{lm:f2_flow} \\
   &x_j \leq \sum_{k: k\neq i, d(i,k,\delta^{ik})\leq d(i,j)} y_{ik} +x_i, &&i\in V, j \in \hat{V}_{i}, \label{lm:f2_flow_delta}  
\end{alignat}
\end{prop}

\ProofNoNL
The proof of valid inequalities is given below: 

The first family of constraints \eqref{ec:flow6}  ensure that an edge is not traversed in both directions on a path from $i$ to $j$, $i<j$. 

The second and the third families of constraints, \eqref{ec:ddvvconection_1} and \eqref{ec:ddvvconection_2}, are based on the fact that a path between two non-adjacent nodes $i$ and $j$ will traverse at least one other node. Therefore, if there exists a path whose length is less than or equal to $R$ that connects a facility at node $j$ ($i$) with demand point $i$ ($j$) and traverses node $k$, then node $k$ will also be assigned to facility $j$ ($i$). More concretely, given $i,j \in V, i<j,$ if $f_{a}^{ij}=1$ for some $a \in A,$ such that $a\in\Gamma_k^{-},$ $k\neq i, k\neq j$ and $\var_{ij}=1$ $( \var_{ji}=1),$ then the constraints impose that $\var_{kj}=1$ $( \var_{ki}=1)$.

Regarding \eqref{lm:f2_flow}, these constraints ensure that a node will be served by the closest service facility that is within the covering distance before upgrading the network (whenever at least one service facility is closer than the coverage radius before upgrading the network). Observe that these constraints eliminate symmetries and are  valid also for formulation \name  because nodes might not be assigned to the closest service facility in the upgraded network. Nevertheless, the situation would be incompatible with constraints \eqref{ec:ddvvconection_1} and \eqref{ec:ddvvconection_2}, as explained in the following remark (Remark~\ref{re:1}). 

Finally, whenever at least one service facility $j$ is closer to a node $i$ than the coverage radius before upgrading the network, constraints \eqref{lm:f2_flow_delta} ensure that this node will either host a facility itself or be assigned to this service facility or to a facility that can be closer after upgrading the network (it considers the distances in the range of
$d(i,j)$ and the most favourable edge length reductions, i.e., $d(i,k, \delta^{ik})$ for any $k \ne i$.).
\EndProofNoNL
\begin{remark}\label{re:1}
\begin{enumerate}
    \item[\textit{i)}] Constraints \eqref{ec:ddvvconection_1} and \eqref{ec:ddvvconection_2} might be incompatible with  \eqref{lm:f2_flow}, i.e., constraints \labelcref{ec:ddvvconection_1,ec:ddvvconection_2,lm:f2_flow} cannot be included in the formulation simultaneously. 
    \item[\textit{ii)}] The family of valid inequalities \eqref{lm:f2_flow} is tighter than \eqref{lm:f2_flow_delta}, but \eqref{lm:f2_flow_delta} are not incompatible with \eqref{ec:ddvvconection_1} and \eqref{ec:ddvvconection_2}. 
\end{enumerate} 
\end{remark}
In the following, we present an example illustrating the first part of Remark~\ref{re:1}. \begin{example}\label{ex:1}
Consider the network depicted in Figure~\ref{fig:Example}. For each edge, its length, its upper bound of reduction, and its cost per unit of reduction, $(\ell_e,u_e,c_e),$ are printed next to the edge. Let $R=1, p=2, B=0.75,$ and the demand of nodes $w_i=1$, $w_j=1$, $w_k=1$, $w_q=1$, $w_r=1000$, $w_s=1000.$  It is straightforward to conclude that the optimal location of the services are the dark nodes, i.e., $x_i^*=1$ and $x_q^*=1$, and that the optimal edge length reduction is  $\delta_{[k,q]}^*=0.75.$  
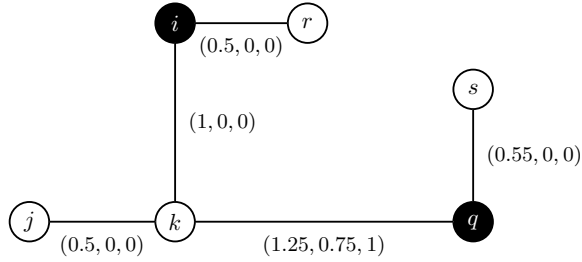
\begin{figure}[htb]
	\centering 
	\resizebox{.35\linewidth}{!}{
	\begin{pspicture}(4.3,1)(11,6)
	\psset{radius=0.2, fillstyle=solid}
	
	\Cnode[radius=1.8ex,fillcolor=black](6.5,5){i} \rput(6.5,5){\small \color{white} $i$}
	\Cnode[radius=1.8ex,](8.5,5){r} \rput(8.5,5){\small $r$}
	\Cnode[radius=1.8ex](4.3,2){j} \rput(4.3,2){\small $j$}
	\Cnode[radius=1.8ex](6.5,2){f} \rput(6.5,2){\small $k$}
	\Cnode[radius=1.8ex,fillcolor=black](11,2){k} \rput(11,2){\small \color{white} $q$}
		\Cnode[radius=1.8ex](11,4){s} \rput(11,4){\small  $s$}
	
	\ncline{-}{i}{f}
	\naput{\footnotesize $(1,0,0)$}
	\ncline{-}{j}{f}
	\nbput{\footnotesize $(0.5,0,0)$}
	\ncline{-}{f}{k}
	\nbput{\footnotesize $(1.25,0.75,1)$}
	\ncline{-}{i}{r}
	\nbput{\footnotesize $(0.5,0,0)$}
	\ncline{-}{k}{s}
	\nbput{\footnotesize $(0.55,0,0)$}
	\end{pspicture}
	}
	\caption{Illustration of incompatibility}
	\label{fig:Example}
\end{figure}

In this case, from constraints \eqref{lm:f2_flow} we obtain that $x_{i}\leq y_{ki}+y_{kj}+x_k.$  Then, $y_{ki}^*=1.$
On the other hand, facility $q$ is the only one that covers node $j,$ then $y_{jq}^*=1.$ Moreover, as the path from node $j$ to node $q$ traverses node $k$, we obtain that $f_{(j,k)}^{jq*}=1.$ Therefore, from constraint \eqref{ec:ddvvconection_1}, we obtain that $y_{kq}^*=1.$ Thus, we have found that these families of constraints are incompatible ($y_{ki}$ and $y_{kq}$ can not take value one simultaneously due to constraints \eqref{ec:zyij}).
\end{example}

Note that the ideas behind constraints \eqref{lm:f2_flow} and \eqref{lm:f2_flow_delta} are practically identical. The reason why constraints \eqref{lm:f2_flow_delta} are not incompatible with \eqref{ec:ddvvconection_1} and \eqref{ec:ddvvconection_2} is that the constraints \eqref{lm:f2_flow_delta} do not force that the node is assigned to the closest service facility before upgrading, instead for given $i,j$ such that $i\neq j$ and $d(i,j)\leq R,$ it enables the node to be assigned to another node whose distance after the most favourable edge length reductions is smaller than or equal to $d(i,j)$. In \cite{EspMarRod12} a detailed description of closest service assignment constraints is given. 

Observe that the variables dropped from the formulation in the preprocessing phase (Propositions \ref{lm:cov} and \ref{prop:nocov}), can also be removed from the valid inequalities presented in this subsection. In the next section, an alternative formulation for this problem is developed. 
	
\section{Path Formulation} \label{sec:Fz}
In this section we present our second formulation for Up-MCLP. It contains fewer variables and constraints than \name. However, this comes at the expense of reducing the scope of preprocessing the model. 

\begin{figure}[htb]
	\centering
	\resizebox{.4\linewidth}{!}{
	\begin{pspicture}(2,0.8)(9.5,3.3)
	\psset{radius=0.2, fillstyle=solid}
	
	\Cnode[radius=1.8ex](2.5,2){i} \rput(2.5,2){\small $i$}
	\Cnode[radius=1.8ex](4.3,2){j} \rput(4.3,2){\small $j$}
	\Cnode[radius=1.8ex,fillcolor=black](6.5,2){f} \rput(6.5,2){\small \color{white} $q$}
	\Cnode[radius=1.8ex](9,2){k} \rput(9,2){\small $k$}
    \rput(2.5,2.5){\small $z_{ij}=1$}
	\rput(4.3,2.5){\small $z_{jq}=1$}
	\rput(6.5,2.5){\small $x_{q}=1$}
	\rput(9,2.5){\small $z_{kq}=1$}
	
	\pcline{|-|}(2.5,1.4)(6.47,1.4) 
	\nbput{\small $\leq R$}
	\pcline{|-|}(6.53,1.4)(8.97,1.4) 
	\nbput{\small $\leq R$}
	
	\ncline{->}{i}{j}
	\ncline{<-}{f}{j}
	\ncline{<-}{f}{k}
	\end{pspicture}
	}
	\caption{Illustration of $z$-variables}
	\label{fig:ProbDescr:Variables}
\end{figure}
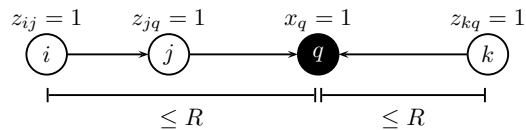

This formulation again models paths of length at most $R$ from a customer node $i$ to a service provider. However, in contrast to \name, the path from $i$ is not modelled as a flow but through the immediate successor of $i$ on a path of length $\le R$ to a facility.
For this purpose, we introduce two new binary variables $z_{ij}$ ($z_{ji}$) for $[i,j]\in E,$ such that $z_{ij}$ ($z_{ji}$) is equal to one if node $j$ ($i$) is the next node on a path of length at most $R$ from $i$ ($j$) to a service facility. In Figure~\ref{fig:ProbDescr:Variables} we illustrate this family of variables where the dark node represents a facility. If $i$ is covered by a facility at $q$, then also $j$ must be covered. Note that a feasible solution resembles a forest rooted at the facilities. 

For the sake of clarity, a description of the decision variables used in the formulation is given next.
\newpage

\noindent\textbf{Decision variables}\\
\begin{tabular}{lp{0.85\textwidth}}
$x_j$&   1, 	if there is a facility at node $j$, and $0$, otherwise, for $j\in V$. \\ 
    $z_{ij}$& 1, if node $j$ is the next node on a path of length $\le R$ from $i$ to a facility, and $0$, otherwise, for $[i,j]\in E.$ \\
	$d_i$&  An upper bound of the length of the built path from node $i$ to its assigned service facility, for $i\in V$.\\ 
	$\delta_e=\delta_{[i,j]}$&  The amount of reduction of the length of edge $e=[i,j]$, for $e\in E.$
\end{tabular}	
\vspace{0.2cm}

The formulation for problem Up-MCLP using these variables, \namet, is as follows:
\begin{alignat}{4}
&\namet&\quad & \max&  \quad & \sum_{i\in V}w_i\left(x_i + \sum_{j:[i,j]\in E}  z_{ij}\right) & & \nonumber \\
&& &\mbox{s.t.} &  \quad & \mbox{\labelcref{ec:B,ec:pfacility,ec:xy_bin}},\nonumber\\
&&& & & \sum_{j:[i,j]\in E}z_{ij} +x_i\leq 1,   & \quad  & i\in V, \label{eq:asig}\\
&&& & & \sum_{j:[i,j]\in E, j\neq k}z_{ij}+x_i\geq z_{ki},   & \quad  &  [k,i]\in E, \label{eq:trans} \\
&&& & & 0 \leq d_i \:\le\: R\sum_{j:[i,j]\in E} z_{ij} & \quad  & i\in V,\label{eq:diRFz} \\
&&& & & d_i\geq d_j+\ell_{[i,j]}z_{ij}-\delta_{[i,j]}-R(1-z_{ij}),   & \quad  & [i,j]\in E, \label{eq:distance} \\
&&& & &  0\leq \delta_{e}\leq u_{e}(z_{ij}+z_{ji}),   & \quad  & e=[i,j]\in E, \label{eq1:ub} \\ 
&&& & & z_{ij}\in\{0,1\}, & \quad  & [i,j]\in E. \label{eq:zent}
\end{alignat} 

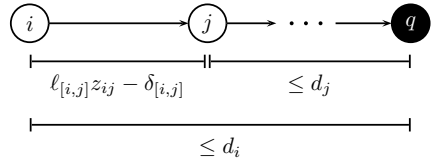
\begin{figure}[htb]
	\centering
	\resizebox{.35\linewidth}{!}{
	\begin{pspicture}(3,0)(10,2.5)
	\psset{radius=0.2, fillstyle=solid}
	
	\Cnode[radius=1.8ex](3.5,2){i} \rput(3.5,2){\small $i$}
	\Cnode[radius=1.8ex](6.3,2){j} \rput(6.3,2){\small $j$}
	\Cnode[radius=1.8ex,fillcolor=black](9.5,2){f} \rput(9.5,2){\small \color{white} $q$}
    \Cnode[radius=0.1ex,linecolor=white](8.3,2){p1}
    \Cnode[radius=0.1ex](8.08,2){p2}
    \Cnode[radius=0.1ex](7.85,2){p3}
    \Cnode[radius=0.1ex](7.625,2){p5}
    \Cnode[radius=0.1ex,linecolor=white](7.4,2){p4}
	\pcline{|-|}(3.5,1.4)(6.27,1.4) \nbput{\small $\ell_{[i,j]}z_{ij}-\delta_{[i,j]}$}
	\pcline{|-|}(9.5,1.4)(6.32,1.4) \naput{\small $\leq d_j$}
	\ncline{->}{i}{j}
		\ncline{->}{j}{p4}
	\ncline{<-}{f}{p1}
	\pcline{|-|}(3.5,0.4)(9.5,0.4) \nbput{\small $\leq d_i$}
	\end{pspicture}
	}
	\caption{Illustration of constraints \eqref{eq:distance}}

	\label{fig:ConstDescr:Distance}
\end{figure}

The family of constraints \eqref{eq:asig} states that each node is assigned to at most one facility or this node is itself a service facility. The family of constraints \eqref{eq:trans} ensures that a node $k$ is not assigned to its service facility through a node $i$, unless node $i$ is also covered or a facility itself. The family of constraints \labelcref{eq:distance,eq:diRFz} set the value of $d_i$, a bound on the distance from node $i$ to its facility, if there exists a path of length at most $R$. We note that \eqref{eq:distance} are equivalent to the well-known Miller-Tucker-Zemlin subtour elimination constraints, extended by our edge length reduction variables. In Figure \ref{fig:ConstDescr:Distance}, an illustration of constraints \eqref{eq:distance} is depicted, in which the dark node represents a facility. Finally, the families of constraints \eqref{ec:B} and \eqref{eq1:ub} establish the bounds on the amount of length edge reductions. 

Note that constraints \eqref{eq:asig} and \eqref{eq:trans} ensure that 
\begin{equation}
    z_{ij}+z_{ji}\leq 1, \quad [i,j]\in E. \label{eq:onlyone}
\end{equation}
The following result presents an improvement to the previous formulation, proving that the integrality condition on the $x$-variables can be relaxed, providing a new family of valid inequalities, and strengthening a family of constraints.  
	\begin{prop}
	The formulation \namet can be enhanced as follows:
\begin{enumerate}[i)]
    \item The binary condition for the $x$-variables can be relaxed.
    \item The following  are valid inequalities for \namet.
    \begin{equation}\label{eq:LBdiz}
    d_i\geq \sum_{j:[i,j]\in E} \left(\ell_{[i,j]}-u_{[i,j]}\right) z_{ij}, \quad i\in V.
    \end{equation}
    \item Constraints \eqref{eq:distance} can be reinforced as follows
\begin{equation}\label{eq:distance_ref}
    d_i\geq d_j+\ell_{[i,j]}z_{ij} -\delta_{[i,j]}-R(1-z_{ij})+z_{ji}(R-\ell_{[i,j]})^+,    \quad [i,j]\in E,
\end{equation}
where $a^+:=\max\left\{a,0\right\}$. 
\end{enumerate}	
	\end{prop}
	\ProofNoNL
	The proof of \textit{i)} is very similar to the proof of the first part of Lemma \ref{lm:x_relax}.
	
	Regarding statement \textit{ii)}, the idea behind these constraints is based on the fact that if a non-facility node is covered, the distance from that node to its assigned facility will be at least the length of the adjacent edge in the path to the service provider, minus the maximally allowed edge length reduction, i.e.,
	\begin{equation}\label{eq:LBdiz_0}
    d_i\geq \left(\ell_{[i,j]}-u_{[i,j]}\right) z_{ij}, \quad [i,j]\in E.
\end{equation}
Moreover, each node is linked to at most one other node in the path to its service facility because of constraint \eqref{eq:asig}.  

In order to prove result \textit{iii)}, we analyse the possible cases. Since the $z$-variables are binary and constraints \eqref{eq:onlyone} are satisfied, 
we get the following four possibilities in the optimal solution for $[i,j]\in E$: a) $z_{ij}^*=z_{ji}^*=0$, b) $z_{ij}^*=1, z_{ji}^*=0$,  c.1) $z_{ij}^*=0,$ $z_{ji}^*=1,$ with $R-\ell_{[i,j]}\leq 0,$ and c.2) $z_{ij}^*=0,$ $z_{ji}^*=1,$ with $R-\ell_{[i,j]}> 0$. In cases a), b), and c.1) 
constraints \eqref{eq:distance} are fulfilled. Hence, \eqref{eq:distance_ref} is valid. Therefore, we focus on case c.2). In this case, since the lower bounds of $d_i$ is only given by \eqref{eq:distance}, we can assume \wloge that \eqref{eq:distance} is satisfied with equality, i.e.: 
$$d_j^*=d_i^*+\ell_{[i,j]}-\delta_{[i,j]}^*.$$
Therefore, 
$$d_i^*=d_j^*-\ell_{[i,j]}+\delta_{[i,j]}^*\geq d_j^*-\ell_{[i,j]}-\delta_{[i,j]}^*.$$ 
Hence, since $z_{ij}^*=0,$ $z_{ji}^*=1$, we have that:
$$d_i^*\geq d_j^* +\ell_{[i,j]}z_{ij}^*-\delta_{[i,j]}^* - R(1-z_{ij}^*) +z_{ji}^*(R-\ell_{[i,j]}),$$
and consequently the family of inequalities \eqref{eq:distance_ref} holds.
Finally, this clearly strengthens the family of constraints \eqref{eq:distance}. 
	\EndProofNoNL

In what follows, we will refer to \namet as the formulation where the above proposition has been applied. Next, we present some valid inequalities linking $x$- and $z$-variables. These inequalities are designed to strengthen the formulation. In \namet,  it is not possible to represent which service is assigned to a given node. Therefore, the ideas of 
Proposition~\ref{prop:nocov} cannot be used. Although it is possible to obtain valid inequalities for this formulation based on constraints~\eqref{lm:f2_flow}. 
\begin{prop}
The following families of constraints are valid inequalities for \namet:
\begin{alignat}{3}
&x_j \leq \sum_{k:[i,k]\in E, \ell_{[i,k]}-u_{[i,k]} \leq d(i,j) }  z_{ik}+x_i,  &\quad \quad  & i\in V, j \in \hat{V}_{i}, \label{lm:F2_fz}  \\
&\sum_{i\in W}\sum_{j\in W:[i,j]\in E}z_{ij}  \leq |W|-1, && W \subset V, 3\leq|W|\leq n-p, \label{eq:ciclo}\\
&d_i \geq  \sum_{j: [i,j]\in S_1}\hspace{-0.3cm} \left(d_j-R(1-z_{ij})\right)  +\sum_{j:[i,j]\in S_2}\hspace{-0.3cm} \left(\ell_{[i,j]}z_{ij}- \delta_{[i,j]}\right), &\; &i \in V, S_1, S_2\subseteq \Gamma_i. \label{lm:cutdiS}
\end{alignat}
\end{prop}
\ProofNoNL
If a facility is open at some node $j$ whose distance to node $i$ before upgrading the network was lower than or equal to the coverage radius (hypothesis of Proposition~\ref{lm:cov}), we can be sure that node $i$ will be covered by some facility. Therefore, in order to eliminate possible symmetries, we assume that either $i$ is a facility itself or the immediate successor of node $i$ on its path to a service facility is a node whose distance to $i$ after upgrading can be smaller than or equal to $d(i,j)$. Using the above argument, the family of constraints \eqref{lm:F2_fz} is obtained.

Secondly, we can include the valid inequalities \eqref{eq:ciclo} to avoid cycles. These inequalities are not required in \namet because the family of constraints \eqref{eq:distance} or equivalently \eqref{eq:distance_ref} avoid cycles in any feasible solution.  However, they can improve the linear relaxation bounds. 
 
Finally, we prove that constraints \eqref{lm:cutdiS} are valid inequalities. Using constraint \eqref{eq:asig}, we know that in any feasible solution, for each $i \in V,$ there is at most one $j_0\in V,$ $[i,j_0]\in E$ such that $z_{ij_0}=1$. 

On the one hand, if $\sum_{j: [i,j]\in E}z_{ij}=0,$ we obtain that $d_i=0$ by \eqref{eq:diRFz}. Furthermore, this latter set of constraints ensures that $d_j\leq R,$ for $j\in V.$ Then, $d_j- R(1-z_{ij})\leq 0,$ for $[i,j]\in E.$ Moreover, since the $\delta$-variables are non-negative and $z_{ij}=0,$ for $j\in V$, we obtain that $\ell_{[i,j]}z_{ij}- \delta_{[i,j]}\leq 0,$ for $[i,j]\in E.$ Hence, it holds that: 
$$0=d_i\geq  \sum_{j: [i,j]\in S_1} \left(d_j-R(1-z_{ij})\right) + \sum_{j: [i,j]\in S_2} \left(\ell_{[i,j]}z_{ij}  - \delta_{[i,j]}\right),\quad S_1,S_2\subseteq \Gamma_i.$$

On the other hand, if exists $j_0\in V,$ $[i,j_0]\in E,$ such that $z_{ij_0}=1,$ using \eqref{eq:distance} we know that: 
$$d_i\geq d_{j_0} -R(1-z_{ij_0}) +\ell_{[i,j_0]}z_{ij_0} -\delta_{[i,j_0]}.$$
Furthermore, $d_j- R(1-z_{ij})\leq 0,$ for $[i,j]\in E,$ such that $j\neq j_0,$ and $\ell_{[i,j]}z_{ij}- \delta_{[i,j]}\leq 0,$ for $[i,j]\in E,$ such that $j\neq j_0.$ Therefore: 
$$d_i\geq \sum_{j: [i,j]\in S_1} \left(d_j-R(1-z_{ij})\right) + \sum_{j: [i,j]\in S_2} \left(\ell_{[i,j]}z_{ij}  - \delta_{[i,j]}\right),\quad S_1,S_2\subseteq \Gamma_i.$$
Thus, we conclude that constraints~\eqref{lm:cutdiS} are valid inequalities.
\EndProofNoNL

Observe that in preliminary computational experiments the addition of the following constraints from family \eqref{eq:ciclo} as 
cuts in the branching tree was quite effective (more details are provided in Section~\ref{sec:ComputationalResults}):    
\begin{equation} \label{DDVV:cycle3}
z_{ij} + z_{ji} + z_{jk} + z_{kj} + z_{ik} + z_{ki}\leq 2, \quad [i,j],[j,k],[i,k] \in E. 
\end{equation}

Next, we solve the separation problem in the family of constraints \eqref{lm:cutdiS}, i.e., given a solution of the LP-relaxation of the formulation, find one or more constraints in family \eqref{lm:cutdiS} that are not satisfied. Hence, sets $S_1$ and $S_2$ that maximises the right-hand-side of the inequality have to be identified.  Let $\bar{d},$ $\bar{\delta}$, and $\bar{z}$ be the optimal vectors of values of the $d$-, $\delta$-, and $z$-variables, respectively, in a node of the branching tree during the resolution of an instance of formulation \namet. Then, it is straightforward to conclude that one of the following constraints maximises the right-hand-side of \eqref{lm:cutdiS}:
\begin{alignat}{3}
&d_i \geq  \sum_{j: [i,j]\in E, \bar{d}_j>R(1-\bar{z}_{ij})} \hspace{-1cm}\left(d_j-R(1-z_{ij})\right)  +\sum_{j: [i,j]\in E, \ell_{[i,j]}\bar{z}_{ij}> \bar{\delta}_{[i,j]}}\hspace{-1cm}\left(\ell_{[i,j]}z_{ij} - \delta_{[i,j]} \right), &\quad &i \in V,  \label{lm:cutdi}\\
&d_i \geq  \sum_{j: [i,j]\in E, \bar{d}_j + \ell_{[i,j]}\bar{z}_{ij}>R(1-\bar{z}_{ij})+\bar{\delta}_{[i,j]}} \hspace{-1.5cm}\left(d_j + \ell_{[i,j]}z_{ij}-R(1-z_{ij}) -\delta_{[i,j]}\right), &\quad &i \in V.  \label{lm:cutdi2}
\end{alignat}

 In Section \ref{sec:ComputationalResults}, the performance of this formulation and the effectiveness of the valid inequalities will be analysed.

\section{Path-Coverage Formulation} \label{sec:fy}
In this section, we introduce a third formulation, which merges components from the first formulation with the second formulation. More precisely, we add the assignment variables $y$ of \name to \namet. For the sake of clarity, all variables of this formulation are explained below. 

\newpage
\noindent\textbf{Decision variables}\\
\begin{tabular}{lp{0.85\textwidth}}
	$x_j$&   1, if there is a facility at node $j,$ and $0$, otherwise, for $j\in V$.\\
    $y_{ij}$& 1, if node $i$ is assigned to a facility at node $j$,  and $0$, otherwise, for $i,j\in V, i\neq j$. \\
    $z_{ij}$& 1, if node $j$ is the next node on a path of length $\le R$ from $i$ to its service facility, and $0$, otherwise, for $[i,j]\in E$. \\
	$d_i$&  An upper bound of the length of the built path from node $i$ to its service facility, for $i\in V$. \\
	$\delta_e=\delta_{[i,j]}$&  The amount of reduction of the length of edge $e=[i,j]$, for $e\in E$.
\end{tabular}
\hspace{0.5cm}

Next, we present the formulation of problem Up-MCLP using the variables described above:
\begin{alignat}{4}
&\hspace{-1.27cm} \named& \; &\max& & \;\sum_{i\in V}w_i\left(x_i + \sum_{j:[i,j]\in E}  z_{ij}\right) & & \nonumber \\
&& &\mbox{s.t.} & & \mbox{\labelcref{ec:pfacility,ec:ybinary,ec:zxij,ec:xy_bin,ec:B,eq1:ub,eq:diRFz,eq:asig,eq:trans,eq:zent,eq:distance_ref,ec:xy_bin}},\nonumber\\ 
&&& & & \sum_{k\in V\setminus\{i\}}y_{ik}=\sum_{j:[i,j]\in E}z_{ij},   &  & i\in V,\label{eq:sumysumz} \\
&&& & &y_{ik}\geq z_{ij}+z_{ji}+ y_{jk}-1,  && k\in V\setminus\{i,j\},[i,j]\in E, \label{eq:zy} \\
 &&& & &y_{ij}\geq z_{ij}+ z_{ji}+x_j-1, && [i,j]\in E. \label{eq:zy_2}
\end{alignat}

The family of constraints \eqref{eq:sumysumz} establishes that if a node is assigned to a service facility, then there is a path from this node to its facility and vice versa. Constraints \labelcref{eq:zy} ensure that if two nodes are on the same path, i.e., $z_{ij}=1$, they must be assigned to the same facility $k$. Constraints \labelcref{eq:zy_2} represent the particular case where node $j$ hosts a service provider. Observe that the objective function can also be expressed as follows:  {\smaller $$\sum_{i\in V}  w_i \left(x_i+ \sum_{j\in V\setminus\{i\}} y_{ij}\right).$$}
But in preliminary computational experiments, we have found that the objective function with the $z$-variables outperforms the one with the $y$-variables. This analysis was again carried out using the data described in Subsection~\ref{sub:data} and the comparison was based on the number of instances  solved to optimality within the time limit, the time employed to obtain the optimal solutions, the MIP relative gaps, the best solution gaps, and the linear relaxation gaps. More details about these performance values can be found in  Section~\ref{sec:ComputationalResults}.

Below, we present a result that proves that the integrality condition of some families of variables can be relaxed. 

\begin{lemma}
		The binary condition on the $x$-variables and the $y$-variables can be relaxed.
	\end{lemma}
	\ProofNoNL
	The proof of the first part of this lemma is very similar to the proof of the first part of Lemma~\ref{lm:x_relax}.
	Regarding the integrality condition on the $y$-variables, since their values are given by the values of the $z$-variables, it is straightforward to conclude that given an optimal solution it is possible to find another optimal solution in which the $y$-variables are integer.
	\EndProofNoNL
	
In contrast to \namet, this formulation controls the service facilities to which the nodes are assigned with the $y$-variables. This information allows us to use a more sophisticated preprocessing phase. 
For doing so, one of the results presented in Subsection~\ref{subsec:preplemas} is used. Under the hypothesis of Proposition~\ref{prop:nocov}, i.e., a facility at node $i$ will never be assigned to a facility located at node $j$, for $i, j \in V,$ and vice versa, variables $y_{ij}$ and $y_{ji}$ are removed from all the constraints in which they are included (fixed to zero and not included in the formulation to save memory).
 Furthermore, using the information obtained in the preprocessing phase we can develop new valid inequalities, which are discussed in the next subsection.

\subsection{Valid inequalities} \label{sec:validinq}
This subsection is devoted to presenting valid inequalities for formulation \named. We start by remarking that the valid inequalities \eqref{lm:f2_flow_delta} obtained for formulation \name can also be implemented in \named.  Similarly, all the valid inequalities obtained for formulation \namet are still valid for \named, namely, the families of constraints \labelcref{eq:LBdiz_0,eq:LBdiz,eq:ciclo,lm:cutdiS,lm:cutdi,lm:cutdi2,DDVV:cycle3,lm:F2_fz}. However, the additional information provided by the $y$-variables in formulation $\named$ can be used to strengthen some of them. The ones that can be enhanced using the covering variables are described below.

First, the lower bound for the $d$-variables can be improved, i.e., constraint \eqref{eq:LBdiz} can be enhanced as: 
\begin{equation}\label{eq:LBdiy}
    d_i\geq \sum_{j\in V\setminus\{i\}} d(i,j,\delta^{ij})y_{ij}, \quad i\in V.
\end{equation}
Recall that $d(i,j,\delta^{ij})$ represents the distance between nodes $i$ and $j$ using the most favourable edge length reductions satisfying the budget constraint \eqref{ec:B}. In what follows, we will refer to \named as the formulation \named in which constraints \eqref{eq:LBdiy} are included.

Finally, we present a new family of valid inequalities that reinforces constraints \eqref{eq:zy}. The objective of this reinforcement is to improve the resolution of the formulation. It is based on the fact that if two nodes are linked (the sum of their $z$-variables is one), then both nodes will be assigned to the same service facility.
\begin{lemma}
The following are valid inequalities for \named:
\begin{equation}
    z_{ij}+z_{ji}\leq \sum_{k\in W, k\neq i}y_{ik}+ x_i\mathcal{I}_{W}(i)+\sum_{k\notin W, k\neq j}y_{jk}+x_j\left(1-\mathcal{I}_{W}(j)\right), \quad [i,j]\in E, W\subseteq V, \label{eq:cutSepa}
\end{equation}
where $\mathcal{I}_{W}(i)$ is the indicator function, i.e., $\mathcal{I}_{W}(i)=1$ if $i\in W$ and 0 otherwise.
\end{lemma}
Note that, for the case $W=\{k\},$ for $k\in V,$ we obtain 
$z_{ij}+z_{ji}\leq y_{ik}+\sum_{t\in V, t\neq k, t\neq j}y_{jt}+x_j,$
using constraints \eqref{eq:asig} and \eqref{eq:sumysumz}, it holds that $z_{ij}+z_{ji}\leq y_{ik}+\sum_{t\in V, t\neq k, t\neq j}y_{jt}+x_j\leq y_{ik}+1-y_{jk}.$
Hence, some constraints of family \eqref{eq:cutSepa} are tighter than \eqref{eq:zy}. 

As the cardinality of \eqref{eq:cutSepa} is exponential, we solve the separation problem in this family of constraints. Therefore, the set $W$ that minimises the right-hand-side of constraints \eqref{eq:cutSepa} has to be identified. Let $\bar{y}$ ($\bar{x}$) be the optimal vector values of $y$-variables ($x$-variables) in a node of the branching tree during the resolution of an instance of formulation \named. Then, it is straightforward to conclude that the following constraints minimise the right-hand-side of \eqref{eq:cutSepa}.
 \begin{equation}\label{eq:cutSepa1}
z_{ij}+z_{ji}\leq \hspace{-0.7cm}\sum_{k\in V:\bar{y}_{ik}\leq \bar{y}_{jk}, \bar{y}_{ik}\leq \bar{x}_{j}}\hspace{-0.7cm}y_{ik}+x_{i}\mathcal{I}_{\{k:\bar{x}_{k}\leq \bar{y}_{jk}\}}(i) +\hspace{-0.7cm}\sum_{k\in V:\bar{y}_{jk}< \bar{y}_{ik},\bar{y}_{jk}< \bar{x}_{i}}\hspace{-0.6cm}y_{jk}+x_{j}\mathcal{I}_{\{k:\bar{x}_{k}< \bar{y}_{ik}\}}(j), \quad [i,j]\in E. \end{equation}

Note that if a pair of nodes satisfies at least one of the conditions of Proposition~\ref{prop:nocov}, their corresponding $y$-variables can be removed from all the constraints including the valid inequalities presented in this subsection.

In the following section, the performance of the three proposed formulations for \mbox{Up-MCLP} are compared.

\section{Computational Results}\label{sec:ComputationalResults}

In this section, we present the results of several computational experiments which compare the performance of the three proposed formulations and show the improvements achieved thanks to the preprocessing phase and the inclusion of the valid inequalities developed throughout the paper. The experiments were conducted on an Intel(R) Xeon(R) W-2135 CPU 3.70 GHz 32 GB RAM, using CPLEX 20.1.0 in  Concert Technology C++ with a time limit of 1800 seconds. We used the default parameter settings for CPLEX.

Regarding the preprocessing phase for formulations \name and \named, aiming to find a balance between preprocessing time and the quality of the $d(i,j,\delta^{ij})$ bounds for each pair $i,j\in V$, the following strategy has been implemented. First, we computed the matrix of pairwise shortest distances without upgrading and the matrix of pairwise shortest distances after upgrading all edges to their full maximum ($d(i,j,u)$) using the Floyd-Warshall algorithm.
Then, we checked if the hypothesis of Proposition~\ref{lm:cov} or if the hypotheses \textit{i)-iii)} of Proposition~\ref{prop:nocov} are fulfilled. If either of these conditions is satisfied for a pair $i,j\in V$, we removed the corresponding variables and constraints and we used $d(i,j,u)$ as $d(i,j,\delta^{ij})$ in the valid inequalities that are required (as e.g. constraints \eqref{lm:f2_flow_delta}). This could be done because $d(i,j,u)$ is a lower bound of $d(i,j,\delta^{ij})$.
If neither of these conditions were fulfilled for a given pair $i,j\in V$, we solved the linear relaxation of \namedist{j}. The minimum between its optimal objective value and $d(i,j,u)$ is the value that we used as $d(i,j,\delta^{ij})$ in the corresponding valid inequalities. As before, this can be done because both values are lower bounds of $d(i,j,\delta^{ij})$. For sake of clarity, we summarise the proccess in Algorithm~\ref{algo:prep}.

\begin{algorithm2e}[h]
	\DontPrintSemicolon \SetAlFnt{\small\sl}
	\SetAlCapFnt{\small\sl} \AlCapFnt
	\caption{Preprocessing phase}
	\label{algo:prep}
	
	\KwIn{Formulation.}
	\KwOut{Preprocessed formulation.}

 \ForEach{$i,j\in V, i<j$}{
			 Compute $d(i,j,u).$
			 
			  \lIf{the hypothesis of Proposition~\ref{lm:cov} or the hypotheses \textit{i)-iii)} of Proposition~\ref{prop:nocov} are fulfilled.}{
			  	\begin{enumerate} 
			  	\item Apply the proposition removing the corresponding variables and constraints.
			  	\item Set $d(i,j,\delta^{ij})=d(i,j,u)$ if it is required in any constraint.
			  	\end{enumerate}}
			  	\vspace*{-3ex}
			  	\nl \Else{
			  	 Solve the linear relaxation of \namedist{j} with objective value $LRP_{d(i,j,\delta^{ij})}$.
			  	 Set $d(i,j,\delta^{ij})= \max\{d(i,j,u), LRP_{d(i,j,\delta^{ij})}\}$ if it is required in any constraint.
			  	}
			  }
	\BlankLine
	
 \Return\ Preprocessed formulation.\;
\end{algorithm2e}

In preliminary computational experiments, we checked the performance of the alternative formulation \name for Up-MCLP and the valid inequalities to identify which ones performed best in each formulation. After this preliminary choice, we conclude that the best formulations are: 
\begin{enumerate}[a)]
    \item Formulation \name where constraints \eqref{ec:flow1}, \eqref{ec:flow2}, and \eqref{ec:alphabounds} have been replaced with \eqref{ec:flow1_gamma}, \eqref{ec:flow2_gamma}, and  \eqref{ec:gammabounds} and the family of constraints \eqref{ec:flow6} is included, named \name for short. 
    \item Formulation \namet where constraints \eqref{eq:LBdiz} are included and constraints \eqref{eq:distance} have been reinforced by constraints \eqref{eq:distance_ref}. In what follows, we call it formulation \namet. 
    \item Formulation \namet with constraints \eqref{lm:F2_fz} as valid inequalities 
and \eqref{DDVV:cycle3} as particular case of \eqref{eq:ciclo} in a pool of user cuts, named \namet~+~VI for short. 
\item Formulation \named where constraints \eqref{eq:LBdiy} are included, we call it formulation \named. 
\item Formulation \named including constraints \eqref{lm:f2_flow_delta} and \eqref{DDVV:cycle3} as particular case of \eqref{eq:ciclo} in a pool of user cuts, named \named~+~VI for short. 
\end{enumerate}

Regarding the families of valid inequalities that we have used, in preliminary investigation, we observe 
that the family of constraints \eqref{lm:F2_fz} had a better performance than the family of constraints \eqref{DDVV:cycle3} for formulation \namet in the majority of the tested cases. Similarly, in the case of \named, the family of constraints \eqref{lm:f2_flow_delta} tended to provide a greater improvement in performance than the family of constraints \eqref{DDVV:cycle3}. Concerning the families of valid inequalities that we have not included in the reported numerical experiments, we would like to remark that including the constraints \eqref{lm:cutdi}, \eqref{lm:cutdi2}, and \eqref{eq:cutSepa1} in the branching tree was effective, as the number of nodes in which the instances were solved decreased. However, this procedure is time-consuming, so that even though the instances were solved on fewer nodes the overall computation time increased. 

The rest of the section is structured as follows. First, the data used in the computational experiments are described. Second, the advantages of the preprocessing phase are shown. Then, the following subsections compare the different formulations with and without valid inequalities in complete graphs and in sparse graphs, respectively. These subsections illustrate the great value of the preprocessing phase and the addition of valid inequalities. 

\subsection{Data}\label{sub:data}
The computational experiments were carried out on two different types of networks. 

First, we generated  instances adapting the procedure used in
\cite{ReSchoWi08, CorFurLju19}, among others. Nodes were given by points whose coordinates followed a uniform distribution over [0,30]. Then, we computed the complete graph where the length of the edges is the Euclidean distance between the nodes.  We named these instances as ``graph" followed by the number of vertices, e.g., ``graph30" is a complete graph with 30 nodes and 435 edges. 

Secondly, we used the uncapacitated $p$-median datasets from the OR-Library, called pmed, see \cite{ORLIBRARY}. As said in the documentation of these datasets, Floyd's algorithm was applied to obtain a symmetric allocation cost matrix that satisfied the triangle inequality. 
The main difference with respect to the previous datasets is that the $p$-median instances are sparser graphs (the number of edges is $n^2/50$). 

The parameters have been chosen as described below. The number of facilities, $p$, was proportional to the number of vertices, i.e., $p\in\{1, n/10, n/20\}$. The node weights or demands, $w_i$ for $i\in V$, were integers uniform randomly generated between 1 and 100. We tested three different coverage radii, $R$, such that we could cover approximately $50\%$, $60\%$, and $70\%$ of the total demand when solving the maximal covering location problem without upgrading $(DT_{\text{\tiny MCLP}})$, i.e., $R\in\{R(50\%DT_{\text{\tiny MCLP}})$,  $R(60\%{DT_{\text{\tiny MCLP}}})$, $R(70\%DT_{\text{\tiny MCLP}})\}$. Upgrading costs, $c_e,$ for $e\in E,$ were uniform randomly generated between 1 and 3. The upper bounds $u_e$, for $e\in E,$ were uniform randomly generated from $(0, 30\%\ell_e),$ for $e\in E$. Then, the length of the edges was modified as $\ell_e + u_e,$ for $e\in E$. This implies that the  instances satisfy the triangle inequality when the full discount is applied in all edges. Finally, the budget $B$ was computed as follows. First, we sorted the upgrade costs $c_eu_e$, for $e\in E$, in non-increasing order. Let $\rho$ be a permutation of set $E$ such that $c_{e_{\rho(1)}}u_{e_{\rho(1)}} \geq c_{e_{\rho(2)}}u_{e_{\rho(2)}} \geq \ldots \geq c_{e_{\rho(m)}}u_{e_{\rho(m)}}.$ Then, since we are constructing a forest with $p$ components (as seen in Section~\ref{sec:Fz}), we can assume that at most $n-p$ edges will be upgraded. Therefore, we computed the maximum  required  budget for upgrading the most expensive edges,  
{\smaller $$B_{max}=\sum_{t=1}^{n-p} u_{e_{\sigma(t)}}{c}_{e_{\sigma(t)}},$$}
and selected $B\in \{0.5\%B_{max}, 1\%B_{max}, 5\%B_{max}\}.$

\subsection{Preprocessing phase}
In this subsection, we show the enhancements provided by the preprocessing, i.e., Propositions~\ref{lm:cov} and \ref{prop:nocov}.  For doing so, we solve \name and \named with and without preprocessing.  

\afterpage{
\begin{landscape}
\begin{table}[p]
 	\resizebox{\linewidth}{!}{
  \begin{tabular}{|c|c|c|r|rrrr|rrrrrrrr|rrrr|rrrrrrrr|}
\cline{1-28}    \multirow{3}[5]{*}{\begin{turn}{90}Data\end{turn}} & \multirow{3}[5]{*}{$B\%$} & \multirow{3}[5]{*}{$p$} & \multicolumn{1}{c|}{\multirow{3}[5]{*}{$R\%$}} & \multicolumn{12}{c|}{\name}                                                                    & \multicolumn{12}{c|}{\named} \\
\cline{5-28}          &       &       &       & \multicolumn{4}{c|}{Without preprocessing} & \multicolumn{8}{c|}{With preprocessing}                        & \multicolumn{4}{c|}{Without preprocessing} & \multicolumn{8}{c|}{With preprocessing} \\
\cline{5-28}          &       &       &       &   \multicolumn{1}{r}{$t_{total}$} & \multicolumn{1}{r}{$G\%$}& \multicolumn{1}{r}{$G_{BS}^t\%$}  & \multicolumn{1}{r|}{$G_{LP}^t$\%} & \multicolumn{1}{l}{$t_{st}$} & \multicolumn{1}{l}{$t_{total}$} & \multicolumn{1}{r}{$G\%$}& \multicolumn{1}{r}{$G_{BS}^t\%$}  & \multicolumn{1}{r}{$G_{LP}^t$\%}& \multicolumn{1}{l}{$R_{c}\%$} & \multicolumn{1}{l}{$R_{v}\%$} & \multicolumn{1}{l|}{$R_{bv}\%$} & \multicolumn{1}{l}{$t_{total}$} &\multicolumn{1}{r}{$G\%$}& \multicolumn{1}{r}{$G_{BS}^t\%$}  & \multicolumn{1}{r|}{$G_{LP}^t$\%} & \multicolumn{1}{l}{$t_{st}$} & \multicolumn{1}{l}{$t_{total}$} & \multicolumn{1}{r}{$G\%$}& \multicolumn{1}{r}{$G_{BS}^t\%$}  & \multicolumn{1}{r}{$G_{LP}^t$\%} & \multicolumn{1}{l}{$R_{c}\%$} & \multicolumn{1}{l}{$R_{v}\%$} & \multicolumn{1}{l|}{$R_{bv}\%$} \\
\hline
    \multicolumn{1}{|c|}{\multirow{27}[17]{*}{\begin{turn}{90}graph40 $|V|=40, |E|=780$\end{turn}}} & \multirow{9}[5]{*}{0.5} & \multirow{3}[1]{*}{1} & 50    & 1042.5 & 0.0(5) & 0.0 & 93.2 & 0.1   & 0.7   & 0.0(5) & 0.0 & 2.3 & 97.0 & 97.0 & 96.9 & 1807.9 & 99.9(0) & 5.7 & 93.2 & 0.1   & 2.3   & 0.0(5) & 0.0 & 2.4 & 69.9 & 46.2 & 52.5 \\
          &       &       & 60    & 1490.0 & 15.4(3) & 0.0 & 52.3 & 0.1   & 3.8   & 0.0(5) & 0.0 & 0.7 & 95.8 & 95.7 & 95.7 & 1804.2 & 74.8(0) & 12.4 & 52.3 & 0.2   & 8.8   & 0.0(5) & 0.0 & 0.7 & 63.1 & 37.6 & 43.5 \\
          &       &       & 70    & 1556.8 & 7.3(4) & 0.0 & 35.4 & 0.2   & 2.0   & 0.0(5) & 0.0 & 1.0 & 96.0 & 96.0 & 95.9 & 1807.5 & 59.6(0) & 13.4 & 35.4 & 0.2   & 724.2 & 1.0(3) & 0.0 & 1.0 & 58.3 & 32.0 & 37.5 \\
\cline{3-28}          &       & \multirow{3}[2]{*}{2} & 50    & 1578.7 & 74.4(1) & 3.9 & 87.0 & 0.1   & 0.4   & 0.0(5) & 0.0 & 1.2 & 98.4 & 98.4 & 98.4 & 61.3  & 0.0(5) & 0.0 & 88.7 & 0.1   & 0.4   & 0.0(5) & 0.0 & 3.3 & 80.8 & 68.7 & 74.3 \\
          &       &       & 60    & 1801.2 & 64.1(0) & 4.2 & 57.0 & 0.1   & 0.4   & 0.0(5) & 0.0 & 0.9 & 97.8 & 97.8 & 97.7 & 1093.2 & 4.3(3) & 0.0 & 57.0 & 0.1   & 0.4   & 0.0(5) & 0.0 & 0.9 & 77.5 & 60.0 & 66.2 \\
          &       &       & 70    & 1801.7 & 40.5(0) & 4.4 & 34.1 & 0.1   & 0.7   & 0.0(5) & 0.0 & 1.8 & 97.9 & 97.8 & 97.8 & 1769.3 & 17.5(1) & 1.2 & 34.1 & 0.1   & 2.4   & 0.0(5) & 0.0 & 1.8 & 74.3 & 53.1 & 59.5 \\
\cline{3-28}          &       & \multirow{3}[2]{*}{4} & 50    & 2.9   & 0.0(5) & 0.0 & 39.4 & 0.1   & 0.1   & 0.0(5) & 0.0 & 1.0 & 99.1 & 99.1 & 99.0 & 0.3   & 0.0(5) & 0.0 & 67.2 & 0.1   & 0.1   & 0.0(5) & 0.0 & 2.8 & 86.7 & 83.4 & 87.3 \\
          &       &       & 60    & 27.5  & 0.0(5) & 0.0 & 48.6 & 0.1   & 0.2   & 0.0(5) & 0.0 & 0.8 & 99.0 & 98.9 & 98.9 & 8.4   & 0.0(5) & 0.0 & 56.8 & 0.1   & 0.2   & 0.0(5) & 0.0 & 1.7 & 84.9 & 78.5 & 83.1 \\
          &       &       & 70    & 1546.9 & 25.6(1) & 0.6 & 34.3 & 0.1   & 1.0   & 0.0(5) & 0.0 & 3.9 & 98.3 & 98.2 & 98.2 & 305.4 & 0.0(5) & 0.0 & 34.3 & 0.1   & 0.6   & 0.0(5) & 0.0 & 4.4 & 81.9 & 70.8 & 76.3 \\
\cline{2-28}          & \multirow{9}[6]{*}{1} & \multirow{3}[2]{*}{1} & 50    & 1384.9 & 18.2(4) & 0.1 & 86.0 & 0.1   & 1.1   & 0.0(5) & 0.0 & 0.4 & 96.0 & 96.0 & 96.0 & 1807.3 & 98.4(0) & 6.6 & 86.0 & 0.1   & 3.6   & 0.0(5) & 0.0 & 0.5 & 68.8 & 45.5 & 51.7 \\
          &       &       & 60    & 1611.6 & 17.1(3) & 0.0 & 51.8 & 0.1   & 1.5   & 0.0(5) & 0.0 & 0.4 & 95.2 & 95.2 & 95.1 & 1809.0 & 71.3(0) & 10.7 & 51.8 & 0.1   & 5.4   & 0.0(5) & 0.0 & 0.4 & 62.5 & 37.2 & 43.1 \\
          &       &       & 70    & 1801.9 & 27.9(0) & 0.0 & 34.0 & 0.2   & 2.6   & 0.0(5) & 0.0 & 1.5 & 94.7 & 94.7 & 94.7 & 1816.8 & 58.6(0) & 13.8 & 34.0 & 0.2   & 135.3 & 0.0(5) & 0.0 & 1.5 & 56.9 & 31.3 & 36.7 \\
\cline{3-28}          &       & \multirow{3}[2]{*}{2} & 50    & 1579.3 & 69.8(1) & 1.9 & 84.9 & 0.1   & 0.4   & 0.0(5) & 0.0 & 0.9 & 98.1 & 98.1 & 98.0 & 95.8  & 0.0(5) & 0.0 & 86.5 & 0.1   & 0.3   & 0.0(5) & 0.0 & 1.4 & 80.4 & 68.4 & 74.1 \\
          &       &       & 60    & 1801.4 & 63.8(0) & 6.1 & 53.2 & 0.1   & 0.8   & 0.0(5) & 0.0 & 2.8 & 96.8 & 96.8 & 96.8 & 1609.6 & 12.4(1) & 0.0 & 53.2 & 0.1   & 1.7   & 0.0(5) & 0.0 & 2.8 & 76.5 & 59.3 & 65.4 \\
          &       &       & 70    & 1801.6 & 43.3(0) & 7.6 & 31.6 & 0.1   & 1.9   & 0.0(5) & 0.0 & 3.0 & 96.7 & 96.7 & 96.6 & 1804.5 & 23.9(0) & 1.2 & 31.6 & 0.1   & 41.9  & 0.0(5) & 0.0 & 3.1 & 73.1 & 52.3 & 58.6 \\
\cline{3-28}          &       & \multirow{3}[2]{*}{4} & 50    & 2.9   & 0.0(5) & 0.0 & 39.4 & 0.1   & 0.1   & 0.0(5) & 0.0 & 1.0 & 99.1 & 99.1 & 99.0 & 0.3   & 0.0(5) & 0.0 & 67.2 & 0.1   & 0.1   & 0.0(5) & 0.0 & 2.8 & 86.7 & 83.4 & 87.3 \\
          &       &       & 60    & 27.6  & 0.0(5) & 0.0 & 48.6 & 0.1   & 0.2   & 0.0(5) & 0.0 & 0.8 & 99.0 & 98.9 & 98.9 & 8.4   & 0.0(5) & 0.0 & 56.8 & 0.1   & 0.2   & 0.0(5) & 0.0 & 1.7 & 84.9 & 78.5 & 83.1 \\
          &       &       & 70    & 1547.3 & 25.7(1) & 0.6 & 34.3 & 0.1   & 1.0   & 0.0(5) & 0.0 & 3.9 & 98.3 & 98.2 & 98.2 & 301.5 & 0.0(5) & 0.0 & 34.3 & 0.1   & 0.6   & 0.0(5) & 0.0 & 4.4 & 81.9 & 70.8 & 76.3 \\
\cline{2-28}          & \multirow{9}[6]{*}{5} & \multirow{3}[2]{*}{1} & 50    & 1805.7 & 88.6(0) & 3.3 & 82.1 & 0.1   & 1.7   & 0.0(5) & 0.0 & 0.3 & 95.0 & 95.0 & 95.0 & 1807.6 & 87.3(0) & 3.3 & 82.1 & 0.1   & 0.6   & 0.0(5) & 0.0 & 0.3 & 67.7 & 44.9 & 51.0 \\
          &       &       & 60    & 1805.2 & 46.(0) & 0.4 & 45.4 & 0.2   & 2.1   & 0.0(5) & 0.0 & 0.5 & 94.3 & 94.3 & 94.3 & 1807.6 & 68.3(0) & 12.8 & 45.4 & 0.1   & 16.3  & 0.0(5) & 0.0 & 0.5 & 61.6 & 36.7 & 42.5 \\
          &       &       & 70    & 1802.5 & 43.1(0) & 7.4 & 31.0 & 0.2   & 1.9   & 0.0(5) & 0.0 & 0.0 & 94.2 & 94.2 & 94.1 & 1811.1 & 44.1(0) & 7.9 & 31.0 & 0.2   & 3.6   & 0.0(5) & 0.0 & 0.0 & 56.3 & 31.0 & 36.3 \\
\cline{3-28}          &       & \multirow{3}[2]{*}{2} & 50    & 1800.6 & 106.7(0) & 12.8 & 79.8 & 0.1   & 0.4   & 0.0(5) & 0.0 & 1.1 & 97.4 & 97.4 & 97.3 & 211.0 & 0.0(5) & 0.0 & 81.6 & 0.1   & 0.4   & 0.0(5) & 0.0 & 1.1 & 79.7 & 67.9 & 73.5 \\
          &       &       & 60    & 1801.3 & 61.7(0) & 8.6 & 47.6 & 0.1   & 1.1   & 0.0(5) & 0.0 & 0.3 & 96.1 & 96.0 & 96.0 & 1804.8 & 14.4(0) & 0.0 & 47.6 & 0.1   & 0.7   & 0.0(5) & 0.0 & 0.3 & 75.7 & 58.7 & 64.8 \\
          &       &       & 70    & 1802.0 & 41.7(0) & 9.9 & 27.6 & 0.1   & 1.5   & 0.0(5) & 0.0 & 0.0 & 95.5 & 95.4 & 95.4 & 1804.3 & 21.3(0) & 1.5 & 27.6 & 0.1   & 2.0   & 0.0(5) & 0.0 & 0.0 & 71.8 & 51.4 & 57.6 \\
\cline{3-28}          &       & \multirow{3}[2]{*}{4} & 50    & 3.1   & 0.0(5) & 0.0 & 36.1 & 0.1   & 0.1   & 0.0(5) & 0.0 & 1.2 & 98.9 & 98.8 & 98.8 & 0.3   & 0.0(5) & 0.0 & 63.4 & 0.1   & 0.1   & 0.0(5) & 0.0 & 1.3 & 86.4 & 83.1 & 87.1 \\
          &       &       & 60    & 512.0 & 0.0(5) & 0.0 & 45.9 & 0.1   & 0.2   & 0.0(5) & 0.0 & 0.0 & 98.5 & 98.4 & 98.4 & 6.4   & 0.0(5) & 0.0 & 54.1 & 0.1   & 0.2   & 0.0(5) & 0.0 & 0.0 & 84.4 & 78.1 & 82.7 \\
          &       &       & 70    & 1801.1 & 35.4(0) & 6.0 & 27.1 & 0.1   & 1.0   & 0.0(5) & 0.0 & 1.5 & 97.4 & 97.3 & 97.3 & 463.7 & 1.4(4) & 0.0 & 27.1 & 0.1   & 0.8   & 0.0(5) & 0.0 & 1.5 & 80.8 & 70.1 & 75.5 \\
   \hline
   \hline
    \end{tabular}%
    }
\caption{Performance of formulations \name and \named with and without preprocessing}
  \label{tab:preprocessing}%
\end{table}%
\end{landscape}
}

As an illustrative example, we include the results for graph40 in Table \ref{tab:preprocessing}. 
The performance was similar in the rest of the datasets. The results are the average over five instances generated with the same procedure, varying only the random seed for the generator. The first column indicates the name of the dataset, the number of nodes and the number of edges. Next, the percentage of the maximal budget ($B\%$), and the number of located facilities are depicted ($p$), followed by the approximate percentage of covered demand in the MCLP without upgrading using this radius ($R\%$).
The following four columns describe information about \name without preprocessing. The first one shows the average  time (in seconds) of solving the corresponding five instances. Observe that if any of these instances is not solved to optimality, 1800 seconds is considered as its solution time to compute this average. Then, the following column of this group depicts the MIP relative gap reported by CPLEX ($G\%$) and in brackets the number of instances solved to optimality within the time limit. Next, it is provided the best solution gap, ($G_{BS}^t$\%), computed as follows:
{\smaller $$\text{G}_{BS}^{t}\%=\dfrac{\text{ BS}^{t}-\text{BS}}{\text{ BS}^{t}}\cdot 100,$$} where BS is the best MIP objective value found within the time limit by the formulation  and $\text{BS}^t$ is the best MIP solution value found within the time limit across all formulations. Finally, it is shown the linear relaxation gap, ($\text{G}_{LP}^{t}$\%), computed as follows: {\smaller$$\text{G}_{LP}^{t}\%=\dfrac{\text{ LP}-\text{ BS}^{t}}{\text{ BS}^{t}}\cdot 100,$$} where LP is the optimal solution value of the linear relaxation of the formulation. Note that $\text{G}_{LP}^{t}$\% enables us to compare the linear relaxation of the formulations with each other (it could be possible that $G\%$ is greater than $\text{G}_{LP}^{t}\%$). 
The following blocks of columns depict information about the rest of formulations, \name with preprocessing and \named with and without preprocessing. Observe that the blocks corresponding to formulations with preprocessing include eight columns. The first column of these blocks reports the average time of the preprocessing phase in seconds, the next one shows the average total time (in seconds) of solving the corresponding five instances including the preprocessing time. The third, the fourth, and the fifth columns report the average $G\%,$ $G_{BS}^t,$ and $G_{LP}^t$ respectively. Then, the average percentage of reduction in the number of constraints ($R_{c}\%$), variables ($R_{v}\%$), and binary variables ($R_{bv}\%$) are depicted. The percentage of reduction in the number of constrains is computed for formulation \name as follows: 
{\smaller $$R_{c}\%\hspace{-0.05cm}=\hspace{-0.05cm}\dfrac{\hspace{-0.1cm}\# \text{constraint of } \name\text{ without prep.} -\# \text{constraint of } \name \text{ with prep.}}{\# \text{constraint of } \name\text{ without prep.}}\cdot 100.$$}
The others are calculated analogously. For interested readers, a similar table that compares the performance of the different phases of the preprocessing can be found in supplementary material.

As can be appreciated in Table~\ref{tab:preprocessing}, the preprocessing phase yields a huge reduction in computation time. For example, using formulation \name without preprocessing, only 53 instances are solved to optimality within the time limit, while using formulation \name with preprocessing, all instances are solved (135) and the average time of solving each instance (including the time of preprocessing) is 1.1 seconds.  Furthermore, the reduction in the number of constraints, variables and binary variables is also very large, approximately 97\% on average. In formulation \named, the reduction of time in the resolution process and the reduction of the size of the problem are also very substantial. 

Based on the above results, we conclude that the preprocessing phase presented in the paper is extremely useful and effective. Therefore, in the subsequent computational experiments, the preprocessing phase is included.

\subsection{Results for complete graphs}

In this subsection, we compare the proposed formulations for complete graphs highlighting the effectiveness of the valid inequalities developed. For interested readers, non-parametric tests (Friedman  test  and  Post-Hoc Holland  adjust) to assess the statistical significance of the comparison among the different formulations can be found in supplementary material. We used the shiny application shinytest\footnote{\url{https://github.com/JacintoCC/shinytests}}, see \cite{Nonparametrictest} for further details. The non-parametric tests show that there are significant differences between the formulations presented.

The results of the smaller datasets (graph30 and graph40) are depicted in Table \ref{tab:g3040}. As before, the provided results are the average over five instances generated with the same procedure, varying only the random seed for the generator. The table describes information about \name, \namet, \namet~+~VI and \named, and its structure is similar to that of Table \ref{tab:preprocessing}. Observe that the blocks corresponding to the \namet and \namet~+~VI formulations have no preprocessing time, since we did not provide a preprocessing for these formulations. Note also that the results of \named+~VI are not included because they are really similar to \named. The differences between these formulations will be shown in instances with a larger number of nodes and edges. Moreover, since several of the valid inequalities are included as cuts in the branching tree, the linear relaxation gap $\text{G}_{LP}^{t}$\% of the formulation with valid inequalities is practically the same as the one for the formulation without valid inequalities. Therefore, they do not appear again in the table. 

Finally, the formulation that provided the smallest average total time is highlighted. If any of the five instances were not solved to optimality, the formulation that solved more instances is shown in bold.  
\begin{table}[htp]
  \centering
   \resizebox{\hsize}{!}{
    \begin{tabular}{|c|c|c|r|rrrrr|rrrr|rrr|rrrrr|}
    \hline
    \multirow{2}[4]{*}{Data} & \multirow{2}[4]{*}{$B\%$} & \multirow{2}[4]{*}{$p$} & \multicolumn{1}{c|}{\multirow{2}[4]{*}{$R\%$}} & \multicolumn{5}{c|}{\name}     & \multicolumn{4}{c|}{\namet}       & \multicolumn{3}{c|}{\namet~+~VI}    & \multicolumn{5}{c|}{\named} \\
\cline{5-21}          &       &       &       & \multicolumn{1}{r}{$t_{st}$} & \multicolumn{1}{r}{$t_{total}$} & \multicolumn{1}{r}{$G\%$}& \multicolumn{1}{r}{$G_{BS}^t\%$}  & \multicolumn{1}{r|}{$G_{LP}^t$\%} &   \multicolumn{1}{r}{$t_{total}$} & \multicolumn{1}{r}{$G\%$} &  \multicolumn{1}{r}{$G_{BS}^t\%$}&\multicolumn{1}{l|}{$G_{LP}^t$\%} &  \multicolumn{1}{r}{$t_{total}$} & \multicolumn{1}{r}{$G\%$}& \multicolumn{1}{r|}{$G_{BS}^t\%$}   & \multicolumn{1}{r}{$t_{st}$} & \multicolumn{1}{r}{$t_{total}$} & \multicolumn{1}{r}{$G\%$}& \multicolumn{1}{r}{$G_{BS}^t\%$} & \multicolumn{1}{r|}{$G_{LP}^t$\%} \\
    \hline
    \hline
    \multicolumn{1}{|c|}{\multirow{27}[18]{*}{\begin{turn}{90} graph30, $|V|=30, |E|=435$\end{turn}}}& \multirow{9}[5]{*}{0.5} & \multirow{3}[1]{*}{1} & 50    & \textbf{0.1} & \textbf{0.5} & \textbf{0.0(5)} & \textbf{0.0} & \textbf{2.2} & 54.6  & 0.0(5) & 0.0 & 84.8 & 13.2  & 0.0(5) & 0.0 & 0.1   & 0.7   & 0.0(5) & 0.0 & 2.2 \\
          &       &       & 60    & \textbf{0.1} & \textbf{0.4} & \textbf{0.0(5)} & \textbf{0.0} & \textbf{0.9} & 533.2 & 6.1(4) & 0.0 & 50.9 & 144.1 & 0.0(5) & 0.0 & 0.1   & 1.4   & 0.0(5) & 0.0 & 0.9 \\
          &       &       & 70    & \textbf{0.1} & \textbf{0.5} & \textbf{0.0(5)} & \textbf{0.0} & \textbf{0.7} & 1162.9 & 14.4(2) & 0.2 & 34.7 & 865.0 & 8.2(3) & 0.0 & 0.1   & 0.8   & 0.0(5) & 0.0 & 0.7 \\
\cline{3-21}          &       & \multirow{3}[2]{*}{2} & 50    & \textbf{0.0} & \textbf{0.1} & \textbf{0.0(5)} & \textbf{0.0} & \textbf{0.6} & 0.3   & 0.0(5) & 0.0 & 83.4 & 0.3   & 0.0(5) & 0.0 & \textbf{0.0} & \textbf{0.1} & \textbf{0.0(5)} & \textbf{0.0} & \textbf{0.6} \\
          &       &       & 60    & \textbf{0.0} & \textbf{0.1} & \textbf{0.0(5)} & \textbf{0.0} & \textbf{0.0} & 0.7   & 0.0(5) & 0.0 & 54.3 & 0.5   & 0.0(5) & 0.0 & \textbf{0.0} & \textbf{0.1} & \textbf{0.0(5)} & \textbf{0.0} & \textbf{0.0} \\
          &       &       & 70    & \textbf{0.1} & \textbf{0.3} & \textbf{0.0(5)} & \textbf{0.0} & \textbf{1.0} & 5.0   & 0.0(5) & 0.0 & 36.0 & 2.9   & 0.0(5) & 0.0 & \textbf{0.0} & \textbf{0.3} & \textbf{0.0(5)} & \textbf{0.0} & \textbf{1.1} \\
\cline{3-21}          &       & \multirow{3}[2]{*}{3} & 50    & \textbf{0.0} & \textbf{0.1} & \textbf{0.0(5)} & \textbf{0.0} & \textbf{0.4} & \textbf{0.1} & \textbf{0.0(5)} & \textbf{0.0} & \textbf{79.2} & \textbf{0.1} & \textbf{0.0(5)} & \textbf{0.0} & \textbf{0.0} & \textbf{0.1} & \textbf{0.0(5)} & \textbf{0.0} & \textbf{3.0} \\
          &       &       & 60    & \textbf{0.0} & \textbf{0.1} & \textbf{0.0(5)} & \textbf{0.0} & \textbf{2.6} & 0.3   & 0.0(5) & 0.0 & 53.8 & 0.3   & 0.0(5) & 0.0 & \textbf{0.0} & \textbf{0.1} & \textbf{0.0(5)} & \textbf{0.0} & \textbf{2.7} \\
          &       &       & 70    & \textbf{0.0} & \textbf{0.1} & \textbf{0.0(5)} & \textbf{0.0} & \textbf{0.7} & 0.8   & 0.0(5) & 0.0 & 32.0 & 0.7   & 0.0(5) & 0.0 & \textbf{0.0} & \textbf{0.1} & \textbf{0.0(5)} & \textbf{0.0} & \textbf{0.7} \\
\cline{2-21}          & \multirow{9}[6]{*}{1} & \multirow{3}[2]{*}{1} & 50    & \textbf{0.1} & \textbf{0.8} & \textbf{0.0(5)} & \textbf{0.0} & \textbf{4.1} & 102.7 & 0.0(5) & 0.0 & 77.2 & 49.7  & 0.0(5) & 0.0 & 0.1   & 1.9   & 0.0(5) & 0.0 & 4.2 \\
          &       &       & 60    & \textbf{0.1} & \textbf{0.6} & \textbf{0.0(5)} & \textbf{0.0} & \textbf{0.7} & 620.2 & 2.9(4) & 0.0 & 45.3 & 479.4 & 1.3(4) & 0.0 & 0.1   & 2.2   & 0.0(5) & 0.0 & 0.7 \\
          &       &       & 70    & \textbf{0.1} & \textbf{0.8} & \textbf{0.0(5)} & \textbf{0.0} & \textbf{2.0} & 1370.0 & 10.0(3) & 0.0 & 32.4 & 967.0 & 11.9(3) & 0.0 & 0.1   & 4.1   & 0.0(5) & 0.0 & 2.1 \\
\cline{3-21}          &       & \multirow{3}[2]{*}{2} & 50    & \textbf{0.0} & \textbf{0.1} & \textbf{0.0(5)} & \textbf{0.0} & \textbf{0.6} & 0.3   & 0.0(5) & 0.0 & 83.4 & 0.3   & 0.0(5) & 0.0 & \textbf{0.0} & \textbf{0.1} & \textbf{0.0(5)} & \textbf{0.0} & \textbf{0.6} \\
          &       &       & 60    & \textbf{0.0} & \textbf{0.1} & \textbf{0.0(5)} & \textbf{0.0} & \textbf{0.3} & 0.7   & 0.0(5) & 0.0 & 52.8 & 0.5   & 0.0(5) & 0.0 & \textbf{0.0} & \textbf{0.1} & \textbf{0.0(5)} & \textbf{0.0} & \textbf{0.3} \\
          &       &       & 70    & \textbf{0.1} & \textbf{0.6} & \textbf{0.0(5)} & \textbf{0.0} & \textbf{1.5} & 11.7  & 0.0(5) & 0.0 & 32.5 & 6.1   & 0.0(5) & 0.0 & 0.0   & 0.9   & 0.0(5) & 0.0 & 1.5 \\
\cline{3-21}          &       & \multirow{3}[2]{*}{3} & 50    & \textbf{0.0} & \textbf{0.1} & \textbf{0.0(5)} & \textbf{0.0} & \textbf{0.4} & \textbf{0.1} & \textbf{0.0(5)} & \textbf{0.0} & \textbf{79.2} & \textbf{0.1} & \textbf{0.0(5)} & \textbf{0.0} & \textbf{0.0} & \textbf{0.1} & \textbf{0.0(5)} & \textbf{0.0} & \textbf{3.0} \\
          &       &       & 60    & \textbf{0.0} & \textbf{0.1} & \textbf{0.0(5)} & \textbf{0.0} & \textbf{2.6} & 0.3   & 0.0(5) & 0.0 & 53.8 & 0.3   & 0.0(5) & 0.0 & \textbf{0.0} & \textbf{0.1} & \textbf{0.0(5)} & \textbf{0.0} & \textbf{2.7} \\
          &       &       & 70    & \textbf{0.1} & \textbf{0.1} & \textbf{0.0(5)} & \textbf{0.0} & \textbf{0.7} & 0.8   & 0.0(5) & 0.0 & 32.0 & 0.7   & 0.0(5) & 0.0 & \textbf{0.0} & \textbf{0.1} & \textbf{0.0(5)} & \textbf{0.0} & \textbf{0.7} \\
\cline{2-21}          & \multirow{9}[6]{*}{5} & \multirow{3}[2]{*}{1} & 50    & 0.1   & 0.6   & 0.0(5) & 0.0 & 0.0 & 254.6 & 0.0(5) & 0.0 & 66.2 & 98.3  & 0.0(5) & 0.0 & \textbf{0.1} & \textbf{0.2} & \textbf{0.0(5)} & \textbf{0.0} & \textbf{0.0} \\
          &       &       & 60    & 0.1   & 0.7   & 0.0(5) & 0.0 & 0.0 & 1097.0 & 5.5(3) & 0.0 & 42.9 & 653.4 & 3.7(4) & 0.0 & \textbf{0.1} & \textbf{0.2} & \textbf{0.0(5)} & \textbf{0.0} & \textbf{0.0} \\
          &       &       & 70    & \textbf{0.1} & \textbf{0.8} & \textbf{0.0(5)} & \textbf{0.0} & \textbf{0.0} & 1384.6 & 13.3(2) & 0.0 & 28.8 & 1329.7 & 14.(2) & 0.0 & 0.1   & 1.5   & 0.0(5) & 0.0 & 0.0 \\
\cline{3-21}          &       & \multirow{3}[2]{*}{2} & 50    & 0.0   & 0.2   & 0.0(5) & 0.0 & 2.2 & 0.4   & 0.0(5) & 0.0 & 70.6 & 0.4   & 0.0(5) & 0.0 & \textbf{0.0} & \textbf{0.1} & \textbf{0.0(5)} & \textbf{0.0} & \textbf{3.6} \\
          &       &       & 60    & \textbf{0.0} & \textbf{0.2} & \textbf{0.0(5)} & \textbf{0.0} & \textbf{1.4} & 1.0   & 0.0(5) & 0.0 & 49.4 & 0.8   & 0.0(5) & 0.0 & \textbf{0.0} & \textbf{0.2} & \textbf{0.0(5)} & \textbf{0.0} & \textbf{1.4} \\
          &       &       & 70    & 0.1   & 0.5   & 0.0(5) & 0.0 & 1.2 & 20.4  & 0.0(5) & 0.0 & 30.7 & 11.1  & 0.0(5) & 0.0 & \textbf{0.1} & \textbf{0.4} & \textbf{0.0(5)} & \textbf{0.0} & \textbf{1.2} \\
\cline{3-21}          &       & \multirow{3}[2]{*}{3} & 50    & \textbf{0.0} & \textbf{0.1} & \textbf{0.0(5)} & \textbf{0.0} & \textbf{1.1} & \textbf{0.1} & \textbf{0.0(5)} & \textbf{0.0} & \textbf{72.1} & \textbf{0.1} & \textbf{0.0(5)} & \textbf{0.0} & \textbf{0.0} & \textbf{0.1} & \textbf{0.0(5)} & \textbf{0.0} & \textbf{1.8} \\
          &       &       & 60    & \textbf{0.1} & \textbf{0.2} & \textbf{0.0(5)} & \textbf{0.0} & \textbf{3.2} & 0.5   & 0.0(5) & 0.0 & 47.4 & 0.5   & 0.0(5) & 0.0 & \textbf{0.0} & \textbf{0.2} & \textbf{0.0(5)} & \textbf{0.0} & \textbf{4.6} \\
          &       &       & 70    & 0.1   & 0.3   & 0.0(5) & 0.0 & 4.0 & 1.3   & 0.0(5) & 0.0 & 26.2 & 0.9   & 0.0(5) & 0.0 & \textbf{0.0} & \textbf{0.2} & \textbf{0.0(5)} & \textbf{0.0} & \textbf{4.2} \\
     \hline
    \hline
    \multicolumn{1}{|c|}{\multirow{27}[18]{*}{\begin{turn}{90} graph40, $|V|=40, |E|=780$\end{turn}}}  & \multirow{9}[6]{*}{0.5} & \multirow{3}[2]{*}{1} & 50    & \textbf{0.1} & \textbf{0.7} & \textbf{0.0(5)} & \textbf{0.0} & \textbf{2.3} & 1395.1 & 31.1(2) & 0.0 & 93.2 & 1089.2 & 39.(3) & 6.1 & 0.1   & 2.3   & 0.0(5) & 0.0 & 2.4 \\
          &       &       & 60    & \textbf{0.1} & \textbf{3.8} & \textbf{0.0(5)} & \textbf{0.0} & \textbf{0.7} & 1801.0 & 48.3(0) & 5.1 & 52.3 & 1802.6 & 45.9(0) & 5.8 & 0.2   & 8.8   & 0.0(5) & 0.0 & 0.7 \\
          &       &       & 70    & \textbf{0.2} & \textbf{2.0} & \textbf{0.0(5)} & \textbf{0.0} & \textbf{1.0} & 1801.6 & 43.6(0) & 8.2 & 35.4 & 1803.7 & 62.1(0) & 16.9 & 0.2   & 724.2 & 1.0(3) & 0.0 & 1.0 \\
\cline{3-21}          &       & \multirow{3}[2]{*}{2} & 50    & \textbf{0.1} & \textbf{0.4} & \textbf{0.0(5)} & \textbf{0.0} & \textbf{1.2} & 6.2   & 0.0(5) & 0.0 & 88.6 & 4.5   & 0.0(5) & 0.0 & \textbf{0.1} & \textbf{0.4} & \textbf{0.0(5)} & \textbf{0.0} & \textbf{3.3} \\
          &       &       & 60    & \textbf{0.1} & \textbf{0.4} & \textbf{0.0(5)} & \textbf{0.0} & \textbf{0.9} & 143.6 & 0.0(5) & 0.0 & 57.0 & 53.9  & 0.0(5) & 0.0 & \textbf{0.1} & \textbf{0.4} & \textbf{0.0(5)} & \textbf{0.0} & \textbf{0.9} \\
          &       &       & 70    & \textbf{0.1} & \textbf{0.7} & \textbf{0.0(5)} & \textbf{0.0} & \textbf{1.8} & 878.8 & 2.3(3) & 0.0 & 34.1 & 549.8 & 1.1(4) & 0.0 & 0.1   & 2.4   & 0.0(5) & 0.0 & 1.8 \\
\cline{3-21}          &       & \multirow{3}[2]{*}{4} & 50    & \textbf{0.1} & \textbf{0.1} & \textbf{0.0(5)} & \textbf{0.0} & \textbf{1.0} & 0.2   & 0.0(5) & 0.0 & 67.2 & 0.2   & 0.0(5) & 0.0 & \textbf{0.1} & \textbf{0.1} & \textbf{0.0(5)} & \textbf{0.0} & \textbf{2.8} \\
          &       &       & 60    & \textbf{0.1} & \textbf{0.2} & \textbf{0.0(5)} & \textbf{0.0} & \textbf{0.8} & 0.8   & 0.0(5) & 0.0 & 56.8 & 0.6   & 0.0(5) & 0.0 & \textbf{0.1} & \textbf{0.2} & \textbf{0.0(5)} & \textbf{0.0} & \textbf{1.7} \\
          &       &       & 70    & 0.1   & 1.0   & 0.0(5) & 0.0 & 3.9 & 6.9   & 0.0(5) & 0.0 & 34.3 & 4.4   & 0.0(5) & 0.0 & \textbf{0.1} & \textbf{0.6} & \textbf{0.0(5)} & \textbf{0.0} & \textbf{4.4} \\
\cline{2-21}          & \multirow{9}[6]{*}{1} & \multirow{3}[2]{*}{1} & 50    & \textbf{0.1} & \textbf{1.1} & \textbf{0.0(5)} & \textbf{0.0} & \textbf{0.4} & 1636.2 & 36.3(1) & 3.0 & 86.0 & 1367.6 & 32.2(3) & 3.2 & 0.1   & 3.6   & 0.0(5) & 0.0 & 0.5 \\
          &       &       & 60    & \textbf{0.1} & \textbf{1.5} & \textbf{0.0(5)} & \textbf{0.0} & \textbf{0.4} & 1800.5 & 48.2(0) & 7.1 & 51.8 & 1801.1 & 55.8(0) & 9.7 & 0.1   & 5.4   & 0.0(5) & 0.0 & 0.4 \\
          &       &       & 70    & \textbf{0.2} & \textbf{2.6} & \textbf{0.0(5)} & \textbf{0.0} & \textbf{1.5} & 1801.0 & 42.2(0) & 8.1 & 34.0 & 1802.9 & 48.6(0) & 11.9 & 0.2   & 135.3 & 0.0(5) & 0.0 & 1.5 \\
\cline{3-21}          &       & \multirow{3}[2]{*}{2} & 50    & 0.1   & 0.4   & 0.0(5) & 0.0 & 0.9 & 10.1  & 0.0(5) & 0.0 & 86.5 & 5.9   & 0.0(5) & 0.0 & \textbf{0.1} & \textbf{0.3} & \textbf{0.0(5)} & \textbf{0.0} & \textbf{1.4} \\
          &       &       & 60    & \textbf{0.1} & \textbf{0.8} & \textbf{0.0(5)} & \textbf{0.0} & \textbf{2.8} & 255.7 & 0.0(5) & 0.0 & 53.2 & 168.2 & 0.0(5) & 0.0 & 0.1   & 1.7   & 0.0(5) & 0.0 & 2.8 \\
          &       &       & 70    & \textbf{0.1} & \textbf{1.9} & \textbf{0.0(5)} & \textbf{0.0} & \textbf{3.0} & 1058.5 & 3.1(3) & 0.0 & 31.6 & 817.9 & 1.8(3) & 0.0 & 0.1   & 41.9  & 0.0(5) & 0.0 & 3.1 \\
\cline{3-21}          &       & \multirow{3}[2]{*}{4} & 50    & \textbf{0.1} & \textbf{0.1} & \textbf{0.0(5)} & \textbf{0.0} & \textbf{1.0} & 0.2   & 0.0(5) & 0.0 & 67.2 & 0.2   & 0.0(5) & 0.0 & \textbf{0.1} & \textbf{0.1} & \textbf{0.0(5)} & \textbf{0.0} & \textbf{2.8} \\
          &       &       & 60    & \textbf{0.1} & \textbf{0.2} & \textbf{0.0(5)} & \textbf{0.0} & \textbf{0.8} & 0.8   & 0.0(5) & 0.0 & 56.8 & 0.6   & 0.0(5) & 0.0 & \textbf{0.1} & \textbf{0.2} & \textbf{0.0(5)} & \textbf{0.0} & \textbf{1.7} \\
          &       &       & 70    & 0.1   & 1.0   & 0.0(5) & 0.0 & 3.9 & 6.9   & 0.0(5) & 0.0 & 34.3 & 4.4   & 0.0(5) & 0.0 & \textbf{0.1} & \textbf{0.6} & \textbf{0.0(5)} & \textbf{0.0} & \textbf{4.4} \\
\cline{2-21}          & \multirow{9}[6]{*}{5} & \multirow{3}[2]{*}{1} & 50    & 0.1   & 1.7   & 0.0(5) & 0.0 & 0.3 & 1549.1 & 37.7(1) & 0.0 & 82.1 & 1488.2 & 36.2(1) & 0.0 & \textbf{0.1} & \textbf{0.6} & \textbf{0.0(5)} & \textbf{0.0} & \textbf{0.3} \\
          &       &       & 60    & \textbf{0.2} & \textbf{2.1} & \textbf{0.0(5)} & \textbf{0.0} & \textbf{0.5} & 1800.6 & 50.4(0) & 8.8 & 45.4 & 1806.5 & 44.6(0) & 7.3 & 0.1   & 16.3  & 0.0(5) & 0.0 & 0.5 \\
          &       &       & 70    & \textbf{0.2} & \textbf{1.9} & \textbf{0.0(5)} & \textbf{0.0} & \textbf{0.0} & 1802.9 & 38.6(0) & 7.2 & 31.0 & 1803.4 & 41.2(0) & 8.6 & 0.2   & 3.6   & 0.0(5) & 0.0 & 0.0 \\
\cline{3-21}          &       & \multirow{3}[2]{*}{2} & 50    & \textbf{0.1} & \textbf{0.4} & \textbf{0.0(5)} & \textbf{0.0} & \textbf{1.1} & 15.7  & 0.0(5) & 0.0 & 81.5 & 10.5  & 0.0(5) & 0.0 & \textbf{0.1} & \textbf{0.4} & \textbf{0.0(5)} & \textbf{0.0} & \textbf{1.1} \\
          &       &       & 60    & 0.1   & 1.1   & 0.0(5) & 0.0 & 0.3 & 499.0 & 0.2(4) & 0.0 & 47.6 & 234.3 & 0.0(5) & 0.0 & \textbf{0.1} & \textbf{0.7} & \textbf{0.0(5)} & \textbf{0.0} & \textbf{0.3} \\
          &       &       & 70    & \textbf{0.1}   & \textbf{1.5}   & \textbf{0.0(5)} & \textbf{0.0} & \textbf{0.0} & 1468.1 & 4.6(1) & 0.1 & 27.6 & 1202.5 & 2.9(2) & 0.1 & 0.1   & 2.0   & 0.0(5) & 0.0 & 0.0 \\
\cline{3-21}          &       & \multirow{3}[2]{*}{4} & 50    & \textbf{0.1} & \textbf{0.1} & \textbf{0.0(5)} & \textbf{0.0} & \textbf{1.2} & 0.2   & 0.0(5) & 0.0 & 63.4 & 0.2   & 0.0(5) & 0.0 & \textbf{0.1} & \textbf{0.1} & \textbf{0.0(5)} & \textbf{0.0} & \textbf{1.3} \\
          &       &       & 60    & \textbf{0.1} & \textbf{0.2} & \textbf{0.0(5)} & \textbf{0.0} & \textbf{0.0} & 0.8   & 0.0(5) & 0.0 & 54.1 & 0.7   & 0.0(5) & 0.0 & \textbf{0.1} & \textbf{0.2} & \textbf{0.0(5)} & \textbf{0.0} & \textbf{0.0} \\
          &       &       & 70    & 0.1   & 1.0   & 0.0(5) & 0.0 & 1.5 & 39.4  & 0.0(5) & 0.0 & 27.1 & 13.5  & 0.0(5) & 0.0 & \textbf{0.1} & \textbf{0.8} & \textbf{0.0(5)} & \textbf{0.0} & \textbf{1.5} \\
    \hline
    \hline
    \end{tabular}%
    }
\caption{Performance of formulations \name, \namet, \namet~+~VI, \named on graph30 and graph40}
  \label{tab:g3040}%
\vspace{1cm}
  \centering
   \resizebox{\hsize}{!}{
    \begin{tabular}{|l|rrrr|rrrr|rrrr|rrrr|}
\cline{2-17}    \multicolumn{1}{r|}{} & \multicolumn{4}{c|}{\name}     & \multicolumn{4}{c|}{\namet}       & \multicolumn{4}{c|}{\namet~+~VI}    & \multicolumn{4}{c|}{\named} \\
\cline{2-17}    \multicolumn{1}{r|}{} & \multicolumn{1}{r}{$c$} & \multicolumn{1}{r}{$v$} & \multicolumn{1}{r}{$bv$} & \multicolumn{1}{r|}{nodes} & \multicolumn{1}{r}{$c$} & \multicolumn{1}{r}{$v$} & \multicolumn{1}{r}{$bv$} & \multicolumn{1}{r|}{nodes} & \multicolumn{1}{r}{$c$} & \multicolumn{1}{r}{$v$} & \multicolumn{1}{r}{$bv$} & \multicolumn{1}{r|}{nodes} & \multicolumn{1}{r}{$c$} & \multicolumn{1}{r}{$v$} & \multicolumn{1}{r}{$bv$} & \multicolumn{1}{r|}{nodes} \\
    \hline
    graph30 & 8551.5 & 6754.0 & 3441.2 & 0.6   & 716.1 & 382.3 & 245.2 & 1124757.3 & 918.3 & 382.3 & 245.2 & 607817.9 & 2632.1 & 588.2 & 451.1 & 2357.8 \\
\hline    graph40 & 25549.9 & 20155.2 & 10180.8 & 9.5   & 1152.4 & 611.7 & 395.1 & 2661373.2 & 1497.2 & 611.7 & 395.1 & 1724565.8 & 5163.7 & 955.7 & 739.2 & 107560.3 \\
    \hline
    \end{tabular}%
    }
  \caption{Number of constraints, variables, and nodes visited in the branching tree of formulations \name, \namet, \namet~+~VI, \named on graph30 and graph40}
  \label{tab:nodesgraph3040}%
\end{table}%
The results in Table \ref{tab:g3040} show that formulations \name and \named outperform \namet and \namet~+~VI (the resolution times, the number of instances solved, the MIP relative gap, the best solution gap, and the linear relaxation gap of these formulations are worse). Moreover, it is clear from these results that the valid inequalities improve the performance of formulation \namet as shown in the number of instances solved to optimality (208 instances with respect to 217 instances) and the average total time (489 seconds with respect to 416 seconds). However, this improvement is not large enough to make this formulation competitive with respect to formulations \name and \named on complete graphs. Nevertheless, it can be seen that formulations \namet and \namet~+~VI in graph40 solve many more instances than formulations \name and \named without preprocessing (85 and 91 instances with respect to 53 and 64 instances), as can be seen in Table~\ref{tab:preprocessing}. 

In Table~\ref{tab:nodesgraph3040}, the average of the number of constraints ($c$), the number of variables ($v$), the number of binary variables ($bv$), and the number of nodes visited in the branching tree (nodes) for each dataset (graph30 and graph40) and each formulation is reported. A detailed table can be found in the supplementary material. It can be seen that the dimension of \name is much larger than the others. Note also that the inclusion of valid inequalities in formulation \namet decreases considerably the number of nodes used. Observe that \named is the one with intermediate size and number of nodes.

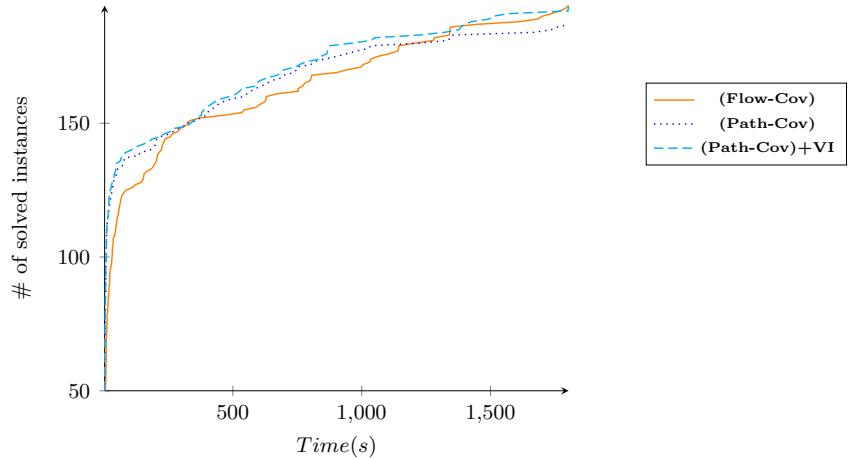
\begin{figure}[hbt] \centering
	\begin{tikzpicture}[scale=0.9,font=\footnotesize]
		\begin{axis}[axis x line=bottom,  axis y line=left,
			xlabel=$Time(s)$,
			ylabel=\# of solved instances,
			legend style={at={(1.6,0.8)}}]
			\addplot[orange,semithick] plot coordinates {		
(5.768,50)
(5.794,51)
(6.045,52)
(6.648,53)
(6.796,54)
(6.851,55)
(7.238,56)
(7.412,57)
(7.675,58)
(7.854,59)
(7.885,60)
(8.021,61)
(8.173,62)
(8.179,63)
(8.201,64)
(8.489,65)
(8.52,66)
(8.542,67)
(8.709,68)
(8.777,69)
(9.318,70)
(10.024,71)
(10.239,72)
(10.699,73)
(10.838,74)
(11.535,75)
(11.559,76)
(11.772,77)
(11.783,78)
(12.774,79)
(14.446,80)
(14.621,81)
(15.245,82)
(15.427,83)
(15.95,84)
(17.594,85)
(18.381,86)
(18.492,87)
(18.841,88)
(19.754,89)
(21.147,90)
(21.893,91)
(22.085,92)
(22.393,93)
(22.843,94)
(23.43,95)
(25.252,96)
(28.952,97)
(29.37,98)
(30.041,99)
(31.039,100)
(31.942,101)
(32.114,102)
(32.479,103)
(33.194,104)
(35.174,105)
(35.215,106)
(35.757,107)
(40.653,108)
(43.348,109)
(45.145,110)
(45.608,111)
(47.519,112)
(49.003,113)
(49.731,114)
(51.875,115)
(55.818,116)
(56.244,117)
(58.381,118)
(60.413,119)
(62.238,120)
(64.154,121)
(66.532,122)
(69.298,123)
(76.34,124)
(85.358,125)
(110.265,126)
(121.644,127)
(145.753,128)
(150.934,129)
(153.924,130)
(154.968,131)
(165.628,132)
(184.081,133)
(192.579,134)
(197.154,135)
(203.893,136)
(204.793,137)
(207.475,138)
(211.994,139)
(222.455,140)
(224.563,141)
(227.048,142)
(232.197,143)
(235.486,144)
(256.63,145)
(263.616,146)
(287.884,147)
(295.972,148)
(313.837,149)
(327.63,150)
(333.683,151)
(371.783,152)
(463.49,153)
(535.709,154)
(543.212,155)
(591.593,156)
(603.132,157)
(623.547,158)
(628.81,159)
(629.303,160)
(668.978,161)
(753.095,162)
(753.313,163)
(771.527,164)
(778.828,165)
(800.548,166)
(803.283,167)
(806.619,168)
(904.41,169)
(943.217,170)
(996.878,171)
(1006.415,172)
(1032.776,173)
(1034.506,174)
(1061.849,175)
(1106.421,176)
(1139.459,177)
(1143.272,178)
(1144.123,179)
(1207.032,180)
(1279.157,181)
(1280.585,182)
(1340.638,183)
(1341.468,184)
(1341.912,185)
(1344.826,186)
(1446.82,187)
(1611.843,188)
(1683.929,189)
(1699.635,190)
(1736.114,191)
(1767.151,192)
(1791.483,193)
(1803,194)
			};
			\addlegendentry{\tiny\textbf{\name}}	
			\addplot[blue,semithick,dotted] plot coordinates {		
(2.401,50)
(2.422,51)
(2.452,52)
(2.462,53)
(2.659,54)
(2.666,55)
(2.703,56)
(2.703,57)
(2.711,58)
(2.767,59)
(3.124,60)
(3.125,61)
(3.19,62)
(3.262,63)
(3.314,64)
(3.379,65)
(3.386,66)
(3.51,67)
(3.525,68)
(3.533,69)
(3.543,70)
(3.65,71)
(3.713,72)
(3.818,73)
(4.024,74)
(4.083,75)
(4.084,76)
(4.233,77)
(4.497,78)
(4.889,79)
(5.048,80)
(5.411,81)
(5.508,82)
(5.647,83)
(5.748,84)
(5.877,85)
(6.029,86)
(6.275,87)
(6.366,88)
(6.419,89)
(6.485,90)
(6.772,91)
(6.803,92)
(6.96,93)
(7.04,94)
(7.48,95)
(7.488,96)
(7.595,97)
(7.988,98)
(7.995,99)
(8.047,100)
(8.087,101)
(8.154,102)
(8.609,103)
(8.672,104)
(9.172,105)
(9.178,106)
(9.477,107)
(9.527,108)
(10.679,109)
(10.96,110)
(11.524,111)
(11.627,112)
(12.877,113)
(14.68,114)
(14.803,115)
(15.389,116)
(17.875,117)
(18.928,118)
(19.566,119)
(23.272,120)
(24.096,121)
(25.375,122)
(28.928,123)
(29.383,124)
(30.279,125)
(32.71,126)
(34.679,127)
(34.741,128)
(41.41,129)
(44.215,130)
(48.677,131)
(49.584,132)
(53.453,133)
(70.5,134)
(70.681,135)
(77.569,136)
(78.99,137)
(136.305,138)
(138.111,139)
(178.933,140)
(189.751,141)
(200.084,142)
(201.685,143)
(203.845,144)
(216.682,145)
(254.297,146)
(269.666,147)
(288.491,148)
(309.861,149)
(312.723,150)
(361.947,151)
(371.108,152)
(399.392,153)
(402.015,154)
(410.005,155)
(430.291,156)
(441.629,157)
(459.91,158)
(489.532,159)
(542.833,160)
(553.959,161)
(568.876,162)
(587.647,163)
(612.038,164)
(633.944,165)
(650.42,166)
(688.665,167)
(696.395,168)
(738.984,169)
(746.411,170)
(748.303,171)
(803.163,172)
(816.746,173)
(853.957,174)
(904.105,175)
(940.551,176)
(977.97,177)
(1026.78,178)
(1035.893,179)
(1228.559,180)
(1336.208,181)
(1343.347,182)
(1350.21,183)
(1667.628,184)
(1735.577,185)
(1759.659,186)
(1803,187)
			};
			\addlegendentry{\tiny\textbf{\named}}	
	\addplot[cyan,semithick,dash pattern=on 4pt off 2pt] plot coordinates {		
(2.316,50)
(2.375,51)
(2.4,52)
(2.406,53)
(2.441,54)
(2.46,55)
(2.466,56)
(2.466,57)
(2.484,58)
(2.518,59)
(2.551,60)
(2.592,61)
(2.631,62)
(2.633,63)
(2.676,64)
(2.857,65)
(2.894,66)
(2.991,67)
(3.104,68)
(3.119,69)
(3.202,70)
(3.243,71)
(3.405,72)
(3.432,73)
(3.545,74)
(3.725,75)
(3.814,76)
(4.027,77)
(4.241,78)
(4.294,79)
(4.318,80)
(4.591,81)
(4.73,82)
(4.823,83)
(4.86,84)
(5.304,85)
(5.434,86)
(5.714,87)
(5.809,88)
(5.883,89)
(6.452,90)
(6.61,91)
(6.73,92)
(6.828,93)
(6.868,94)
(6.917,95)
(7.104,96)
(7.407,97)
(7.607,98)
(7.761,99)
(8.56,100)
(8.951,101)
(8.955,102)
(9.21,103)
(9.281,104)
(9.286,105)
(9.334,106)
(9.509,107)
(10.256,108)
(10.55,109)
(10.606,110)
(11.159,111)
(11.531,112)
(13.897,113)
(14.156,114)
(15.174,115)
(16.227,116)
(16.229,117)
(16.32,118)
(17.42,119)
(18.541,120)
(18.711,121)
(19.423,122)
(19.912,123)
(22.803,124)
(24.104,125)
(28.139,126)
(29.182,127)
(32.406,128)
(34.21,129)
(38.702,130)
(39.985,131)
(44.159,132)
(45.574,133)
(46.392,134)
(47.656,135)
(64.226,136)
(65.879,137)
(80.445,138)
(81.588,139)
(106.185,140)
(118.388,141)
(145.672,142)
(174.185,143)
(184.201,144)
(215.085,145)
(233.054,146)
(258.556,147)
(265.417,148)
(325.429,149)
(335.987,150)
(346.783,151)
(359.49,152)
(377.849,153)
(378.773,154)
(387.914,155)
(397.725,156)
(414.257,157)
(428.414,158)
(448.978,159)
(486.339,160)
(513.292,161)
(526.074,162)
(536.559,163)
(591.66,164)
(592.948,165)
(606.341,166)
(632.572,167)
(672.516,168)
(687.589,169)
(696.724,170)
(741.672,171)
(759.934,172)
(783.066,173)
(837.375,174)
(840.787,175)
(864.666,176)
(867.01,177)
(873.099,178)
(876.081,179)
(964.887,180)
(1036.777,181)
(1049.529,182)
(1224.334,183)
(1320.536,184)
(1377.58,185)
(1385.793,186)
(1398.056,187)
(1426.252,188)
(1495.952,189)
(1507.768,190)
(1569.083,191)
(1800.761,192)
(1801.388,193)
(1803,194)
			};
			\addlegendentry{\tiny\textbf{\named\hspace{-0.1cm}+VI}}			

		\end{axis}
	\end{tikzpicture}
	\caption{Performance profile graph of \#solved instances using \name, \named, and \named~+~VI  formulations on graph100 and graph120}  \label{performanceProfilesbig}
\end{figure}

In Table~\ref{tab:g100120}, a second set of computational experiments is reported. Here, datasets of larger size (graph100 and graph120) are solved so that formulations \name, \named and \named~+~VI can be compared.  Table \ref{tab:g100120} has the same structure as Table \ref{tab:g3040}, but now \named~+~VI is included whereas \namet and \namet~+~VI are not. Moreover, a similar analysis to the one in Table~\ref{tab:nodesgraph3040} is reported in Table~\ref{tab:nodesgraph100120} for graph100 and graph120. A detailed table can be found in the supplementary material. For the purpose of a clearer comparison of these formulations, the performance profile graph of the number of solved instances is depicted in Figure~\ref{performanceProfilesbig}.
\begin{table}[htbp]
  \centering
   \resizebox{\hsize}{!}{
    \begin{tabular}{|c|c|c|r|rrrrr|rrrrr|rrrr|}
    \hline
    \multirow{2}[4]{*}{Data} & \multirow{2}[4]{*}{$B\%$} & \multirow{2}[4]{*}{$p$} & \multicolumn{1}{c|}{\multirow{2}[4]{*}{$R\%$}} & \multicolumn{5}{c|}{\name}     & \multicolumn{5}{c|}{\named}       & \multicolumn{4}{c|}{\named~+~VI}     \\
\cline{5-18}          &       &       &       & \multicolumn{1}{r}{$t_{st}$} & \multicolumn{1}{r}{$t_{total}$} & \multicolumn{1}{r}{$G\%$} &  \multicolumn{1}{r}{$G_{BS}^t\%$}&\multicolumn{1}{r|}{$G_{LP}^t$\%} &  
\multicolumn{1}{r}{$t_{st}$}&  \multicolumn{1}{r}{$t_{total}$} & \multicolumn{1}{r}{$G\%$} &  \multicolumn{1}{r}{$G_{BS}^t\%$}& \multicolumn{1}{r|}{$G_{LP}^t$\%} & \multicolumn{1}{r}{$t_{st}$} & \multicolumn{1}{r}{$t_{total}$} & \multicolumn{1}{r}{$G\%$} & \multicolumn{1}{r|}{$G_{BS}^t\%$} \\
    \hline
    \hline
    \multicolumn{1}{|c|}{\multirow{27}[18]{*}{\begin{turn}{90} graph100, $|V|=100, |E|=4950$\end{turn}}} & \multirow{9}[6]{*}{0.5} & \multirow{3}[2]{*}{1} & 50    & \textbf{6.8} & \textbf{1512.9} & \textbf{5.7(1)} & \textbf{0.0} & \textbf{6.9} & 6.8   & 1807.1 & 10.6(0) & 3.2 & 7.0 & 6.8   & 1564.4 & 7.9(1) & 0.8  \\
          &       &       & 60    & \textbf{7.1} & \textbf{1716.5} & \textbf{1.9(1)} & \textbf{0.0} & \textbf{2.4} & 7.3   & 1807.5 & 3.3(0) & 0.8 & 2.4 & 7.2   & 1807.5 & 2.6(0) & 0.2  \\
          &       &       & 70    & \textbf{8.1} & \textbf{1770.6} & \textbf{1649.1(1)} & \textbf{14.7} & \textbf{3.2} & 7.9   & 1808.3 & 4.4(0) & 1.1 & 3.3 & 7.9   & 1808.4 & 4.8(0) & 1.4  \\
\cline{3-18}          &       & \multirow{3}[2]{*}{5} & 50    & 1.2   & 12.0  & 0.0(5) & 0.0 & 1.6 & 1.2   & 7.5   & 0.0(5) & 0.0 & 1.7 & \textbf{1.2} & \textbf{5.7} & \textbf{0.0(5)} & \textbf{0.0}  \\
          &       &       & 60    & \textbf{1.3} & \textbf{455.3} & \textbf{0.3(4)} & \textbf{0.0} & \textbf{4.8} & 1.3   & 710.2 & 0.0(4) & 0.0 & 5.1 & 1.3   & 654.6 & 0.0(4) & 0.0  \\
          &       &       & 70    & \textbf{1.6} & \textbf{624.0} & \textbf{0.0(4)} & \textbf{0.0} & \textbf{2.5} & 1.6   & 1101.1 & 1.5(2) & 0.3 & 2.5 & 1.6   & 1096.8 & 1.5(2) & 0.5  \\
\cline{3-18}          &       & \multirow{3}[2]{*}{10} & 50    & 1.1   & 2.0   & 0.0(5) & 0.0 & 2.8 & \textbf{1.1} & \textbf{1.4} & \textbf{0.0(5)} & \textbf{0.0} & \textbf{3.8} & \textbf{1.1} & \textbf{1.4} & \textbf{0.0(5)} & \textbf{0.0}  \\
          &       &       & 60    & 1.2   & 8.4   & 0.0(5) & 0.0 & 3.5 & 1.2   & 3.0   & 0.0(5) & 0.0 & 4.0 & \textbf{1.2} & \textbf{2.1} & \textbf{0.0(5)} & \textbf{0.0}\\
          &       &       & 70    & 1.2   & 19.7  & 0.0(5) & 0.0 & 3.3 & 1.2   & 6.5   & 0.0(5) & 0.0 & 3.6 & \textbf{1.2} & \textbf{6.2} & \textbf{0.0(5)} & \textbf{0.0}  \\
\cline{2-18}          & \multirow{9}[6]{*}{1} & \multirow{3}[2]{*}{1} & 50    & 7.2   & 1285.7 & 3.3(2) & 0.0 & 3.6 & \textbf{7.0} & \textbf{1262.4} & \textbf{4.1(2)} & \textbf{0.4} & \textbf{3.7} & 6.9   & 1265.4 & 3.7(2) & 0.0 \\
          &       &       & 60    & \textbf{7.1} & \textbf{1383.9} & \textbf{0.6(3)} & \textbf{0.0} & \textbf{0.6} & 7.4   & 1557.2 & 1.5(2) & 0.8 & 0.6 & 7.3   & 1378.9 & 0.9(2) & 0.3  \\
          &       &       & 70    & \textbf{8.1} & \textbf{1643.4} & \textbf{16.3(2)} & \textbf{11.1} & \textbf{0.9} & 8.1   & 1692.6 & 1.6(1) & 0.7 & 0.9 & 8.1   & 1710.9 & 2.0(1) & 1.1  \\
\cline{3-18}          &       & \multirow{3}[2]{*}{5} & 50    & 1.2   & 16.0  & 0.0(5) & 0.0 & 1.3 & \textbf{1.2} & \textbf{6.7} & \textbf{0.0(5)} & \textbf{0.0} & \textbf{1.3} & 1.2   & 6.9   & 0.0(5) & 0.0 \\
          &       &       & 60    & \textbf{1.3} & \textbf{413.9} & \textbf{0.3(4)} & \textbf{0.0} & \textbf{3.9} & 1.3   & 473.8 & 0.3(4) & 0.0 & 4.1 & 1.3   & 448.4 & 0.3(4) & 0.0 \\
          &       &       & 70    & \textbf{1.6} & \textbf{244.2} & \textbf{0.0(5)} & \textbf{0.0} & \textbf{1.8} & 1.7   & 1091.4 & 0.8(2) & 0.1 & 1.9 & 1.6   & 1087.7 & 1.5(2) & 0.4 \\
\cline{3-18}          &       & \multirow{3}[2]{*}{10} & 50    & 1.1   & 1.9   & 0.0(5) & 0.0 & 2.8 & 1.1   & \textbf{1.4} & \textbf{0.0(5)} & \textbf{0.0} & \textbf{3.8} & \textbf{1.1} & \textbf{1.4} & \textbf{0.0(5)} & \textbf{0.0} \\
          &       &       & 60    & 1.2   & 8.7   & 0.0(5) & 0.0 & 3.5 & 1.2   & 2.8   & 0.0(5) & 0.0 & 3.7 & \textbf{1.2} & \textbf{2.2} & \textbf{0.0(5)} & \textbf{0.0}  \\
          &       &       & 70    & 1.2   & 22.1  & 0.0(5) & 0.0 & 2.7 & 1.2   & 14.8  & 0.0(5) & 0.0 & 2.8 & \textbf{1.2} & \textbf{14.1} & \textbf{0.0(5)} & \textbf{0.0}  \\
\cline{2-18}          & \multirow{9}[6]{*}{5} & \multirow{3}[2]{*}{1} & 50    & 6.9   & 747.2 & 0.0(5) & 0.0 & 0.0 & 7.1   & 284.7 & 0.0(5) & 0.0 & 0.0 & \textbf{7.0} & \textbf{278.3} & \textbf{0.0(5)} & \textbf{0.0}  \\
          &       &       & 60    & 7.3   & 1064.1 & 0.0(5) & 0.0 & 0.0 & 7.3   & 432.2 & 0.0(5) & 0.0 & 0.0 & \textbf{7.2} & \textbf{368.9} & \textbf{0.0(5)} & \textbf{0.0} \\
          &       &       & 70    & 7.9   & 1148.0 & 0.0(5) & 0.0 & 0.0 & \textbf{8.1} & \textbf{513.4} & \textbf{0.0(5)} & \textbf{0.0} & \textbf{0.0} & 8.1   & 632.4 & 0.0(5) & 0.0  \\
\cline{3-18}          &       & \multirow{3}[2]{*}{5} & 50    & 1.2   & 4.0   & 0.0(5) & 0.0 & 0.0 & \textbf{1.2} & \textbf{1.7} & \textbf{0.0(5)} & \textbf{0.0} & \textbf{0.0} & \textbf{1.2} & \textbf{1.7} & \textbf{0.0(5)} & \textbf{0.0}  \\
          &       &       & 60    & 1.3   & 8.6   & 0.0(5) & 0.0 & 0.0 & \textbf{1.3} & \textbf{4.9} & \textbf{0.0(5)} & \textbf{0.0} & \textbf{0.0} & 1.2   & 5.3   & 0.0(5) & 0.0  \\
          &       &       & 70    & \textbf{1.7} & \textbf{21.1} & \textbf{0.0(5)} & \textbf{0.0} & \textbf{0.1} & 1.7   & 403.0 & 0.0(4) & 0.0 & 0.1 & 1.6   & 22.9  & 0.0(5) & 0.0  \\
\cline{3-18}          &       & \multirow{3}[2]{*}{10} & 50    & 1.1   & 1.7   & 0.0(5) & 0.0 & 0.9 & 1.2   & 1.4   & 0.0(5) & 0.0 & 1.0 & \textbf{1.1} & \textbf{1.3} & \textbf{0.0(5)} & \textbf{0.0}  \\
          &       &       & 60    & 1.2   & 3.3   & 0.0(5) & 0.0 & 0.3 & \textbf{1.2} & \textbf{1.5} & \textbf{0.0(5)} & \textbf{0.0} & \textbf{0.3} & \textbf{1.1} & \textbf{1.5} & \textbf{0.0(5)} & \textbf{0.0}  \\
          &       &       & 70    & 1.2   & 7.6   & 0.0(5) & 0.0 & 0.5 & 1.2   & 4.7   & 0.0(5) & 0.0 & 0.5 & \textbf{1.2} & \textbf{5.2} & \textbf{0.0(5)} & \textbf{0.0}  \\
    \hline
    \hline
    \multicolumn{1}{|c|}{\multirow{27}[18]{*}{\begin{turn}{90} graph120, $|V|=120, |E|=7140$\end{turn}}}& \multirow{9}[6]{*}{0.5} & \multirow{3}[2]{*}{1} & 50    & \textbf{10.6} & \textbf{1718.2} & \textbf{6503.6(1)} & \textbf{16.9} & \textbf{3.3} & 10.8  & 1820.2 & 3.6(0) & 0.3 & 3.3 & 10.9  & 1811.3 & 4.2(0) & 0.8  \\
          &       &       & 60    & 12.2  & 1814.7 & 9004.2(0) & 24.8 & 3.3 & 12.1  & 1812.6 & 4.0(0) & 0.6 & 3.3 & 12.2  & 1812.5 & 7.5(0) & 3.3  \\
          &       &       & 70    & 13.3  & 1815.4 & 10993.9(0) & 28.8 & 3.4 & 13.4  & 1814.1 & 3.6(0) & 0.1 & 3.4 & 13.1  & 1814.0 & 4.7(0) & 1.2  \\
\cline{3-18}          &       & \multirow{3}[2]{*}{6} & 50    & 2.1   & 462.0 & 0.7(4) & 0.0 & 4.5 & \textbf{2.1} & \textbf{396.3} & \textbf{0.5(4)} & \textbf{0.0} & \textbf{4.7} & 2.1   & 405.6 & 1.1(4) & 0.1  \\
          &       &       & 60    & \textbf{2.2} & \textbf{1353.5} & \textbf{0.1(3)} & \textbf{0.0} & \textbf{3.5} & 2.2   & 1247.1 & 0.8(2) & 0.1 & 3.6 & 2.2   & 1309.7 & 0.8(2) & 0.1  \\
          &       &       & 70    & \textbf{3.9} & \textbf{1248.9} & \textbf{0.9(2)} & \textbf{0.0} & \textbf{3.7} & 3.9   & 1592.4 & 5.(1) & 1.7 & 3.8 & 3.9   & 1744.7 & 4.9(1) & 1.7 \\
\cline{3-18}          &       & \multirow{3}[2]{*}{12} & 50    & 2.0   & 4.2   & 0.0(5) & 0.0 & 1.7 & \textbf{2.0} & \textbf{2.3} & \textbf{0.0(5)} & \textbf{0.0} & \textbf{2.2} & \textbf{2.0} & \textbf{2.3} & \textbf{0.0(5)} & \textbf{0.0}  \\
          &       &       & 60    & 2.0   & 76.1  & 0.0(5) & 0.0 & 5.1 & 2.0   & 6.1   & 0.0(5) & 0.0 & 5.7 & \textbf{2.0} & \textbf{6.0} & \textbf{0.0(5)} & \textbf{0.0}  \\
          &       &       & 70    & 2.0   & 355.5 & 0.0(5) & 0.0 & 3.4 & 2.1   & 215.1 & 0.0(5) & 0.0 & 3.8 & \textbf{2.0} & \textbf{65.4} & \textbf{0.0(5)} & \textbf{0.0}  \\
\cline{2-18}          & \multirow{9}[6]{*}{1} & \multirow{3}[2]{*}{1} & 50    & 10.7  & 1798.6 & 2019.5(1) & 18.4 & 0.5 & 11.0  & 1796.0 & 1.0(1) & 0.4 & 0.5 & \textbf{10.9} & \textbf{1768.8} & \textbf{0.8(1)} & \textbf{0.3}  \\
          &       &       & 60    & 12.3  & 1815.1 & 9023.2(0) & 30.3 & 1.7 & 12.3  & 1719.2 & 2.2(1) & 0.5 & 1.7 & \textbf{12.3} & \textbf{1728.1} & \textbf{2.2(2)} & \textbf{0.4}  \\
          &       &       & 70    & 13.2  & 1815.4 & 11005.2(0) & 30.0 & 1.8 & 13.7  & 1803.5 & 1.9(1) & 0.1 & 1.8 & \textbf{13.7} & \textbf{1734.1} & \textbf{2.3(2)} & \textbf{0.6} \\
\cline{3-18}          &       & \multirow{3}[2]{*}{6} & 50    & 2.1   & 456.6 & 0.5(4) & 0.0 & 3.9 & 2.1   & 529.4 & 0.5(4) & 0.0 & 4.1 & \textbf{2.1} & \textbf{412.0} & \textbf{0.8(4)} & \textbf{0.2}  \\
          &       &       & 60    & \textbf{2.2} & \textbf{1076.9} & \textbf{0.1(3)} & \textbf{0.0} & \textbf{3.0} & 2.2   & 1274.0 & 0.9(2) & 0.1 & 3.1 & 2.2   & 1323.6 & 0.8(2) & 0.1 \\
          &       &       & 70    & \textbf{4.1} & \textbf{1480.6} & \textbf{1.4(2)} & \textbf{0.0} & \textbf{3.8} & 4.1   & 1806.9 & 5.3(0) & 1.9 & 3.9 & 4.1   & 1655.6 & 5.5(1) & 1.9 \\
\cline{3-18}          &       & \multirow{3}[2]{*}{12} & 50    & 2.0   & 4.2   & 0.0(5) & 0.0 & 1.7 & 2.0   & 2.3   & 0.0(5) & 0.0 & 2.2 & \textbf{1.9} & \textbf{2.2} & \textbf{0.0(5)} & \textbf{0.0}  \\
          &       &       & 60    & 2.0   & 263.1 & 0.0(5) & 0.0 & 5.3 & 2.0   & 7.0   & 0.0(5) & 0.0 & 5.6 & \textbf{2.0} & \textbf{6.2} & \textbf{0.0(5)} & \textbf{0.0}  \\
          &       &       & 70    & 2.1   & 780.9 & 0.2(4) & 0.0 & 3.2 & 2.0   & 380.8 & 0.2(4) & 0.0 & 3.6 & \textbf{2.0} & \textbf{335.9} & \textbf{0.0(5)} & \textbf{0.0}  \\
\cline{2-18}          & \multirow{9}[6]{*}{5} & \multirow{3}[2]{*}{1} & 50    & 10.8  & 1530.5 & 2019.3(2) & 18.5 & 0.0 & \textbf{10.9} & \textbf{367.1} & \textbf{0.0(5)} & \textbf{0.0} & \textbf{0.0} & 11.1  & 506.9 & 0.0(5) & 0.0  \\
          &       &       & 60    & 12.4  & 1816.0 & 7167.1(0) & 31.3 & 0.0 & 12.5  & 971.9 & 5.8(4) & 4.5 & 0.0 & \textbf{12.2} & \textbf{776.0} & \textbf{0.0(5)} & \textbf{0.0} \\
          &       &       & 70    & 13.2  & 1817.0 & 6993.(0) & 31.4 & 0.0 & \textbf{13.6} & \textbf{985.4} & \textbf{0.0(5)} & \textbf{0.0} & \textbf{0.0} & 13.6  & 1057.8 & 0.0(5) & 0.0  \\
\cline{3-18}          &       & \multirow{3}[2]{*}{6} & 50    & 2.1   & 30.9  & 0.0(5) & 0.0 & 0.2 & \textbf{2.1} & \textbf{5.3} & \textbf{0.0(5)} & \textbf{0.0} & \textbf{0.2} & 2.1   & 7.7   & 0.0(5) & 0.0  \\
          &       &       & 60    & 2.3   & 102.7 & 0.0(5) & 0.0 & 0.4 & \textbf{2.2} & \textbf{44.7} & \textbf{0.0(5)} & \textbf{0.0} & \textbf{0.4} & 2.2   & 48.6  & 0.0(5) & 0.0  \\
          &       &       & 70    & \textbf{4.0} & \textbf{281.8} & \textbf{0.0(5)} & \textbf{0.0} & \textbf{0.4} & 4.1   & 1484.1 & 0.4(1) & 0.1 & 0.4 & 4.1   & 1451.1 & 0.2(1) & 0.0  \\
\cline{3-18}          &       & \multirow{3}[2]{*}{12} & 50    & 2.0   & 3.6   & 0.0(5) & 0.0 & 0.6 & \textbf{2.0} & \textbf{2.2} & \textbf{0.0(5)} & \textbf{0.0} & \textbf{0.7} & \textbf{2.0} & \textbf{2.2} & \textbf{0.0(5)} & \textbf{0.0}  \\
          &       &       & 60    & 2.0   & 24.3  & 0.0(5) & 0.0 & 1.5 & 2.0   & 4.4   & 0.0(5) & 0.0 & 1.6 & \textbf{2.0} & \textbf{4.3} & \textbf{0.0(5)} & \textbf{0.0}  \\
          &       &       & 70    & 2.0   & 52.5  & 0.0(5) & 0.0 & 0.7 & 2.0   & 10.3  & 0.0(5) & 0.0 & 0.7 & \textbf{2.0} & \textbf{10.0} & \textbf{0.0(5)} & \textbf{0.0}  \\
 
      \hline
    \hline
    \end{tabular}%
    }
\caption{Performance of formulations \name, \named, and \named~+~VI on graph100 and graph120}
  \label{tab:g100120}%
%
  \centering
   \resizebox{\hsize}{!}{
    \begin{tabular}{|l|rrrr|rrrr|rrrr|}
\cline{2-13}    \multicolumn{1}{r|}{} & \multicolumn{4}{c|}{\name}     & \multicolumn{4}{c|}{\named}      & \multicolumn{4}{c|}{\named~+~VI} \\
\cline{2-13}    \multicolumn{1}{r|}{} & \multicolumn{1}{r}{$c$} & \multicolumn{1}{r}{$v$} & \multicolumn{1}{r}{$bv$} & \multicolumn{1}{r|}{nodes} & \multicolumn{1}{r}{$c$} & \multicolumn{1}{r}{$v$} & \multicolumn{1}{r}{$bv$} & \multicolumn{1}{r|}{nodes} & \multicolumn{1}{r}{$c$} & \multicolumn{1}{r}{$v$} & \multicolumn{1}{r}{$bv$} & \multicolumn{1}{r|}{nodes} \\
    \hline
    graph100 & 796249.9 & 632308.2 & 316639.7 & 575.3 & 55293.7 & 4581.3 & 3605.5 & 95947.3 & 55310.7 & 4581.3 & 3605.5 & 71610.8 \\
    graph120 & 1556268.4 & 1236714.5 & 619029.0 & 1776.1 & 88816.1 & 6384.0 & 5035.1 & 166633.4 & 88847.8 & 6384.0 & 5035.1 & 154253.2 \\
    \hline
    \end{tabular}%
    }
  \caption{Number of constraints, variables, and nodes visited in the branching tree of formulations \name, \named, and \named~+~VI on graph100 and graph120}
  \label{tab:nodesgraph100120}%
\end{table}%

Analysing the results shown in Table~\ref{tab:g100120}, we can conclude that the difficulty of solving the instances is highly dependent on the parameters ($B,R$ and $p$). It can be observed that the instances become more difficult as $B$ and $p$ decrease and $R$ increases. As shown in Table~\ref{tab:g100120} and Figure~\ref{performanceProfilesbig}, the performance of these formulations is quite similar. An interesting observation is that they complement each other. In other words, there are instances in which formulation \name did not find the optimal solution within the time limit but \named~+~VI did and vice versa. Nevertheless, \named~+~VI found the optimal solution in more instances than \named and all the instances solved to optimality by \named were also solved to optimality by \named~+~VI. It can be appreciated that the average preprocessing time is lower than fourteen seconds in all cases. Furthermore, the linear relaxation of the formulations ($G_{LP}^t$\%) is quite good, being on average 2.1\% for \name and 2.3\% for \named and \named~+~VI. 

In view of the results reported in this subsection, we conclude that the best formulations for solving Up-MCLP on complete graphs are \name and \named~+~VI. In the next subsection, sparser graphs will be analysed.  

\subsection{Results on sparse graphs}

The aim of this subsection is to compare the proposed formulations on sparse graphs. In particular, we used the uncapacited $p$-median  datasets from the OR-Library. For interested readers, non-parametric tests (Friedman  test  and  Post-Hoc Holland  adjust) to assess the statistical significance of the comparison among the different formulations can be found in supplementary material. They show that there are significant differences between the formulations presented. 
\afterpage{
\begin{landscape}
\begin{table}[p]
  \centering
 	\resizebox{\linewidth}{!}{
    \begin{tabular}{|c|c|c|r|rrrrr|rrrr|rrr|rrrrr|rrrr|}
    \hline
    \multirow{2}[4]{*}{Data} & \multirow{2}[4]{*}{$B\%$} & \multirow{2}[4]{*}{$p$} & \multicolumn{1}{c|}{\multirow{2}[4]{*}{$R\%$}} & \multicolumn{5}{c|}{\name}     & \multicolumn{4}{c|}{\namet}       & \multicolumn{3}{c|}{\namet~+~VI}    & \multicolumn{5}{c|}{\named} &\multicolumn{4}{c|}{\named~+~VI} \\
\cline{5-25}          &       &       &       & \multicolumn{1}{r}{$t_{st}$} & \multicolumn{1}{r}{$t_{total}$} & \multicolumn{1}{r}{$G\%$} & \multicolumn{1}{r}{$G_{BS}^t\%$}& \multicolumn{1}{r|}{$G_{LP}^t$\%} &   \multicolumn{1}{r}{$t_{total}$} & \multicolumn{1}{r}{$G\%$} & \multicolumn{1}{r}{$G_{BS}^t\%$}& \multicolumn{1}{r|}{$G_{LP}^t$\%} &  \multicolumn{1}{r}{$t_{total}$} & \multicolumn{1}{r}{$G\%$} &  \multicolumn{1}{r|}{$G_{BS}^t\%$} & \multicolumn{1}{r}{$t_{st}$} & \multicolumn{1}{r}{$t_{total}$} & \multicolumn{1}{r}{$G\%$} & \multicolumn{1}{r}{$G_{BS}^t\%$}& \multicolumn{1}{r|}{$G_{LP}^t$\%} &\multicolumn{1}{r}{$t_{st}$} & \multicolumn{1}{r}{$t_{total}$} & \multicolumn{1}{r}{$G\%$} & \multicolumn{1}{r|}{$G_{BS}^t\%$}\\
    \hline
    \hline
    \multicolumn{1}{|c|}{\multirow{27}[18]{*}{\begin{turn}{90} pmeds, $|V|=100, |E|=195.2$\end{turn}}} & \multirow{9}[5]{*}{0.5} & \multirow{3}[1]{*}{1} & 50    & 1.4   & 26.1  & 0.0(5) & 0.0 & 4.3 & 918.8 & 17.5(4) & 2.3 & 87.5 & 474.8 & 23.1(4) & 4.5 & \textbf{1.4} & \textbf{20.5} & \textbf{0.0(5)} & \textbf{0.0} & \textbf{4.3} & 1.4   & 24.2  & 0.0(5) & 0.0 \\
          &       &       & 60    & \textbf{1.5} & \textbf{15.5} & \textbf{0.0(5)} & \textbf{0.0} & \textbf{2.1} & 1764.7 & 39.5(1) & 8.7 & 55.3 & 1168.2 & 31.9(3) & 7.8 & 1.5   & 38.7  & 0.0(5) & 0.0 & 2.1 & 1.5   & 25.7  & 0.0(5) & 0.0 \\
          &       &       & 70    & \textbf{1.5} & \textbf{32.6} & \textbf{0.0(5)} & \textbf{0.0} & \textbf{4.0} & 1801.9 & 44.0(0) & 8.9 & 35.8 & 1580.2 & 75.0(1) & 16.9 & 1.5   & 991.7 & 1.5(3) & 0.0 & 4.1 & 1.5   & 572.9 & 0.3(4) & 0.0 \\
\cline{3-25}          &       & \multirow{3}[2]{*}{5} & 50    & \textbf{1.2} & \textbf{1.4} & \textbf{0.0(5)} & \textbf{0.0} & \textbf{1.1} & 14.5  & 0.0(5) & 0.0 & 75.8 & 13.5  & 0.0(5) & 0.0 & \textbf{1.2} & \textbf{1.4} & \textbf{0.0(5)} & \textbf{0.0} & \textbf{2.7} & \textbf{1.2} & \textbf{1.4} & \textbf{0.0(5)} & \textbf{0.0} \\
          &       &       & 60    & \textbf{1.2} & \textbf{1.8} & \textbf{0.0(5)} & \textbf{0.0} & \textbf{1.3} & 575.0 & 0.0(5) & 0.0 & 57.9 & 589.1 & 0.0(5) & 0.0 & 1.2   & 2.1   & 0.0(5) & 0.0 & 3.4 & 1.2   & 2.3   & 0.0(5) & 0.0 \\
          &       &       & 70    & 1.2   & 3.9   & 0.0(5) & 0.0 & 1.7 & 803.4 & 3.2(3) & 0.0 & 38.8 & 979.5 & 2.5(3) & 0.0 & 1.2   & 3.9   & 0.0(5) & 0.0 & 2.6 & \textbf{1.2} & \textbf{3.7} & \textbf{0.0(5)} & \textbf{0.0} \\
\cline{3-25}          &       & \multirow{3}[2]{*}{10} & 50    & 1.2   & 1.2   & 0.0(5) & 0.0 & 0.1 & 0.8   & 0.0(5) & 0.0 & 42.8 & \textbf{0.7} & \textbf{0.0(5)} & \textbf{0.0} & 1.1   & 1.2   & 0.0(5) & 0.0 & 3.7 & 1.1   & 1.2   & 0.0(5) & 0.0 \\
          &       &       & 60    & \textbf{1.2} & \textbf{1.2} & \textbf{0.0(5)} & \textbf{0.0} & \textbf{0.7} & 1.6   & 0.0(5) & 0.0 & 46.9 & 1.8   & 0.0(5) & 0.0 & 1.1   & 1.3   & 0.0(5) & 0.0 & 2.9 & 1.1   & 1.3   & 0.0(5) & 0.0 \\
          &       &       & 70    & \textbf{1.2} & \textbf{1.4} & \textbf{0.0(5)} & \textbf{0.0} & \textbf{1.6} & 9.2   & 0.0(5) & 0.0 & 35.4 & 28.2  & 0.0(5) & 0.0 & 1.2   & 1.5   & 0.0(5) & 0.0 & 3.0 & 1.2   & 1.6   & 0.0(5) & 0.0 \\
\cline{2-25}          & \multirow{9}[6]{*}{1} & \multirow{3}[2]{*}{1} & 50    & 1.7   & 570.6 & 0.0(5) & 0.0 & 7.2 & 1496.5 & 31.6(1) & 2.7 & 77.1 & 1137.1 & 83.1(3) & 13.1 & 1.7   & 202.6 & 0.0(5) & 0.0 & 7.3 & \textbf{1.6} & \textbf{68.7} & \textbf{0.0(5)} & \textbf{0.0} \\
          &       &       & 60    & 1.7   & 1460.1 & 2.4(2) & 0.0 & 7.7 & 1800.6 & 40.7(0) & 4.9 & 49.6 & 1488.5 & 38.2(1) & 6.7 & 1.7   & 960.6 & 3.1(4) & 1.3 & 7.8 & \textbf{1.8} & \textbf{850.1} & \textbf{3.1(4)} & \textbf{1.3} \\
          &       &       & 70    & \textbf{1.8} & \textbf{425.1} & \textbf{0.0(5)} & \textbf{0.0} & \textbf{4.8} & 1801.2 & 29.3(0) & 5.7 & 30.1 & 1800.7 & 349.5(0) & 51.1 & 1.8   & 1425.3 & 2.7(2) & 0.0 & 4.9 & 1.8   & 1513.8 & 3.8(1) & 0.0 \\
\cline{3-25}          &       & \multirow{3}[2]{*}{5} & 50    & 1.2   & 1.8   & 0.0(5) & 0.0 & 2.8 & 84.1  & 0.0(5) & 0.0 & 75.7 & 77.7  & 0.0(5) & 0.0 & 1.2 & 1.8 & 0.0(5) & 0.0 & 4.9 & \textbf{1.2} & \textbf{1.7} & \textbf{0.0(5)} & \textbf{0.0} \\
          &       &       & 60    & \textbf{1.2} & \textbf{3.4} & \textbf{0.0(5)} & \textbf{0.0} & \textbf{5.0} & 918.9 & 0.7(4) & 0.1 & 55.8 & 1116.2 & 3.9(2) & 0.0 & 1.2   & 22.7  & 0.0(5) & 0.0 & 6.6 & 1.3   & 18.2  & 0.0(5) & 0.0 \\
          &       &       & 70    & 1.3   & 62.1  & 0.0(5) & 0.0 & 5.2 & 1718.3 & 8.3(1) & 0.2 & 36.3 & 1602.7 & 9.6(1) & 0.3 & 1.3   & 76.6  & 0.0(5) & 0.0 & 5.3 & \textbf{1.3} & \textbf{39.9} & \textbf{0.0(5)} & \textbf{0.0} \\
\cline{3-25}          &       & \multirow{3}[2]{*}{10} & 50    & 1.1   & 1.2   & 0.0(5) & 0.0 & 1.5 & 0.8   & 0.0(5) & 0.0 & 43.1 & \textbf{0.7} & \textbf{0.0(5)} & \textbf{0.0} & 1.1   & 1.2   & 0.0(5) & 0.0 & 4.4 & 1.1   & 1.2   & 0.0(5) & 0.0 \\
          &       &       & 60    & \textbf{1.2} & \textbf{1.3} & \textbf{0.0(5)} & \textbf{0.0} & \textbf{1.7} & 2.6   & 0.0(5) & 0.0 & 46.3 & 2.4   & 0.0(5) & 0.0 & 1.2   & 1.9   & 0.0(5) & 0.0 & 4.2 & 1.2   & 1.7   & 0.0(5) & 0.0 \\
          &       &       & 70    & \textbf{1.2} & \textbf{1.7} & \textbf{0.0(5)} & \textbf{0.0} & \textbf{3.2} & 18.0  & 0.0(5) & 0.0 & 34.6 & 25.3  & 0.0(5) & 0.0 & 1.2   & 2.6   & 0.0(5) & 0.0 & 5.6 & 1.2   & 2.5   & 0.0(5) & 0.0 \\
\cline{2-25}          & \multirow{9}[6]{*}{5} & \multirow{3}[2]{*}{1} & 50    & 1.5   & 1015.6 & 0.0(5) & 0.0 & 3.0 & 1288.4 & 24.(2) & 2.7 & 55.1 & 1155.6 & 25.7(2) & 3.2 & \textbf{1.5} & \textbf{64.8} & \textbf{0.0(5)} & \textbf{0.0} & \textbf{3.0} & 1.5   & 85.5  & 0.0(5) & 0.0 \\
          &       &       & 60    & 1.6   & 1407.8 & 3.4(2) & 0.5 & 4.5 & 1800.6 & 25.1(0) & 4.8 & 33.7 & 1690.5 & 27.4(1) & 6.2 & \textbf{1.6} & \textbf{869.4} & \textbf{2.5(3)} & \textbf{0.1} & \textbf{4.6} & 1.6   & 1140.9 & 2.5(2) & 0.1 \\
          &       &       & 70    & 1.6   & 1205.0 & 1.8(2) & 0.0 & 2.7 & 1801.7 & 22.6(0) & 5.7 & 19.3 & 1758.5 & 268.6(1) & 35.2 & 1.6   & 1077.4 & 1.8(3) & 0.0 & 2.8 & \textbf{1.7} & \textbf{946.9} & \textbf{1.9(3)} & \textbf{0.0} \\
\cline{3-25}          &       & \multirow{3}[2]{*}{5} & 50    & \textbf{1.2} & \textbf{2.3} & \textbf{0.0(5)} & \textbf{0.0} & \textbf{4.1} & 102.8 & 0.0(5) & 0.0 & 66.8 & 132.1 & 0.0(5) & 0.0 & 1.2   & 3.0   & 0.0(5) & 0.0 & 4.2 & 1.2   & 2.9   & 0.0(5) & 0.0 \\
          &       &       & 60    & \textbf{1.2} & \textbf{5.2} & \textbf{0.0(5)} & \textbf{0.0} & \textbf{2.7} & 686.2 & 0.3(4) & 0.0 & 46.2 & 738.7 & 0.0(5) & 0.0 & 1.2   & 21.0  & 0.0(5) & 0.0 & 2.7 & 1.2   & 17.6  & 0.0(5) & 0.0 \\
          &       &       & 70    & 1.3   & 207.3 & 0.0(5) & 0.0 & 3.0 & 1418.4 & 2.8(3) & 0.1 & 28.8 & 1458.9 & 4.2(2) & 0.0 & 1.3   & 48.5  & 0.0(5) & 0.0 & 3.0 & \textbf{1.3} & \textbf{48.1} & \textbf{0.0(5)} & \textbf{0.0} \\
\cline{3-25}          &       & \multirow{3}[2]{*}{10} & 50    & 1.2   & 1.2   & 0.0(5) & 0.0 & 0.8 & 1.0   & 0.0(5) & 0.0 & 38.4 & \textbf{0.8} & \textbf{0.0(5)} & \textbf{0.0} & 1.1   & 1.3   & 0.0(5) & 0.0 & 1.1 & 1.2   & 1.3   & 0.0(5) & 0.0 \\
          &       &       & 60    & \textbf{1.2} & \textbf{1.5} & \textbf{0.0(5)} & \textbf{0.0} & \textbf{3.1} & 4.3   & 0.0(5) & 0.0 & 40.7 & 4.6   & 0.0(5) & 0.0 & 1.2   & 1.9   & 0.0(5) & 0.0 & 4.0 & 1.2   & 1.8   & 0.0(5) & 0.0 \\
          &       &       & 70    & \textbf{1.2} & \textbf{2.2} & \textbf{0.0(5)} & \textbf{0.0} & \textbf{2.9} & 61.8  & 0.0(5) & 0.0 & 29.4 & 38.3  & 0.0(5) & 0.0 & 1.2   & 4.1   & 0.0(5) & 0.0 & 3.0 & 1.2   & 3.9   & 0.0(5) & 0.0 \\
   \hline
   \hline
    \end{tabular}%
    }
\caption{Performance of formulations \name, \namet, \namet~+~VI, \named, and \named~+~VI on pmeds}
  \label{tab:pmed1}%
\vspace{1cm}
  \resizebox{\linewidth}{!}{
    \begin{tabular}{|l|rrrr|rrrr|rrrr|rrrr|rrrr|}
\cline{2-21}    \multicolumn{1}{r|}{} & \multicolumn{4}{c|}{\name}     & \multicolumn{4}{c|}{\namet}       & \multicolumn{4}{c|}{\namet~+~VI}    & \multicolumn{4}{c|}{\named}      & \multicolumn{4}{c|}{\named~+~VI} \\
\cline{1-21}    \multicolumn{1}{|r|}{Data} & \multicolumn{1}{r}{$c$} & \multicolumn{1}{r}{$v$} & \multicolumn{1}{r}{$bv$} & \multicolumn{1}{r|}{nodes} & \multicolumn{1}{r}{$c$} & \multicolumn{1}{r}{$v$} & \multicolumn{1}{r}{$bv$} & \multicolumn{1}{r|}{nodes} & \multicolumn{1}{r}{$c$} & \multicolumn{1}{r}{$v$} & \multicolumn{1}{r}{$bv$} & \multicolumn{1}{r|}{nodes} & \multicolumn{1}{r}{$c$} & \multicolumn{1}{r}{$v$} & \multicolumn{1}{r}{$bv$} & \multicolumn{1}{r|}{nodes} & \multicolumn{1}{r}{$c$} & \multicolumn{1}{r}{$v$} & \multicolumn{1}{r}{$bv$} & \multicolumn{1}{r|}{nodes} \\
\cline{1-21}    pmeds & 166204.5 & 120362.9 & 60983.3 & 412.4 & 1023.2 & 606.4 & 376.0 & 1296509.9 & 2333.6 & 606.4 & 376.0 & 786929.8 & 7368.9 & 2248.1 & 2017.7 & 29238.2 & 7396.4 & 2248.1 & 2017.7 & 25113.6 \\
\hline
    \end{tabular}%
    }
  \caption{Number of constraints, variables, and nodes visited in the branching tree of formulations \name, \namet, \namet~+~VI, \named, and \named~+~VI on pmeds}
  \label{tab:nodespmeds}%
\end{table}%
\end{landscape}
}

Table~\ref{tab:pmed1} has the same structure as Table~\ref{tab:g3040}. In this table, we reported the results of formulations \name, \namet, \namet~+~VI, \named, and \named~+~VI on the smallest pmed datasets (pmed1-pmed5), named pmeds. These networks contain 100 nodes and 195.2 edges on average (the smallest have 190 and the largest 198). The provided results are the average over the five datasets where the other parameters were randomly generated following the procedure described in Subsection~\ref{sub:data}. Moreover, a similar analysis to Table~\ref{tab:nodesgraph100120} is reported in Table~\ref{tab:nodespmeds} for pmeds instances. A detailed table can be found in the supplementary material.

Similarly to complete graphs, it can be seen in Table~\ref{tab:pmed1} that the difficulty of the instances is highly parameter-dependent. In addition, it is shown that the preprocessing time is small (less than two seconds in all instances). Furthermore, the number of instances solved to optimality by \name, \named,  and \named~+~VI is considerably higher than the ones solved by \namet and \namet~+~VI. Observe that the average of the linear relaxation gaps of \name (3.1\%), \named and \named~+~VI (4.0\%) are better than the linear relaxation gap of \namet and \namet~+~VI (47.5\%). In Table~\ref{tab:nodespmeds} the different sizes of the problem for each formulation can be appreciated.

In Table~\ref{tab:pmed2}, we report the result obtained on datasets of larger size, named pmedb. More concretely, the table shows the results of formulations \name, \named, and \named~+~VI on larger pmed datasets (pmed6-pmed10). These networks contain 200 nodes and  777.8 edges on average (the smallest have 774 and the largest 785). Note that the results of \namet, \namet~+~VI are not reported because very few instances were solved to optimality. In Table~\ref{tab:nodespmedb}, a similar analysis to Table~\ref{tab:nodespmeds} is depicted. A detailed table can be found in the supplementary material. In Figure~\ref{performanceProfilespmed}, the performance profile of the number of solved instances using these formulations is depicted.
\afterpage{
\begin{table}[htbp]
  \centering
 	\resizebox{\linewidth}{!}{
    \begin{tabular}{|c|c|c|r|rrrrr|rrrrr|rrrr|}
    \hline
    \multirow{2}[4]{*}{Data} & \multirow{2}[4]{*}{$B\%$} & \multirow{2}[4]{*}{$p$} & \multicolumn{1}{c|}{\multirow{2}[4]{*}{$R\%$}} & \multicolumn{5}{c|}{\name}     & \multicolumn{5}{c|}{\named} &\multicolumn{4}{c|}{\named~+~VI} \\
\cline{5-18}          &       &       &       & \multicolumn{1}{r}{$t_{st}$} & \multicolumn{1}{r}{$t_{total}$} & \multicolumn{1}{r}{$G\%$} & \multicolumn{1}{r}{$G_{BS}^t\%$}& \multicolumn{1}{r|}{$G_{LP}^t$\%} &    \multicolumn{1}{r}{$t_{st}$} & \multicolumn{1}{r}{$t_{total}$} & \multicolumn{1}{r}{$G\%$} & \multicolumn{1}{r}{$G_{BS}^t\%$}& \multicolumn{1}{r|}{$G_{LP}^t$\%} &\multicolumn{1}{r}{$t_{st}$} & \multicolumn{1}{r}{$t_{total}$} & \multicolumn{1}{r}{$G\%$} & \multicolumn{1}{r|}{$G_{BS}^t\%$}\\
    \hline
    \hline
    \multicolumn{1}{|c|}{\multirow{27}[18]{*}{\begin{turn}{90} pmedb, $|V|=200, |E|=777.8$\end{turn}}}& \multirow{9}[5]{*}{0.5} & \multirow{3}[1]{*}{1} & 50    & 12.5  & 1815.3 & 5881.2(0) & 18.7 & 9.1 & 12.4  & 1812.9 & 9.3(0) & 0.1 & 9.9 & 12.3  & 1812.7 & 8.9(0) & 0.0  \\
          &       &       & 60    & 18.3  & 1821.2 & 6383.5(0) & 33.1 & 7.4 & 18.5  & 1819.1 & 7.5(0) & 0.0 & 8.0 & 18.4  & 1818.7 & 7.7(0) & 0.1  \\
          &       &       & 70    & 36.3  & 1839.9 & 13701.3(0) & 24.8 & 5.8 & 37.6  & 1848.8 & 5.8(0) & 0.0 & 6.4 & 37.1  & 1837.5 & 6.1(0) & 0.3 \\
\cline{3-18}          &       & \multirow{3}[2]{*}{10} & 50    & 9.3   & 22.4  & 0.0(5) & 0.0 & 2.1 & \textbf{9.2} & \textbf{14.8} & \textbf{0.0(5)} & \textbf{0.0} & \textbf{4.3} & 9.1   & 15.3  & 0.0(5) & 0.0 \\
          &       &       & 60    & 9.5   & 1451.5 & 1.2(2) & 0.0 & 4.4 & 9.3   & 1150.8 & 1.1(2) & 0.2 & 6.5 & \textbf{9.3} & \textbf{975.3} & \textbf{1.3(3)} & \textbf{0.3}  \\
          &       &       & 70    & 10.0  & 1811.3 & 2.8(0) & 0.6 & 3.9 & 9.8   & 1815.2 & 2.7(0) & 1.0 & 4.4 & 9.8   & 1811.3 & 2.6(0) & 1.0 \\
\cline{3-18}          &       & \multirow{3}[2]{*}{20} & 50    & \textbf{8.9} & \textbf{9.1} & \textbf{0.0(5)} & \textbf{0.0} & \textbf{1.8} & 8.9   & 9.2   & 0.0(5) & 0.0 & 3.6 & \textbf{8.8} & \textbf{9.1} & \textbf{0.0(5)} & \textbf{0.0} \\
          &       &       & 60    & 9.0   & 45.4  & 0.0(5) & 0.0 & 2.5 & \textbf{9.0} & \textbf{13.3} & \textbf{0.0(5)} & \textbf{0.0} & \textbf{3.8} & 8.9   & 14.4  & 0.0(5) & 0.0 \\
          &       &       & 70    & 9.4   & 785.8 & 0.1(4) & 0.1 & 2.7 & \textbf{9.2} & \textbf{617.7} & \textbf{0.4(4)} & \textbf{0.0} & \textbf{5.2} & 9.2   & 771.0 & 0.3(3) & 0.0  \\
\cline{2-18}          & \multirow{9}[6]{*}{1} & \multirow{3}[2]{*}{1} & 50    & 11.2  & 1814.1 & 6101.7(0) & 23.3 & 8.2 & 11.2  & 1811.5 & 8.4(0) & 0.0 & 8.4 & 11.1  & 1811.5 & 8.4(0) & 0.0 \\
          &       &       & 60    & 16.2  & 1819.9 & 13905.8(0) & 25.2 & 7.4 & 16.3  & 1819.4 & 7.5(0) & 0.0 & 7.6 & 16.1  & 1816.8 & 7.5(0) & 0.0  \\
          &       &       & 70    & 33.7  & 1836.7 & 17155.8(0) & 26.3 & 5.1 & 33.7  & 1834.1 & 5.4(0) & 0.2 & 5.3 & 33.5  & 1834.1 & 5.5(0) & 0.4  \\
\cline{3-18}          &       & \multirow{3}[2]{*}{10} & 50    & 9.3   & 299.3 & 0.0(5) & 0.0 & 3.7 & 9.0   & 57.5  & 0.0(5) & 0.0 & 4.7 & \textbf{9.0} & \textbf{53.7} & \textbf{0.0(5)} & \textbf{0.0} \\
          &       &       & 60    & 9.4   & 1810.0 & 3.2(0) & 0.1 & 5.5 & 9.5   & 1566.1 & 2.2(1) & 0.5 & 5.8 & \textbf{9.3} & \textbf{1526.1} & \textbf{1.7(1)} & \textbf{0.2}  \\
          &       &       & 70    & 9.5   & 1811.1 & 3.9(0) & 0.7 & 3.8 & 9.6   & 1811.5 & 2.7(0) & 0.4 & 3.8 & 9.5   & 1810.6 & 4.1(0) & 1.9 \\
\cline{3-18}          &       & \multirow{3}[2]{*}{20} & 50    & 9.2   & 9.5   & 0.0(5) & 0.0 & 3.2 & 9.1   & 9.4   & 0.0(5) & 0.0 & 4.2 & \textbf{8.9} & \textbf{9.3} & \textbf{0.0(5)} & \textbf{0.0}  \\
          &       &       & 60    & 9.0   & 141.0 & 0.0(5) & 0.0 & 3.6 & 9.3   & 23.2  & 0.0(5) & 0.0 & 4.5 & \textbf{9.0} & \textbf{20.7} & \textbf{0.0(5)} & \textbf{0.0} \\
          &       &       & 70    & 9.3   & 1613.6 & 1.0(1) & 0.1 & 3.2 & 9.2   & 1158.1 & 0.5(3) & 0.0 & 4.1 & \textbf{9.1} & \textbf{959.8} & \textbf{0.3(4)} & \textbf{0.0}  \\
\cline{2-18}          & \multirow{9}[6]{*}{5} & \multirow{3}[2]{*}{1} & 50    & 11.2  & 1546.2 & 6102.5(2) & 25.3 & 0.0 & \textbf{11.2} & \textbf{181.1} & \textbf{0.0(5)} & \textbf{0.0} & \textbf{0.0} & 11.0  & 185.9 & 0.0(5) & 0.0  \\
          &       &       & 60    & 16.6  & 1742.8 & 6551.2(1) & 24.6 & 0.1 & 16.0  & 917.2 & 0.1(3) & 0.0 & 0.1 & \textbf{16.0} & \textbf{803.1} & \textbf{0.1(4)} & \textbf{0.0}  \\
          &       &       & 70    & 33.2  & 1837.5 & 11417.1(0) & 29.4 & 0.2 & 32.6  & 853.6 & 0.2(4) & 0.0 & 0.2 & \textbf{32.6} & \textbf{769.9} & \textbf{0.2(4)} & \textbf{0.0}  \\
\cline{3-18}          &       & \multirow{3}[2]{*}{10} & 50    & 9.1   & 342.8 & 0.0(5) & 0.0 & 1.7 & 9.2   & 34.8  & 0.0(5) & 0.0 & 1.7 & \textbf{9.0} & \textbf{30.0} & \textbf{0.0(5)} & \textbf{0.0} \\
          &       &       & 60    & 9.6   & 1527.2 & 1.1(1) & 0.2 & 1.7 & 9.2   & 696.4 & 0.0(5) & 0.0 & 1.7 & \textbf{9.2} & \textbf{434.5} & \textbf{0.0(5)} & \textbf{0.0} \\
          &       &       & 70    & 9.5   & 1811.4 & 1.2(0) & 0.8 & 0.5 & 9.5   & 733.6 & 0.2(4) & 0.0 & 0.5 & \textbf{9.4} & \textbf{687.4} & \textbf{0.5(4)} & \textbf{0.3}  \\
\cline{3-18}          &       & \multirow{3}[2]{*}{20} & 50    & \textbf{9.0} & \textbf{9.4} & \textbf{0.0(5)} & \textbf{0.0} & \textbf{1.2} & 8.9   & 9.6   & 0.0(5) & 0.0 & 1.3 & 8.8   & 9.5   & 0.0(5) & 0.0 \\
          &       &       & 60    & 9.2   & 44.5  & 0.0(5) & 0.0 & 1.3 & 9.0   & 17.1  & 0.0(5) & 0.0 & 1.4 & \textbf{8.9} & \textbf{15.5} & \textbf{0.0(5)} & \textbf{0.0}  \\
          &       &       & 70    & 9.3   & 1324.6 & 0.4(2) & 0.0 & 0.8 & 9.1   & 730.6 & 0.0(4) & 0.0 & 0.8 & \textbf{9.3} & \textbf{493.2} & \textbf{0.1(4)} & \textbf{0.0}  \\
    \hline
    \hline
    \end{tabular}%
    }
\caption{Performance of formulations \name, \named, and \named~+~VI on pmedb}
  \label{tab:pmed2}%
  \vspace{1cm}
	\resizebox{\linewidth}{!}{
    \begin{tabular}{|l|rrrr|rrrr|rrrr|}
\cline{2-13}    \multicolumn{1}{r|}{} & \multicolumn{4}{c|}{\name}     & \multicolumn{4}{c|}{\named}      & \multicolumn{4}{c|}{\named~+~VI} \\
    \hline
    Data  & \multicolumn{1}{r}{$c$} & \multicolumn{1}{r}{$v$} & \multicolumn{1}{r}{$bv$} & \multicolumn{1}{r|}{nodes} & \multicolumn{1}{r}{$c$} & \multicolumn{1}{r}{$v$} & \multicolumn{1}{r}{$bv$} & \multicolumn{1}{r|}{nodes} & \multicolumn{1}{r}{$c$} & \multicolumn{1}{r}{$v$} & \multicolumn{1}{r}{$bv$} & \multicolumn{1}{r|}{nodes} \\
    \hline
    pmedb & 2147283.6 & 1599445.8 & 802394.7 & 1809.3 & 35199.8 & 6849.2 & 6350.2 & 97218.6 & 35336.9 & 6849.2& 6350.2 & 80667.2 \\
    \hline
    \end{tabular}%
    }
      \caption{Number of constraints, variables, and nodes visited in the branching tree of formulations \name, \named, and \named~+~VI on pmedb}
  \label{tab:nodespmedb}%

\end{table}%
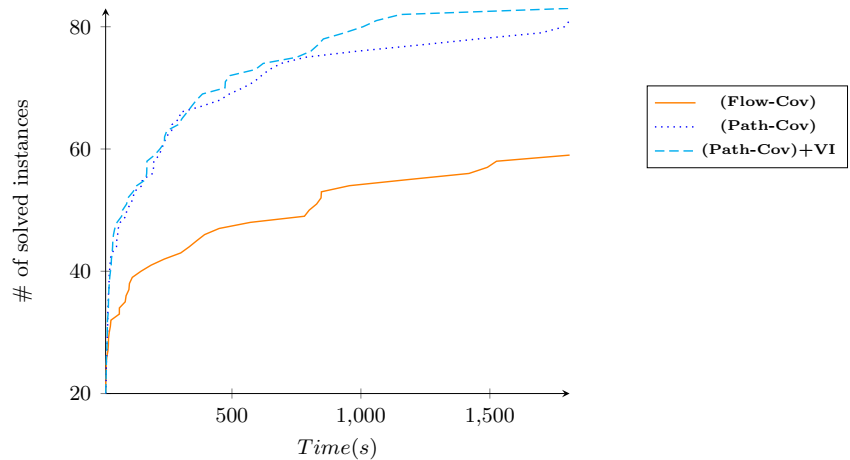
\begin{figure}[hbtp] \centering
	\begin{tikzpicture}[scale=0.9,font=\footnotesize]
		\begin{axis}[axis x line=bottom,  axis y line=left,
			xlabel=$Time(s)$,
			ylabel=\# of solved instances,
			legend style={at={(1.6,0.8)}}]
			\addplot[orange,semithick] plot coordinates {		
(10.197,20)
(10.262,21)
(10.625,22)
(11.739,23)
(12.216,24)
(12.234,25)
(14.84,26)
(19.109,27)
(20.791,28)
(21.303,29)
(24.124,30)
(28.544,31)
(30.739,32)
(63.368,33)
(63.496,34)
(86.081,35)
(89.534,36)
(100.496,37)
(102.452,38)
(113.503,39)
(147.101,40)
(186.478,41)
(237.758,42)
(301.463,43)
(334.271,44)
(363.515,45)
(394.008,46)
(450.256,47)
(571.684,48)
(781.014,49)
(799.59,50)
(828.305,51)
(845.439,52)
(846.855,53)
(955.62,54)
(1183.734,55)
(1419.093,56)
(1490.352,57)
(1525.848,58)
(1809.253,59)
			};
			\addlegendentry{\tiny\textbf{\name}}	
						\addplot[blue,semithick,dotted] plot coordinates {		
(10.3,20)
(10.943,21)
(11.72,22)
(12.282,23)
(12.652,24)
(13.011,25)
(14.266,26)
(14.754,27)
(15.326,28)
(16.072,29)
(17.835,30)
(17.969,31)
(18.076,32)
(18.459,33)
(19.418,34)
(20.597,35)
(20.738,36)
(23.033,37)
(23.257,38)
(23.831,39)
(24.957,40)
(27.79,41)
(27.884,42)
(34.497,43)
(52.221,44)
(53.24,45)
(56.738,46)
(59.099,47)
(66.922,48)
(89.467,49)
(91.859,50)
(106.605,51)
(117.508,52)
(119.742,53)
(150.99,54)
(152.574,55)
(193.032,56)
(196.095,57)
(196.89,58)
(217.932,59)
(226.842,60)
(236.686,61)
(245.32,62)
(265.93,63)
(267.707,64)
(296.894,65)
(299.324,66)
(382.974,67)
(453.377,68)
(490.054,69)
(542.525,70)
(586.407,71)
(614.484,72)
(645.313,73)
(690.624,74)
(783.911,75)
(985.558,76)
(1225.401,77)
(1453.738,78)
(1698.609,79)
(1789.298,80)
(1809.253,81)

			};
			\addlegendentry{\tiny\textbf{\named}}	
	\addplot[cyan,semithick,dash pattern=on 4pt off 2pt] plot coordinates {		
(10.37,20)
(11.06,21)
(11.746,22)
(12.095,23)
(12.986,24)
(12.999,25)
(13.023,26)
(14.062,27)
(15.704,28)
(15.785,29)
(16.626,30)
(16.722,31)
(17.234,32)
(19.698,33)
(20.836,34)
(21.161,35)
(21.838,36)
(22.765,37)
(24.259,38)
(26.114,39)
(28.929,40)
(31.727,41)
(32.62,42)
(35.552,43)
(37.149,44)
(38.049,45)
(40.685,46)
(45.697,47)
(54.45,48)
(72.26,49)
(75.87,50)
(91.474,51)
(96.53,52)
(116.372,53)
(130.678,54)
(165.772,55)
(169.318,56)
(170.071,57)
(170.343,58)
(199.758,59)
(220.563,60)
(237.923,61)
(238.807,62)
(247.638,63)
(292.395,64)
(306.846,65)
(325.742,66)
(343.691,67)
(362.338,68)
(385.677,69)
(472.973,70)
(474.391,71)
(491.625,72)
(590.605,73)
(621.227,74)
(750.7,75)
(805.74,76)
(832.598,77)
(854.449,78)
(935.942,79)
(1008.203,80)
(1059.905,81)
(1153.567,82)
(1809.253,83)
			};
			\addlegendentry{\tiny\textbf{\named \hspace{-0.1cm}+VI}}	

		\end{axis}
	\end{tikzpicture}
	\caption{Performance profile graph of \#solved instances using \name, \named, and \named~+~VI  formulations on pmedb dataset}  \label{performanceProfilespmed}
\end{figure}
}

As can be appreciated in Table~\ref{tab:pmed2} and Figure~\ref{performanceProfilespmed}, \named~+~VI outperforms \name and \named. Observe that although the gap of the linear relaxation of \name ($3.4\%$) is smaller than the gap of \named and \named~+~VI ($4\%$), the latter is the one that solves to optimality within the time limit the largest number of instances. 

Based on the results presented in this subsection, we conclude that the best formulation for solving Up-MCLP on pmedian graphs was \named~+~VI. Furthermore, the results included in this subsection show the usefulness of including the valid inequalities discussed in the paper. 

\section{Conclusions and outlook}
\label{sec:Conclusions}%
In this paper, we have tackled an interesting problem: the upgrading maximal covering location problem with edge length modifications, Up-MCLP. As far as we know, it is the first time that this problem is discussed in the literature. 

Since we were dealing with a new problem, we proposed three different mixed-integer formulations to model the situation from various perspectives. Moreover, we developed an effective preprocessing phase, which fixed many variables and reduced the size of the problem considerably. Then, for each formulation, we provided several sets of valid inequalities. These constraints allowed us to strengthen the formulations and to reduce the symmetries contained in the problem, shortening the time to solve the formulations. The performance of the three formulations and the improvement provided by the preprocessing phase and the valid inequalities can be appreciated in the computational results included in the paper. In these experiments, it can be seen that the most efficient formulations for solving Up-MCLP are \name and \named. In complete graphs, there is little difference in performance between these two formulations, while in sparse graphs \named performs better than \name. In both types of graphs, the addition of valid inequalities allows us to optimally solve a larger number of instances within the time limit.

We believe that this paper is an encouraging starting point that opens up many opportunities for further research. While our preprocessing techniques and valid inequalities allow us to solve significantly larger instances, and in shorter time, the size of the solvable instances still falls short of what you would encounter in practice. This leads to the necessity of developing heuristics that can obtain good solutions (even if they are not optimal) in shorter time than the exact algorithms. An additional line of future work is to further study and develop formulations for other upgrading versions of location problems. For example, interesting and similar problems could be obtained by modifying the covering criterion (e.g. gradual coverage, cooperative coverage, etc.), the location criterion (e.g. center problems, set covering problems, etc.), and the upgrading assumptions (e.g. non-linear upgrading cost or additional requirements on the sets of edges to be reduced, e.g., they have to form a connected set).


\vspace*{2ex}
{\footnotesize
\section*{Acknowledgements}

\textbf{Marta Baldomero-Naranjo} was partially supported by research project PID2020-114594GB-C22 (Agencia Estatal de Investigaci\'on, Spain and the European Regional Development's funds (ERDF)), research projects FEDER-UCA18-106895 and P18-FR-1422 (Regional Government of Andaluc\'ia, Spain and ERDF), the BritishSpanish Society Scholarship Programme 2018 (Telef\'onica and the BritishSpanish Society), and PhD grant UCA/REC01VI/2017 (Universidad de C\'adiz). 

\textbf{Alfredo Mar\'in} was partially supported by research project PID2019-110886RB-100 (Ministerio de Ciencia e Innovaci\'on, Spain). Part of this research were conducted while the author was on sabbatical at the University of C\'adiz.

\textbf{Antonio M. Rodr\'iguez-Ch\'ia} was partially supported by research project PID2020-114594GB-C22 (Agencia Estatal de Investigaci\'on, Spain and ERDF), research projects FEDER-UCA18-106895 and P18-FR-1422 (Regional Government of Andaluc\'ia, Spain and ERDF), and research project NetmeetData (Fundaci\'on BBVA).    

The authors would like
to thank the anonymous reviewers for their comments and
suggestions.

\bibliographystyle{abbrvnat}
\bibliography{./references}
}

\end{document}